\documentclass[10pt,a4paper,twoside]{article}

\usepackage{amsfonts, amssymb, amsmath, amsthm}
\usepackage{mathrsfs} 
\usepackage{latexsym}
\usepackage{enumerate}

\usepackage{url}
\usepackage{csquotes}
\usepackage{graphicx}
\usepackage{placeins}
\usepackage{epstopdf}
\usepackage{subfigure}

\usepackage[colorlinks=true, pdftex]{hyperref}
\hypersetup{citecolor=blue, linkcolor=blue}
\usepackage{mathabx} 

\newtheorem*{quest}{Question}
\newtheorem{thm}{Theorem}[section]
\newtheorem{lem}[thm]{Lemma}
\newtheorem{cor}[thm]{Corollary}

\theoremstyle{remark}
\newtheorem{remark}[thm]{Remark}

\newcommand{\Z}[1]{\mathbb{Z}/#1\mathbb{Z}}
\newcommand{\N}[1]{\{1\cdots{}#1\}}
\let\ve=\varepsilon
\let\vp=\varphi

\def\goes{\mathrel{\rightarrow}}
\DeclareMathOperator{\Id}{Id}
\DeclareMathOperator{\Pure}{\mathcal{R}}
\DeclareMathOperator{\Base}{\mathscr{E}}
\newcommand{\Transform}[1]{#1^\star}
\newcommand{\lowm}{1}
\newcommand{\upm}{2}
\newcommand{\PartialTransform}[1]{\mathscr{F}(#1)}
\def\Vscr{\mathscr{V}}

\def\Dcal{\mathcal{D}}

\def\Ocal{\mathcal{O}}
\def\1{1\!\!\!1}
\def\myS{\mathscr{S}}
\def\gL{\mathfrak{L}}
\def\ooWo{I_0( W)}
\def\LL{\Gamma^*}
\def\va{\goes}
\def\mode{\mathbin{\,\textrm{mod}^*}}
\def\mod{\mathbin{\,\textrm{mod}\,}}
\let\oldotimes=\otimes
\renewcommand{\otimes}{\mathop{\oldotimes}}

\title{A spectral resolution of the large sieve
}

\author{O. Ramar\'e}

\begin{document}
\maketitle

\begin{abstract}
  The quadratic form
  $V(\varphi,Q)=\sum_{q\sim Q}\sum_{a\mode q}|S(\varphi,a/q)|^2$ and
  its eigenvalues are well understood when $Q=o(\sqrt{N})$, while
  $V(\varphi,Q)$ is expected to behave like a Riemann sum when
  $N=o(Q)$. The behavior in the range $Q\in[\sqrt{N},100 N]$ is still
  mysterious. In the present work we present a full spectral analysis
  when $Q\ge N^{7/8}$ in terms of the eigenvalues of a one-parameter
  family of nuclear difference operators. We show in particular that
  (a smoothed version of) the quadratic form $V(\varphi,Q)$ may stay
  \emph{away} from $(6/\pi^2)Q\sum_n|\varphi_n|^2$ when $Q\asymp N$,
  though only on a vector space of positive but small dimension.
  An improved version of this paper, with the same title, will appear
(2024 or 2025) in the Bulletin of the French Mathematical Society.
\end{abstract}


\section{Introduction and results}

\subsubsection*{Main consequence}
We are interested in this paper in the quantity $\sum_{q\sim
  Q}\sum_{a\mode q}|S(\varphi,a/q)|^2$ where $(\varphi_n)_{n\le N}$ is any sequence of complex numbers
and $ S(\varphi,\alpha)=\sum_{n\le N}\varphi_n e(n\alpha)$.
It is this quantity that we analyze. Our main steps in this analysis are
Theorem~\ref{Precise}, Formula~\eqref{formula11bb-1} and Theorem~\ref{yoddle}.
One of the
main consequence of our work is the next theorem.
\begin{thm}
  \label{lowerbound}
  There exists $c>0$ such that for every~$N$ large enough and
  $Q\in[cN/\sqrt{\log N},20N]$, we have
  \begin{equation*}
    \sum_{\lowm<q/Q\le \upm}\sum_{a\mode q}
    |S(\varphi, a/q)|^2\ge Q^2 e^{-c N/Q}\sum_m|\varphi_m|^2.
  \end{equation*}
\end{thm}
This is to be compared with the lower bound given by W. Duke \&
H. Iwaniec in \cite{Duke-Iwaniec*92}.  Note that the summation therein
extends over all classes $a$ modulo~$q$ rather than over the reduced
classes, see the remark following \cite[Theorem 2.7]{Ramare*06} on
this issue.  In particular, the principal character is included
(i.e. $q=1$) with a definite influence.  J.-C. Schlage-Puchta
in~\cite{SchlagePuchta*09} gives, for some random sequences, a lower bound of a large sieve
quantity under the sole assumption that $Q^2/N$ goes to infinity.
Read also the papers of P. Erd\"os \&
A. Renyi \cite{Erdos-Renyi*68} and of D. Wolke \cite{Wolke*74}.  

The proof of Theorem~\ref{lowerbound} will unfold in four steps:
\begin{itemize}
\item By appealing to the $\delta$-symbol technique, we relate the
  above sum to a sum of similar kind but where the moduli $h$ are much
  smaller, namely $h\le H$ for some $H$ of size roughly $N/Q$.
\item We then interpret, for each $h$, the intervening quantity as a
  scalar product of some function $\Pure_{N,h}(\varphi)$ together with
  the value of a difference operator applied at this same vector.
\item After analyzing the one-parameter family of compact operators that intervene,
  we use their eigenvalues to derive a spectral decomposition of the
  large sieve quantity we are interested in.
\item When $W$ is non-negative and $N/Q$ is small
  enough, we prove that these eigenvalues are $<1$ by using the harmonic analysis
  uncertainty principle. Theorem~\ref{lowerbound} is a consequence of that.
\end{itemize}
\subsubsection*{Setting the horizon for a lower bound}
\begin{quest}
  Do we have
  $\displaystyle \sum_{\lowm<q/Q\le \upm}\sum_{a\mode q}
    |S(\varphi, a/q)|^2\gg N\sum_m|\varphi_m|^2$ when $Q\ge
    N^{1/2+\ve}$ for some positive $\ve$?
\end{quest}
When
$N=\sum_{q\le Q}\phi(q)$, we gave in \cite[Theorem 1.2]{Ramare*06} the
(rather weak) lower bound
$\|\varphi\|_2^2\exp(\frac{-1+o(1)}{2} N\log N)$ for the quantity
$\sum_{q\le Q}\sum_{a\mode q} |S(\varphi,
a/q)|^2$. Theorem~\ref{Vic} implies that the better
lower bound $Q^2\|\varphi\|_2^2$ holds true as soon as $\varphi$ oscillates
enough along small arithmetic progressions in intervals of length
about~$Q$. The main result of \cite[Theorem 2.4]{Conrey-Iwaniec-Soundararajan*12}
by B. Conrey, H. Iwaniec and K. Soundararajan implies a similar lower
bound for functions $\varphi$ that are the convolution product of an
oscillating factor supported on $[1,Q^{1-\ve}]$ and a rather general
sequence.
\subsubsection*{Some functional transforms of our weight function}
The $\delta$-symbol technique involves some functional
transforms of our weight function~$W$ that we better treat before
starting the
analysis proper.
Assumptions $W$ being as above, we define
$\Transform{W}$ in \eqref{defWstar}, but the following expression valid
for $z\in\mathbb{R}$ is better:
\begin{equation*}
  \label{explicitTransform-1}
  \Transform{W}(z)
  =
  -2\sum_{n\ge1}\frac{\phi(n)}{n}
  \int_0^{\infty}\cos(2\pi n y)W(z/y)dy/y.
\end{equation*}
By Lemma~\ref{boundWstarC}, the function $\Transform{W}$ is even, twice
differentiable outside $z=0$ where it vanishes, and is of bounded
variations over $[0,1]$ and decreases like $1/z^{\frac72-\ve}$ at infinity.  The expression for its Mellin transform,
valid when $\Re s\in[0,3/2)$ is simply
$\check{\Transform{W}}(s)=\check{W}(s)\zeta(1-s)/\zeta(1+s)$, see
Lemma~\ref{MTWstar}, where $\check{W}(s)$ is the Mellin transform of
$W$.  We finally mention the following expression for its Fourier
transform, valid for $u\neq0$ and obtained in
Lemma~\ref{ExplicitWhatstar}:
\begin{equation}
  \label{neat}
  \hat{W}^\star(u)=
    \frac{6}{\pi^2}\int_0^\infty W(t)dt
    -
    \frac{1}{|u|}\sum_{n\ge1}\frac{\phi(n)}{n}W(n/|u|).
\end{equation}
This Fourier transform satisfies $\hat
W^\star(u)=\frac{6}{\pi^2}\int_0^\infty W(t)dt$ when $|u|\le 1/\upm$ and
$|u\hat W^\star(u)|\ll \exp-c_0\sqrt{\log|u|}$ otherwise, for some
positive constant $c_0$, ensuring that $\hat W^\star(u)$ belongs to
$L^1$. It is worth specifying that $\hat W^\star(u)$ varies in sign
when $W$ is non-negative\footnote{Such a sign-change may be detected by
using~\eqref{neat} for $u\in[1/2,1]$. The positivity of $\hat
W^\star(u)$ implies that $\frac{6}{v\pi^2}\int_1^2 W(u)du\ge
W(v)$ when $v\in[1,2]$, leading to a contradiction.}.
\subsubsection*{A smoothed setup}
Our analysis revolves around the quantity
\begin{equation}
  \label{eq:19}
  \sum_{q\ge 1}\frac{W(q/Q)}{q}\sum_{a\mode q}
  |S(\varphi,a/q)|^2
\end{equation}
for some weight function $W$ satisfying:

\begin{itemize}
\item[$(W_1)\,\bullet$] The function $W$ is C$\mkern1.5mu{}^3$ over $]-\infty,\infty[$ and
  C$\mkern1.5mu{}^4$ per pieces.
\item[$(W_2)\,\bullet$] It is even and its support lies inside
  $[-\upm,-\lowm]\cup[\lowm,\upm]$.  
\item[$(W_3)\,\bullet$]  We have $\int_0^\infty W(u)du\neq0$.
\end{itemize}
We do not need $W$ to be non-negative, though nothing is made to avoid this
natural condition.
We do not seek generality but on the reverse to restrict ourselves
to as smooth a situation as necessary. 

We define
\begin{equation}
  \label{defIW}
  \ooWo=\sum_{q}\frac{\phi(q)W(q/Q)}{qQ}
  =\frac{6}{\pi^2}\int_0^\infty W(u)du+\Ocal((\log Q)/Q).
\end{equation}
The quantity $\ooWo$ depends on $Q$, but in a very mild manner.

\subsubsection*{First step: an equality via $\delta$-symbol}
The proof of Theorem~\ref{lowerbound} will unfold in four steps.
We start our journey with the following essential formula that is of
independent interest.
\begin{thm}
  \label{Precise}
  When $1/2\le H\le \sqrt{N}/(\log N)^5$ and $\log Q\ll \log N$, we have
  \begin{multline*}
    \sum_{q}\frac{W(q/Q)}{qQ}\sum_{a\mode q} |S(\varphi, a/q)|^2 = \bigl(\ooWo +
    \Ocal(N(QH)^{-1})\bigr) \sum_m|\varphi_m|^2
    \\
    - \sum_{\substack{h\le H}}\frac{1}{h} \sum_{a\mode h}
    \int_{-\infty}^{\infty} \hat W^\star(u)\Bigl|S\Bigl(\varphi,
    \frac{a}{h}+\frac{u}{hQ}\Bigr)\Bigr|^2 du.
\end{multline*}
\end{thm}
The reader will find a refined version for primes in Theorem~\ref{PrecisePrimes}. 
Please note that the
factor $N(QH)^{-1}$ is not polluted by any power of $\log N$ and that
$\hat W^\star(u)$ belongs to $L^1$. The proof shows clearly that a
polarized version is accessible of the same strength, namely:
\begin{multline*}
  \sum_{q}\frac{W(q/Q)}{qQ}\sum_{a\mode q} S(\varphi, a/q)\overline{S(\psi,
    a/q)} 
  = \ooWo\sum_m\varphi_m\overline{\psi_m}
  \\
  - \sum_{\substack{h\le H}}\frac{1}{hQ} \sum_{a\mode h}
  \int_{-\infty}^{\infty} \hat W^\star(u)
  S\Bigl(\varphi,
  \frac{a}{h}+\frac{u}{hQ}\Bigr)
  \overline{S\Bigl(\psi,
    \frac{a}{h}+\frac{u}{hQ}\Bigr)}
  du
  \\
  +\Ocal\bigl(N(QH)^{-1})\|\varphi\|_2\|\psi\|_2\bigr)
\end{multline*}
where $\|\varphi\|_2=\sqrt{\sum_m|\varphi_m|}$ and similarly for
$\|\psi\|_2$. Similar polarized versions are true for
Theorems~\ref{Vic},~\ref{yoddle} and Corollary~\ref{easyy}.  The
beginning of our proof follows closely the one
of B. Conrey \& H. Iwaniec~\cite{Conrey-Iwaniec*03} (which has been for the most part
incorporated in \cite{Conrey-Iwaniec-Soundararajan*12} by B. Conrey,
H. Iwaniec \& K. Soundararajan) and can be
considered as an additive analogue of their result.  Our main new
ingredient at this stage, with respect to this proof, is the use of a
maximal large sieve inequality.  To introduce this part, we got
inspired from another try at a large sieve equality due to W. Duke \& H. Iwaniec and contained
in~\cite{Duke-Iwaniec*92}. The treatment of the finite parts
(meaning: for $h\le H$) diverges from~\cite{Conrey-Iwaniec*03},
and in particular we show that what may appear like two main terms in
the first coarse formula we get in fact cancels out in their leading
contribution. This part of the treatment is similar to what happens
for the $\delta$-symbol of W. Duke, J. Friedlander \& H. Iwaniec in \cite{Duke-Friedlander-Iwaniec*93} (see
also \cite[Section~20.5]{Iwaniec-Kowalski*04} by H. Iwaniec \& E. Kowalski. A more precise
version of this remark is documented Section~\ref{Cancel}).

Since $\hat W^\star(u)$ has its main contribution around~$u=0$, the
sum over~$h$ contributes to the main term only when the sequence
$(\varphi_n)$ accumulates in some arithmetic progression of modulus
$\le H$. When it does \emph{not}, we have the following result that
implies a conditional large sieve equality.

\begin{thm}
  \label{Vic}
  When $\frac12\le H\le \sqrt{N}/(\log N)^5$ and $\log Q\ll \log N$, we have
  \begin{multline*}
  \sum_{q}\frac{W(q/Q)}{qQ}\sum_{a\mode q}
  |S(\varphi, a/q)|^2
  =
  \bigl(\ooWo
  + \Ocal(N(QH)^{-1})\bigr)
  \sum_m|\varphi_m|^2
  \\
  +\Ocal\biggl(
  \sum_{\substack{h\le H}}\frac{N+hQ}{hQ^2}
  \max_{u<v<u+2hQ}
  \sum_{c\mod h}
  \biggl|\sum_{\substack{u<n\le v,\\ n\equiv c[h]}}\varphi_n\biggr|^2
  \biggr).
\end{multline*}
\end{thm}
Recall that the size condition $u,v\le N$ is included in the condition
on the support of $\varphi$.  See Theorem~\ref{betterVic} for a
sharper remainder term.  See also the work \cite{Friedlander-Iwaniec*92} 
of J. Friedlander \& H. Iwaniec, as well as
 \cite[Theorem 2.6]{Ramare*06} for a large sieve equality for
coefficients of a special form (convolution of a shortly supported
sequence with a smooth sequence).  The case $H=1/2$ has also an
interesting methodological consequence.
\begin{cor}
  \label{easyy}
  When $\log Q\ll\log N$, we have
  \begin{equation*}
    \sum_{q}\frac{W(q/Q)}{qQ}\sum_{a\mode q}
  |S(\varphi, a/q)|^2
  =
  (1+\Ocal(N/Q))\ooWo\sum_m|\varphi_m|^2.
  \end{equation*}
\end{cor}

\subsubsection*{Second step: Functional rephrasing}
Corollary~\ref{easyy} describes the situation satisfactorily when $\tau=N/Q$
goes to zero. When $\tau$ is larger, we show that the situation is controlled
by a family of embeddings $(\Pure_{N,h})_h$ of $L^2(\N{N})$ and a family
of self-adjoined nuclear operators $\Vscr_{\tau,h}$ on the subspace
$L^2_*(X_h)$ of $L^2(X_h) $: we endow $X_h=\Z{h}\times[0,1]$ with the
natural probability measure; the space $L^2_*(X_h)$ is the one of functions
from $L^2(X_h)$ whose Fourier transform with respect to the first variable is
supported by $(\Z{h})^*\times[0,1]$, see Section~\ref{infinity} for more
details. We denote by $U_{\tilde h\goes h}$ the orthonormal
projection on this subspace.

Let us define the local embedding $\Pure_{N,h}$. We start by defining
the (nearly) unitary (see Lemma~\ref{iso}) embedding
$\Gamma_{N,h}$ of $L^2(\N{N})$ in $L^2(X_h)$ by:
\begin{equation}
  \label{defFNh}\arraycolsep=1.4pt
  \begin{array}{rcl}
    \Gamma_{N,h}:L^2(\N{N})&\rightarrow&
    \displaystyle L^2(X_h)\\[0.5em]
    \displaystyle \varphi=(\varphi_n)_{1\le n\le N}&\mapsto&
    \begin{array}[t]{rcl}
      \Gamma_{N,h}(\varphi):\Z{h}\times[0,1]&\rightarrow&\mathbb{C}\\
      (b,y)&\mapsto&\displaystyle \varphi_{\sigma_h(b)+h[N'y/h]}
    \end{array}
  \end{array}
\end{equation}
where $\sigma_h(z)$ is the unique integer $b$ in $\N{h}$ that is
congruent to $z$ modulo~$h$; we have set $\varphi_n=0$ when 
the index~$n$ is (strictly) larger than $N$ and
\begin{equation}
  \label{defN'}
  N'=N+\sqrt{N}.
\end{equation}
The embedding we
need is given by
\begin{equation}
  \label{defPure}
  \Pure_{N,h}=U_{\tilde h\va h}\circ \Gamma_{N,h}.
\end{equation}
This is to be compared with the case of integers where we send $\mathbb{Z}$
inside $\mathbb{Z}_p$ for every prime~$p$, though we have here an ``infinite
place'' for each modulus~$h$ (this is the factor $[0,1]$) and that we may not
rely on multiplicativity. It would be interesting to show that the diagonal
embedding $\varphi\mapsto(\Pure_{N,h}(\varphi))_h$ has a dense range, as in
the adelic case. The situation is somewhat more intricate because of
the dependence in $N$.
We next define the one-parameter family of operators~$\Vscr_{\tau,h}$ by
\begin{equation}
  \label{defVscrtauh}
  \Vscr_{\tau,h}(G)(b,y)=\int_0^1 
  G(b,y')W^\star\Bigl(\frac{\tau(y-y')}{h}\Bigr)dy'.
\end{equation}
They are shown to be compact symmetric nuclear operators in
Theorem~\ref{nuclear} and to verify a Mercer like theorem (see
Theorem~\ref{Mercer}).  The fundamental formula
is~\eqref{formula11bb-1} which we repeat here: \smallskip

\noindent\fbox{\vbox{%
\begin{multline*}
  \tag{\ref{formula11bb-1}}
  \sum_{q}\frac{W(q/Q)}{qQ}\sum_{a\mode q}
  |S(\varphi, a/q)|^2
  =
  \ooWo\|\varphi\|_2^2(1+\Ocal(\tau/H))
  \\
  -
  N\sum_{h\le H}\frac{\tau}{h}
  \bigl[\Pure_{N,h}(\varphi)|\Vscr_{\tau,h} \Pure_{N,h}(\varphi)\bigr]_{h\times[0,1]}
  \\
  (H\ll N^{1/8}(\log N)^{-3/2},\tau=N/Q\ll H, Q\ll N^2)
  .
\end{multline*}
}}
\smallskip

\noindent
Allowing $H$ to be as large as a power of $N$ requires quite some
efforts and we have to rely on te moe technical
formula~\eqref{formula51} rather than on the simplified form given in
Theorem~\ref{Precise}. Ideally, we should be able to allow $H$ roughly as large as
$\sqrt{N}$.

\subsubsection*{Analysis of a class of difference operators}
We treat in Section~\ref{ACDO} the analysis of the intervening family of operators in
an abstracted setting.
For a function $V$ satisfying the regularity assumptions $(R_1)$,
$(R_2)$ and $(R_3)$, we define
\begin{equation}
  \label{defV0}
  \Vscr_0:\quad G\in L^2([0,1])\mapsto \biggl(y\mapsto\int_0^1 G(y')V(y-y')dy'\biggr)
\end{equation}
Assumptions $(R_1)$, $(R_2)$ and $(R_3)$ indeed hold when
$V(y)=W^\star(\tau y/h)$. 
It is classical theory that $\Vscr_0$ is a compact Hilbert-Schmidt operator,
see for instance \cite[Theorem 7.7]{Gohberg-Goldberg-Krupnik*00}. Let
$(\lambda_\ell,G_\ell)_\ell$ be a complete orthonormal system of
eigenvalues / eigenfunctions, ordered with non-increasing
$|\lambda_\ell|$. 
The Fredholm equation $\lambda G(y')=\int_0^1 K(y',y)G(y)dy$ has been
intensively studied. It is not the purpose of this paper to introduce  to this
theory, a task for which it is better to read the complete and
classical~\cite{Gohberg-Krein*69}, or the more modern
\cite{Gohberg-Goldberg-Krupnik*00}. 
Kernel of type $V(y'-y)$ are often called
\emph{difference kernel}, and lead to operators that are distinct from
convolution operators as the integration and definition interval is \emph{not} the
whole real line. The book \cite{Sakhnovich*15} is dedicated to the operators
built from such kernels. The book \cite{Cochran*72} contains also many useful
informations. 

Here is a summary of what we prove in Section~\ref{ACDO}.
\begin{thm}
  \label{recapdivertimento}
  The operator $\Vscr_0$ is nuclear.  Given a complete collection
  $(\lambda_\ell,G_{\ell})_\ell$ of non-zero eigenvalues /
  eigenvectors, arranged with non-increasing $|\lambda_\ell|$ and
  normalized by
  $\int_0^1|G_\ell(t)|^2dt=1$, we have
  the three following properties:
  \begin{itemize}
  \item (Explicit nuclearity) $\displaystyle
    \sum_{\ell\ge 1}|\lambda_\ell|\ll
    \|V\|_2e^{-c'\sqrt{1+\|V\|_\infty\|V'\|_1/\|V\|_2^2}}$ for some
    positive constant $c'$ depending only on $A$, $B$ and $c$. The
    notation $\|V'\|_1$ stands for the total variation.
  \vspace*{-10pt}
  \item (Mercer like property) $\displaystyle
    V(y'-y)=\sum_{\ell\ge 1}\lambda_\ell G_\ell(y')G_\ell(y)$ uniformly.
  \vspace*{-10pt}
  \item (Lidskii's Theorem) $\displaystyle \sum_{\ell\ge1}\lambda_\ell=0$.
  \end{itemize}
\end{thm}
This is proved in Theorem~\ref{nuclear} and~\ref{Mercer}. These
properties shows that this class of operators is indeed very
regular. We recall that the Mercer Theorem concerns similar operators
but having a \emph{non-negative} reproducing kernel. On integrating
the case $y=y'$ of the Mercer like property, we recover the third property.

\subsubsection*{Third step: Spectral decomposition of the large sieve}
\noindent 
\begin{thm}
  \label{yoddle}
  Assume that $\sqrt{N}\le Q\le N$. There exist two positive constants $c_0$
  and $c_3$
  such that the following holds. For each $\tau=N/Q$ and integer
  $h\ge1$, let $(G_{\ell, \tau/h},\lambda_{\ell}(\tau/h))_\ell$ be a
  complete family of two by two orthonormal eigenfunctions of \eqref{defVscrtauh} coupled
  with their respective non-zero eigenvalues.
  These eigenfunctions are all continuous and of bounded variations.
  The sequence $(\lambda_{h,\ell}(\tau))_{\ell\ge1}$
  is arranged in non-increasing absolute
  value, and satisfies
  $\lambda_{\ell}(\tau/h)\ll 1/\sqrt{\ell}$ uniformly in~$h$
  and~$\tau$.
  We also have
  \begin{equation}
  \label{t2}
    \sum_{\ell\ge1}\lambda_{\ell}(\tau/h)=0,\
    \sum_{\ell\ge1}|\lambda_{\ell}(\tau/h)|<\infty,\ 
    \sum_{\ell\ge1}|\lambda_{\ell}(\tau/h)|^2=
     2\int_0^1
    W^\star\Bigl(\frac{\tau y}{h}\Bigr)^2(1-y)dy
  \end{equation}
  and this last value is bounded uniformly in $\tau$. Under
  the Riemann Hypothesis, we also have
  $\sum_{\ell\ge1}|\lambda_{\ell}(\tau/h)|^p<\infty$ for any $p>4/5$.
  For any sequence of complex numbers $\varphi$, any $L\ge1$, any
  $H\ll N^{1/8}(\log N)^{-3/2}$ and any $\xi\in[0,1]$, we have
  \begin{multline*}
  \sum_{q}\frac{W(q/Q)}{qQ}\sum_{a\mode q}
  |S(\varphi, a/q)|^2
  =
  \ooWo\|\varphi\|_2^2
  \\
  -
  \frac{1}{N}\mkern-20mu\mathop{\sum_{h\le H}
    \sum_{\substack{\ell\le L}}}\limits_{|\lambda_{\ell}(\tau/h)|\ge
      \xi\eta_0(N)}
  \mkern-20mu
  \frac{\tau}{h}\lambda_{\ell}(\tau/h)
  \sum_{a\mod^*h}\biggl|
  \sum_{n\le
    N}\varphi_n G_{\ell,\tau/h}\Bigl(\frac{n}{N}\Bigr)e\Bigl(\frac{na}{h}\Bigr)
  \biggr|^2
  \\+\Ocal\biggl( 
  \Bigl(
  \frac{\log H}{L}+\frac{1}{H} 
  +
  \xi \eta_0(N)
  \Bigr)
  \tau
  \|\varphi\|_2^2\biggr)
  \end{multline*}
  where  $\eta_0(N)=\exp-c_3\sqrt{\log N}$. We have furthermore
  \begin{equation*}
  \mathop{\sum_{h\le H}\sum_{\ell\le L}\sum_{a\mod^*h}}_{|\lambda_{\ell}(\tau/h)|\ge
      \eta_0(N)}
  \biggl|
  \sum_{n\le
    N}\varphi_n
  G_{\ell,\tau/h}\Bigl(\frac{n}{N}\Bigr)e\Bigl(\frac{na}{h}\Bigr)
  \biggr|^2
  \le N \|\varphi\|^2_2\bigl(1+ H^2L\eta_0(N)\bigr).
  \end{equation*}
  When $W$ is non-negative, the one-sided inequality
  $(\tau/h)\lambda_{h,\ell}\le \ooWo+o(1)$ holds true, where
  $o(1)$ is here a function of $Q$ that goes to~0 with~$1/Q$.
\end{thm}

We prove that infinitely many
$\lambda_{\ell}(\tau/h)$ are positive (resp. negative),
once $h$ is also allowed to vary; see end of
Subsection~\ref{QF}. When $W$ is further assumed to be non-negative,
Theorem~\ref{bounds} shows that
$(\tau/h)\lambda_{\ell}(\tau/h)\le
\ooWo+o(1)$.
The parameter $\xi$
above has only been introduced for flexibility purpose, in case one
needs a lower bound that is independent on~$N$. 

\subsubsection*{Fourth step: Uncertainty principle and eigenvalues properties}

A closer study of the eigenvalues that uses F.I. Nazarov's version
\cite{Nazarov*92} of the uncertainty principle combined with some positivity
argument leads to the following.
\begin{thm}\label{lowerls}
  For any non-negative $W$ satisfying the above conditions there exist
  $c_4,c_6,c_7>0$ such that we have, for any $H\le
  \exp(c_6\sqrt{\log N})$ and any $Q\in[N\exp(-c_6\sqrt{\log
    N}),N^2]$, 
  \begin{equation*}
    (\tau/\ooWo)|\lambda_{h,\ell}|/h\le
  1-c_6e^{-c_4\tau/h}+\Ocal(\exp-c_7\sqrt{\log N})
  \end{equation*}
  for any $h\le H$, any
  $\ell\ge1$ and with $\tau=N/Q$.
\end{thm}
P. Jaming tells me that he believes $c_4=120$ to be an admissible choice.

\subsubsection*{Arithmetical consequences}

\begin{cor}
  \label{ees}
  For every $\epsilon>0$, and every $N\ge1$ and $Q\ge1$, there exist a
  constant $c_4$ and 
  a subspace of dimension $\Ocal(\tau^2/[\epsilon^2\log(1/\epsilon)])$ such that we have,
  for any $(\vp_n)$ orthogonal to this subspace and when $\log Q>c_4\log^2 (N/Q)$, 
  \begin{equation*}
    (1-\epsilon)\sum_m|\varphi_m|^2\le
    \sum_{q}\frac{W(q/Q)}{qQ\ooWo}\sum_{a\mode q}
    |S(\varphi, a/q)|^2\le (1+\epsilon)\sum_m|\varphi_m|^2.
  \end{equation*}
  Moreover, when $\tau\asymp 1$ and for every integer $K\ge1$, there exist
  $\epsilon_0>0$ depending only on 
  $\tau$ and $K$, and $2K$ unitary sequences $(\alpha_k)_{k\le K}$ and
  $(\beta_k)_{k\le K}$, two by two almost orthogonal in the sense that
  \begin{equation*}
    \forall\gamma,\gamma'\in\{\alpha_k\}\cup\{\beta_k\},\quad
    [\gamma,\gamma']_N
    =
    \delta_{\gamma=\gamma'}
    +
    \Ocal\bigl(
    \exp\bigl(-c_4\sqrt{\log N}\bigr)
    \bigr),
  \end{equation*}
  and such that, on one side, we have
  \begin{equation*}
    \sum_{q}\frac{W(q/Q)}{qQ\ooWo}\sum_{a\mode q}
    |S(\alpha_k, a/q)|^2 > (1+\epsilon_0)\sum_m|\alpha_{k,m}|^2
  \end{equation*}
  while on the other side, we have
  \begin{equation*}
    \sum_{q}\frac{W(q/Q)}{qQ\ooWo}\sum_{a\mode q}
    |S(\beta_k, a/q)|^2 < (1-\epsilon_0)\sum_m|\beta_{k,m}|^2.
  \end{equation*}
\end{cor}
The orthogonality is according to the hermitian product defined by
\begin{equation}
  \label{defscalarN}
  [\varphi,\psi]_N=\frac{1}{N}\sum_{1\le n\le N}\varphi_n\overline{\psi_n}.
\end{equation}
The sequences $(\alpha_k)$ and $(\beta_k)$ are pull-backs of
eigenvectors. Note that the pulling-back process depends on $N$ but that the
eigenvectors do not. They are very regular and do not result from some
exotic construction; in particular they are uniformly
bounded and there exists $\epsilon>0$ such that $\{n\le
N,|\alpha_{k,n}|\ge\epsilon\}$ is a set of density (in short: their
``essential support'' is a set of density).

\subsection*{Notation}
We note the Mellin transform by $\check{W}(s)=\int_0^\infty
W(t)t^{s-1}dt$ and the Fourier transform by
$\hat{W}(u)=\int_{-\infty}^{\infty}W(t)e(-ut)dt$. Several other
transforms of $W$ will be used, $W^\sharp$, $W$, $\tilde{W}$,
$W^\star$ and $W^{\star\star}$; they are described in
section~\ref{SomeTransforms}.  We note here that the transform
$W^\sharp$ is very close to what appears in~\cite[section 20.5,
(20.145)]{Iwaniec-Kowalski*04} provided the changes of notation is
incorporated: our $W(y)$ is their $w(y/C)$.
 We recall that
$\|\sigma_\infty'(\psi,\cdot)\|_{1,N}=\int_1^N|\sigma_\infty'(\psi,t)|dt$.
We denote by $a_{|_t}=(a_n)_{n\le t}$ the truncated sequence. We also
define
\begin{equation*}
  \gL(u)=\exp\sqrt{\log(2+u)}.
\end{equation*}
We denote the Euler totient function by $\phi$ and distinguish it from
the sequence by using a different script for the latter,
namely~$\varphi$.  We use the following norms:
\begin{equation}
  \label{defNorms}
  \|f\|_{1,N}=\int_1^N|f(t)|dt, \quad
  \|f\|_{\infty,N}=\max_{1\le t\le N}|f(t)|.
\end{equation}

\subsection*{Comment}
An improved version of this paper, with the same title, will appear
(2024 or 2025) in the Bulletin of the French Mathematical Society.
\section{Related works}

\subsubsection*{Influence of the Riemann Hypothesis}

Under the Riemann Hypothesis (and not the Generalized one as
one may believe), the proof we present allows to select $Q$ as small
as $N/(\log N)^{1-\varepsilon}$ for any positive~$\varepsilon$.  The
coefficient $e^{-c N/Q}$ may be questioned and may well be
superfluous in this range. 
\subsubsection*{Eigenvalues considerations when $Q\ll\sqrt{N}$}

The eigenvalues of the quadratic form
$\sum_{q\le Q}\sum_{a\mode q}|S(\varphi,a/q)|^2$ are well understood
when $Q=o(\sqrt{N})$, see the paper of I. Kobayashi
\cite{Kobayashi*73} and this quantity is expected to behave like a
Riemann sum when $N=o(Q)$ (Corollary~\ref{easyy} below gives a precise
form to this statement), but the behavior in the range
$Q\in[\sqrt{N},cN]$ (for any positive constant $c$) is still
mysterious. When $Q\sim\sqrt{N}$, F. Boca and M. Radziwi\l\l{} have
shown in \cite{Boca-Radziwill*16} by a very delicate analysis that the
distribution of the eigenvalues of this quadratic form tend to a
limiting distribution, henceforth proving a conjecture made
in~\cite{Ramare*07a}. In fact, though this
went unnoticed by the authors, the paper \cite{Chan-Kumchev*12} of
T.H. Chan \& A.V. Kumchev can be read as also providing some
informations on the eigenvalues in the case $Q\sim\sqrt{N}$.
The values for the even moments of this limit
distribution reveals that it is \emph{not} a classical distribution,
confirming
what the (rather limited) computations from~\cite{Ramare*09a}.

\subsubsection*{Eigenvalues considerations when $Q\ge N$}
H. Niederreiter evaluated in
\cite{Niederreiter*73} the discrepancy of the Farey
sequence, a study refined by F. Dress in~\cite{Dress*99}, and this, together
with the Koksma-Hlawka's inequality, proves immediately that
\begin{equation*}
  \sum_{q\le Q}\sum_{a\mode q}
  |S(\varphi, a/q)|^2= \sum_{q\le Q}\phi(q)\sum_n|\varphi_n|^2(1+\Ocal(N/Q))
\end{equation*}
in very much the same way P. Gallagher in \cite{Gallagher*67} derived the large
sieve inequality. Note that the arithmeticity of the Farey sequence is only
mildly used: a discrepancy estimate is enough.

\part{A large sieve equality}

\section{Large sieve ingredients}

We adapt here the proof of S. Uchiyama \cite{Uchiyama*72}
concerning the maximal large sieve to get a result which is a (weak)
additive analogue of a result of~P.D.T.A. Elliott \cite{Elliott*91}. This is 
\cite[Lemma 1]{Elliott*85b} or 
\cite[Chapter 29,
exercise 3, page 254]{Elliott*97}.

\begin{lem}
  Let $(x_d)_{d\le D}$ be a $\delta$-spaced sequence of points of
  $\mathbb{R}/\mathbb{Z}$. We have
  \begin{equation*}
    \sum_{d\le D}
    \max_{u<v\le u+L}\bigl|\sum_{u<n\le v}\varphi_me(m x_d)\bigr|^2
    \le (L+2\delta^{-1}\log(e/\delta)) \sum_m|\varphi_m|^2.
  \end{equation*}
\end{lem}

Here is the version we shall use.
\begin{lem}\label{LSM}
  We have
  \begin{equation*}
    \sum_{q\le Q}\sum_{a\mode q}
    \max_{u<v\le u+L}\bigl|\sum_{u<n\le v}\varphi_me(ma/q)\bigr|^2
    \ll (L+Q^2\log Q) \sum_m|\varphi_m|^2.
  \end{equation*}
\end{lem}

\section{A functional transform}

The transform we investigate here is given by
\begin{equation}
  \label{defTransform}
  \Transform{W}(z)=
  \frac{-1}{2i\pi}\int_{-i\infty}^{i\infty}
  \check{W}(s)\frac{\zeta(1-s)}{\zeta(1+s)} z^{-s}ds.
\end{equation}
Please note that $|{\zeta(1-s)}/{\zeta(1+s)}|=1$ on the line $\Re s=0$.
This transform of $\check W$ is already the one the occurs in
\cite{Ramare*07a}, see for instance equation numbered~$(48)$ there, and in
\cite{Chan-Kumchev*12}, see their equation~$(4.19)$. We keep the same
hypothesis as before for $W$. In particular, it is compactly supported and
$\check{W}(s)\ll (1+|s|)^{-4}$. We follow \cite[Section 9]{Ramare*07a} pretty closely.
We start by recalling a handy form of the complex Stirling formula.
\begin{lem}[Uniform complex Stirling formula]
  \label{Stirling}
  Let $\varepsilon\in]0,1]$ and a compact subset $\mathcal{A}$ of $\mathbb{C}$
  be fixed. In the domain $|\arg z|\le \pi-\ve$ and $|z|\ge1$, we have
  \begin{equation*}
    \Gamma(z+a)=\sqrt{2\pi}e^{-z}z^{z+a-1/2}
    \bigl(1+\Ocal(1/|z|)\bigr).
  \end{equation*}
  uniformly for $a\in\mathcal{A}$.
\end{lem}
As a (classical) conclusion and taking $z=it$ in the above, we find that
\begin{equation}
  \label{csf}
  |\cos(\sigma+it)\Gamma(\sigma+it)|=\sqrt{\pi/2}|t|^{\sigma-\1/2}\bigl(1+\Ocal(1/|t|)\bigr)
\end{equation}
uniformly in any domain $\sigma_1\le \sigma\le \sigma_2$ and $|t|\ge1$. 

\subsection*{Isolating the arithmetical behavior}

We proceed as in \cite{Ramare*07a} and appeal to the functional equation of
the Riemann $\zeta$-function (see
\cite{Titchmarsh*51} or \cite{Iwaniec-Kowalski*04}) which may be written as
\begin{equation}
  \label{func}
  \zeta(1-s)=2^{1-s}\pi^{-s}\cos(\pi s/2)\Gamma(s)\zeta(s).
\end{equation}
To do so we first shift the line of integration in \eqref{defTransform} to
$\Re s=9/8$. Since $|\zeta(-\sigma+it)|\ll_\varepsilon
(1+|t|)^{(1+\sigma)/2+\varepsilon}$ when $\sigma\ge0$ and for any
$\varepsilon>0$, it is enough to assume that $\check{W}(s)\ll (1+|s|)^{-2}$ to
ensure the convergence of our integrals. Since the line shifting does not meet
any pole, we get
\begin{align}
  \Transform{W}(z)
  &=\notag
  \frac{-1}{i\pi}\int_{\frac98-i\infty}^{\frac98+i\infty}
  \check{W}(s)\frac{\cos(\pi s/2)\Gamma(s)\zeta(s)}{\zeta(1+s)}
  (2\pi z)^{-s}ds,
  \\&=
  2\sum_{n\ge1}\frac{\phi(n)}{n}
  \PartialTransform{W}(2\pi nz)
\end{align}
where
\begin{equation}
  \label{defPartialTransform}
  \PartialTransform{W}(u)=
  \frac{-1}{2i\pi}\int_{\frac98-i\infty}^{\frac98+i\infty}
  \check{W}(s)\cos(\pi s/2)\Gamma(s)
  u^{-s}ds.
\end{equation}

\subsection*{A bound at infinity}

We infer from the estimate \eqref{csf}
that the line of integration in \eqref{defPartialTransform}
can be pushed up to $\Re s=7/2-\varepsilon$ and thus
\begin{equation}
  \label{FirstBound}
  \PartialTransform{W}(2\pi nz)\ll_\varepsilon (nz)^{-7/2+\varepsilon}.
\end{equation}
Here is the main conclusion of this part.
\begin{lem}
  \label{eq:17}
  We have $\Transform{W}(z)\ll_{\varepsilon} z^{-7/2+\varepsilon}$, for any $\varepsilon>0$.
\end{lem}

\subsection*{A real-valued formula}

The next step is to proceed as in section~9 of \cite{Ramare*07a}, which we
only sketch here. We employ equation $(35)$ therein:
\begin{equation}
  \label{revMT}
  \cos\frac{\pi s}{2}\,\Gamma(s)
  =\int_0^{\infty}\cos(y)y^{s-1}dy
  =\int_0^{\infty}\cos(y)y^{s}dy/y
\end{equation}
valid for $0<\Re s< 1$ to infer that
\begin{align*}
  \PartialTransform{W}(u)
  &=-
  \frac{1}{2i\pi}\int_{\frac12-i\infty}^{\frac12+i\infty}
  \check{W}(s)\cos(\pi s/2)\Gamma(s)
  u^{-s}ds
  \\&=-
  \int_0^{\infty}\cos(y)\frac{1}{2i\pi}\int_{\frac12-i\infty}^{\frac12+i\infty}
  \check{W}(s)(u/y)^{-s}\frac{dsdy}{y}
  \\&=-
  \int_0^{\infty}\cos(y)W(u/y)dy/y
\end{align*}
by Mellin inversion formula. This yilds formula~\eqref{explicitTransform-1}.

\section{More auxiliary functional transforms}
\label{SomeTransforms}

Several functional transforms of our bump-function $W$ will occur. We
have already seen $\Transform{W}$ and $\hat{W}^\star$
at~\eqref{explicitTransform-1} and~\eqref{neat}. These two functions
are central in our work, but it is expedient to introduce several
others. We start with the couple
\begin{equation}
      W^\sharp(y)=  \label{defWsharp}
      \sum_{k\ge1}\frac{W(y/k)}{k},\quad
      W^\flat(y)=
      \sum_{f\ge1}\frac{W(yf)}{f}.
\end{equation}
We show in Lemma~\ref{evalWflat} that $W^\flat(y)=J(W)+ \Ocal(y)$ where
\begin{equation}
   \label{defTW}
   J(W)=\int_{0}^\infty \frac{W(u)du}{u}.
 \end{equation}
 When $y$ is small as in our case of application, the approximation
 of $W^\flat$ by~$J(W)$ is efficient. The proof will then lead us to understand $W^\sharp-J(W)$, a quantity we call $-\tilde{W}$, i.e.
 \begin{equation}
   \label{defWtildebis}
   W^\sharp(y)=J(W)-\tilde{W}(y).
 \end{equation}
 The situation is there more difficult than with $W^\flat$, in
 particular because $W^\flat(y)$ is \emph{not} small when $y$ is small
 but takes the constant value~$J(W)$! See Lemma~\ref{boundWtilde}.  As it turns
 out, we do not need to grasp $\tilde{W}$ but the average
 \begin{equation}
  \label{defWstarC}
  W^\star_C(z)=\sum_{1\le c\le C}\frac{\mu(c)}{c}\tilde W(cz).
\end{equation}
The value for small $z$, i.e. when $|z|\le \lowm$, is now
$J(W)\sum_{1\le c\le C}\mu(c)/c$ which tends to~0 when $C$ is
large. The rate of convergence is fast enough on the Riemann
Hypothesis, but rather slow otherwise. As a consequence, we have to
treat this point with care. In particular, we want to replace $C$ by
$\infty$ and still save a power of $C$.  We have already defined
$W^\star$ at \eqref{defTransform} and Lemma~\ref{formulaWstar} will
show that both definitions coincide.  Let us start our journey.

\subsection{Approximating $W^\flat$}

The transform $W^\flat$ is also studied in \cite[section
20.5]{Iwaniec-Kowalski*04}: the function $V(z)$ defined there in $(20.143)$
corresponds to $J(W)-W^\flat(z)$ where one should change $W(y)$ into $w(y/C)$
(albeit the trivial facts that $w$ is supported on $[C,2C]$, while our $W$ is
supported on $[\lowm,\upm]$ and extended to the negative real axis by evenness).

\begin{lem}
  \label{evalWflat}
  Assume that $|\hat{W}(u)|\ll1/(1+|u|)^2$.
  We have, when $z>0$,
 \begin{equation}
   \label{myWflat}
   W^\flat(z)=\sum_{\substack{ f\ge1}}\frac{W(zf)}{f}
   =J(W)+ \Ocal(z).
 \end{equation}
 with $J(W)$ being defined at~\eqref{defTW}.
\end{lem}
In practice, $z$ is small ($\le CE/Q$). The proof we present uses the Fourier
transform but one could also use the Mellin transform.

\begin{proof}
  We introduce Fourier transforms to write
  \begin{align*}
    W^\flat(z)
    &=\int_{-\infty}^{\infty}\hat{W}(u)\sum_{f\ge1}\frac{e(fuz)}{f}du
    =-\int_{-\infty}^{\infty}\hat{W}(u)\log(1-e(zu))du
    \\&=-\int_{-\infty}^{\infty}\hat{W}(u)
    \bigl(\log|2\sin(\pi zu)|+i\pi(\{zu\}-\tfrac12)\bigr)du
    \\&=-\int_{-\infty}^{\infty}\hat{W}(u)
    \bigl(\log|2\sin(\pi zu)|+i\pi B_1(zu)\bigr)du.
  \end{align*}
  For the sake of the evaluation next to $z=0$, it is better to adopt the
  expression
  \begin{equation*}
    W^\flat(z)
    =
    -\int_{-\infty}^{\infty}\hat{W}(u)
    \bigl(\log\frac{|2\sin(\pi zu)|}{\pi z|u|}
    +\log(\pi z|u|)+i\pi (\{zu\}-\tfrac12)\bigr)du.
  \end{equation*}
  which we may simplify, with $\int_{-\infty}^{\infty}\hat{W}(u)du=W(0)=0$, into
  \begin{equation*}
    W^\flat(z)
    =
    -2\int_{0}^{\infty}\hat{W}(u)
    \log |u|\,du
    -\int_{-\infty}^{\infty}\hat{W}(u)
    \Bigl(\log\frac{|2\sin(\pi zu)|}{\pi z|u|}+i\pi \{zu\}\Bigr)du.
  \end{equation*}
  We split the integral according to whether $|u|\le 1/z$ or not. In both
  cases we use $|\hat{W}(u)|\ll1/(1+|u|)^2$ and bound $\log\frac{|\sin(\pi
    zu)|}{\pi zu}$ by $\Ocal(zu)$ when $|u|\le 1/z$ and by $\log(|zu|+1)$
  otherwise.

  We proceed by getting a simpler form for $-2\int_{0}^{\infty}\hat{W}(u)
  \log |u|\,du$. We readily check that
  \begin{align*}
    \int_{0}^{L}\hat{W}(u)\log |u|\,du
    &=
    2\int_{0}^{L}\int_1^2 W(t)\cos(2\pi ut)dt\log |u|\,du
    \\&=
    2\int_1^2 W(t)
    \biggl(
    \Bigl[\frac{\sin(2\pi ut)}{2\pi t}\log |u|\Bigr]_0^L
    -\frac{1}{2\pi t}\int_{0}^{L}\frac{\sin(2\pi ut)}{u}du
    \biggr)
    \\&\rightarrow\frac{-1}{2}\int_1^2 \frac{W(t)dt}{t}
  \end{align*}
  therefore concluding the proof of our lemma.
\end{proof}

\subsection{From $W^\sharp$ to $\tilde W$}

In this part, we start from the definition of $\tilde{W}$ provided by~\eqref{defWtildebis} and we reach the definition~\eqref{defWtilde} given below.
With $v>0$ fixed, we define 
\begin{equation}
  \label{eq:2}
  f(t)=W\bigl({v}/{t}\bigr)/t.
\end{equation}
We simply write  when $v>0$
\begin{align*}
  \sum_{g\ge1}f(g)
  &=
  -\sum_{g\ge1}\int_g^\infty f'(t)dt
  =
  -\int_0^\infty [t]f'(t)dt
  \\&=
  \int_0^\infty f(t)dt
  +
  \int_0^\infty \{t\}f'(t)dt
  \\&=
  \int_0^\infty W(u)\frac{du}{u}
  -
  \int_0^\infty \{t\}
  \biggl(\frac{vW'(v/t}{t^3}+\frac{W(v/t)}{t^2}\biggr)dt
  \\&=
  J(W)
  -
  \frac{1}{v}\int_0^\infty \{v/u\}(W'(u)u+W(u))du.
\end{align*}
This establishes Eq. \eqref{defWtilde}.
The condition $v>0$ has been used on the last line: when $v<0$, we should
reverse the integration path, or divide by $|v|$ instead of by~$v$. 

\subsection{Treatment of $\tilde W$}

Define
\begin{equation}
  \label{defWtilde}
  \tilde W(z)
  =\frac{1}{|z|}\int_0^\infty \{z/u\}(uW'(u)+W(u))du.
\end{equation}
The expression $\tilde W(z)=\int_0^\infty \{1/v\}(vzW'(vz)+W(vz))dv$ shows
that $\tilde W$ is an even function\footnote{Still reading \cite[Section
  20.5]{Iwaniec-Kowalski*04} by H. Iwaniec \& E. Kowalski, we find that our $\tilde W$ satisfies $\tilde
  W(z)=(C/|z|)\int_0^\infty\{vz/C\}(W(C/v)/v)'dv$, and is thus like their
  $W(C/z)$.}.

\begin{lem}
  \label{boundWtilde}
  The function $\tilde W$ is C${}^1$ and C${}^2$ per pieces, and both
  derivatives are bounded.

  When $|z|\le \lowm$, we have $\tilde W(z)= J(W)$.

  When $|z|\ge \lowm$, we have $\tilde W(z)\ll 1/z^2$.
\end{lem}

\begin{proof}
  Eq.~\eqref{defWtilde} shows that the first part of the Lemma, by
  distinguishing whether $|z|>\lowm$ or not.

  When $z\in[0,\lowm)$, then $z/u\in[0,\lowm)$ when $u$ lies in the support of $W$,
  which implies that $\{z/u\}=z/u$ in this case. Hence the first equality.  We
  can furthermore write, when $z\neq0$, and with $t=z/u$, and with
  $B_2^*(t)=\int_0^t B_1(v)dv$:
\begin{align*}
  \tilde W(z)
  &=\int_0^\infty (\{t\}-\tfrac12)(zt^{-1}W'(z/t)+W(z/t))dt/t^2
  \\&=\int_0^\infty B_2^*(t)
  (4zt^{-2}W'(z/t)+2W(z/t)t^{-1}+z^2t^{-3}W''(z/t))dt/t^2
  \\&=z^{-2}\int_0^\infty B_2^*(z/u)(4u^2W'(u)+2uW(u)+u^3W''(u))du
\end{align*}
from which the bound claimed in the lemma follows readily.
\end{proof}

\subsection{Study of $W^\star_C$ and $W^\star$}

The function $W^\star_C(z)$ is even since so is $\tilde W$.
Lemma~\ref{boundWtilde} tells us that this function is constant when $|z|\le
\lowm/C$, with value $J(W)\sum_{1\le c\le C}{\mu(c)}/{c}$. We can even select
$C=\infty$ in which case we write simply $W^\star$:
\begin{equation}
  \label{defWstar}
  W^\star_\infty(z)=W^\star(z)=\sum_{c\ge1}\frac{\mu(c)}{c}\tilde W(cz).
\end{equation}
The next expression of $W^\star$ will in particular establish that
$W^\star$ is continuous at $z=0$ where we have~$W^\star(0)=0$.

\begin{lem}
  \label{formulaWstar}
  We assume that $W$ is at least C$\mkern1.5mu{}^2$. We have, when $\ve>0$ and
  $z>0$,
  \begin{equation*}
    W^\star_C(z)=
    J(W)\sum_{c\le C}\frac{\mu(c)}{c}
    -\frac{1}{2i\pi}\int_{\ve-i\infty}^{\ve+i\infty}\check{W}(s)
    \zeta(1-s)\sum_{c\le C}\frac{\mu(c)}{c^{1+s}}z^{-s}ds
  \end{equation*}
  where $\check{W}(s)=\int_0^{\infty}W(x)x^{s-1}dx$ is the Mellin
  transform of $W$.  When $C=\infty$, the expression above is correct
  provided we select $\ve=0$ and replace $\sum_{c\ge
    1}{\mu(c)}/{c^{1+s}}$ by $1/\zeta(1+s)$.
\end{lem}

\begin{proof}
  We first reduce the case $C=\infty$ to the case $C$ finite. On using
  $\{z\}=z-[z]$, we get
  \begin{align*}
    W^\star(z)
    &=\sum_{c\ge1}\frac{\mu(c)}{c^2z}
    \int_0^\infty \{zc/u\}(uW'(u)+W(u))du
    \\&=
    \lim_{C'\rightarrow\infty}
    \biggl(
    \sum_{c\le C'}\frac{\mu(c)}{cz}
    \int_0^\infty \frac{uW'(u)+W(u)}{u}du
    \\&\qquad\qquad-
    \sum_{c\le C'}\frac{\mu(c)}{c^2z}
    \int_0^\infty \sum_{d\le zc/u}1\,(uW)'(u)du
    \biggr)
    \\&=
    -\lim_{C'\rightarrow\infty}
    \sum_{c\le C'}\frac{\mu(c)}{c}\sum_{d\ge1}\frac{W(zc/d)}{d}.
  \end{align*}
  We introduce the Mellin transform of $W$ and write
  \begin{align*}
    \sum_{d\ge1}\frac{W(zc/d)}{d}
    &=\frac{1}{2i\pi}\int_{-1-i\infty}^{-1+i\infty}\check{W}(s)
    \zeta(1-s)(zc)^{-s}ds
    \\&=\check{W}(0)+
    \frac{1}{2i\pi}\int_{\ve-i\infty}^{\ve+i\infty}\check{W}(s)
    \zeta(1-s)(zc)^{-s}ds
  \end{align*}
  which gives us (note that $J(W)=\check{W}(0)$)
  \begin{equation*}
    W_C^\star(z)=
    J(W)\sum_{c\le C}\frac{\mu(c)}{c}
    -\frac{1}{2i\pi}\int_{\ve-i\infty}^{\ve+i\infty}\check{W}(s)
    \zeta(1-s)\sum_{c\le C}\frac{\mu(c)}{c^{1+s}}z^{-s}ds
  \end{equation*}
  hence the expression given, seeing that the pole of~$\zeta(1-s)$ cancels out
  with the zero of~$1/\zeta(1+s)$ at~$s=0$ and that $\check{W}(s)$ is
  $\Ocal(1/(1+|s|)^2)$.
\end{proof}

\begin{lem}
  \label{form}
  For $\Re s\in(-1,0)$, we have
  \begin{equation*}
    \int_0^1z^s\cos(2\pi z)dz-s\int_1^{\infty}
    \frac{z^{s-1}\sin(2\pi z)dz}{2\pi}
    =
    (2\pi)^{-s-1}\Gamma(s+1)\cos\frac{\pi (s+1)}{2}.
  \end{equation*}
\end{lem}

\begin{proof}
  We call the left-hand side $j(s)$. 
  It is not difficult to see that (this is how is occurs below)
  \begin{equation*}
    j(s)=\int_0^\infty z^s\cos(2\pi z)dz
  \end{equation*}
  and is thus the Mellin transform of $\cos(2\pi z)$.
  On looking at \cite[(21), page
  319]{Erdelyi-Magnus-Oberhettinger-Tricomi*54}, we readily discover that,
  when $\Re s\in(-1,0)$ (note the shift or $+1$ between the $s$ variable
  $j(s)$ and the one of the table we refer to), the above formula follows.
  Giving a full proof is not difficult by using $\cos w=(e^{iw}+e^{-iw})/2$.
\end{proof}

We define, when $C<\infty$,
\begin{align}
  \label{defWstarstar}
  W^{\star\star}_C(u)
  &=W^{\star}_C(u)-W^\star_C(0)
  =W^{\star}_C(u)-J(W)\sum_{c\le C}\frac{\mu(c)}{c}
  \\&=\notag
  -\sum_{c\le C}\frac{\mu(c)}{c}W^\sharp(cu)
\end{align}
on recalling~\eqref{defWtildebis} and \eqref{defWstarC}. Note also that $W^{\star\star}_\infty=W^{\star}_\infty=W^\star$ by~\eqref{defWstar}.
We recall that $W^\sharp$ is defined at~\eqref{defWsharp}.

\begin{lem}
  \label{formulacheckW}
  When $\check{W}(s)\ll 1/(1+|s|)^3$, we have, when $u>0$,
  \begin{equation*}
    \hat{W}^{\star\star}_C(u)=\frac{-1}{2i\pi}
    \int_{-i\infty}^{i\infty}
    \frac{\check{W}(s)\zeta(s)}{u^{1-s}}
    \sum_{c\le C}\frac{\mu(c)}{c^{1+s}}
    ds.
  \end{equation*}
  When $C=\infty$, we replace $\sum_{c\ge
    1}{\mu(c)}/{c^{1+s}}$ by $1/\zeta(1+s)$.  As a consequence, when
  $C<\infty$ and for any real number $k<3/2$, we have $\hat{W}^{\star\star}_C(u)\ll
  (1+|u|)^{-1}(1+|u|/C)^{-k}$. Moreover, in the sense of distribution, we have
  $\hat{W}^\star_C(u)=J(W)\sum_{c\le C}\mu(c)/c\cdot\delta_{u=0}+ \hat{W}^{\star\star}_C(u)$ where
  $\delta_{u=0}$ is the Dirac mass at $u=0$.
\end{lem}

\begin{proof}
  The value $\hat{W}^{\star\star}_C(u)$ is the limit, as $Z$ goes to infinity, of
  \begin{equation*}
    2\int_0^Z W^{\star\star}_C(z) \cos(2\pi uz)dz.
  \end{equation*}
  We employ Lemma~\ref{formulaWstar} and reach the expression
  \begin{equation*}
    \frac{-1}{i\pi}
    \int_{-\ve-i\infty}^{-\ve+i\infty}
    \check{W}(-s)\zeta(1+s)\sum_{c\le C}\frac{\mu(c)}{c^{1-s}}
    \int_0^Z z^s\cos(2\pi zu)dz\, ds
  \end{equation*}
  which is also
  \begin{equation*}
    \frac{-1}{i\pi}
    \int_{-\ve-i\infty}^{-\ve+i\infty}
    \frac{\check{W}(-s)\zeta(1+s)}{u^{1+s}}\sum_{c\le C}\frac{\mu(c)}{c^{1-s}}
    \int_0^{uZ} z^s\cos(2\pi z)dz\, ds.
  \end{equation*}
  When $C=\infty$, we start with $\ve=0$ and shift the line of integration in $s$ just to the left-hand side of
  $\Re s=0$ but still within the zero-free region of $\zeta(1-s)$.
  Concerning the inner integral, we write
  \begin{multline*}
    \int_0^{uZ} z^s\cos(2\pi z)dz
    =
    \int_0^{1} z^s\cos(2\pi z)dz
    \\+
    (uZ)^s\frac{\sin 2\pi uZ}{2\pi}
    -\frac{s}{2\pi}
    \int_1^{uZ} z^{s-1}\sin(2\pi z)dz.
  \end{multline*}
  It is then enough to use the Lebesgue dominated convergence Theorem to send
  $Z$ to infinity (when $u>0$). We next appeal to Lemma~\ref{form} to get that
  \begin{align*}
    \int_0^{\infty} z^s\cos(2\pi z)dz
    &=(2\pi)^{-s-1}\Gamma(s+1)\cos(\pi (s+1)/2)
    \\&=-(2\pi)^{-1-s}\Gamma(1+s)\sin(\pi s/2)
    =
    \frac{1}{2}\frac{\zeta(-s)}{\zeta(1+s)}
  \end{align*}
  by using the functional equation of the Riemann zeta-function.
  This gives us
  \begin{equation}
    \label{eq:13}
    \hat{W}^{\star\star}_C(u)=\frac{-1}{2i\pi}
    \int_{-i\infty}^{i\infty}
    \frac{\check{W}(-s)\zeta(-s)}{u^{1+s}}\sum_{c\le C}\frac{\mu(c)}{c^{1-s}}
     ds
  \end{equation}
  The bound on $\hat{W}^{\star\star}_C(u)$ comes
  by separating the cases $|u|\le C$ and $|u|> C$ and in the latter case in
  shifting the line of integration to~$\Re s=k$ and using $|\zeta(-s)|\ll_\ve
  (1+|s|)^{k+1/2+\ve}$ (for any positive $\ve$) there.
\end{proof}

Let us mention the following consequence of Lemma~\ref{formulaWstar} together with
Mellin inversion formula.
\begin{lem}
  \label{MTWstar}
  The hypothesis on $W$ being as above, we have 
  \begin{equation*}
    \widecheck{W^{\star\star}_C}(s)
    =
    -\check{W}(s)\zeta(1-s)\sum_{c\le C}\mu(c)/c^{1+s}.
  \end{equation*}
  for $\Re s\in(0,3/2)$.
\end{lem}

\begin{lem}
  \label{ExplicitWhatstar}
  When $u>0$ and for $C\le \infty$, we have
  \begin{equation*}
    \hat{W}_C^{\star\star}(u)
    =
    \sum_{c\le C}\frac{\mu^2(c)}{c^2}\hat{W}(0)
    -
    \frac{1}{u}\sum_{n\ge1}\frac{\phi_C(n)}{n}W(n/u)
  \end{equation*}
  where $\phi_C(n)/n=\sum_{d|n,d\le C}\mu(d)/d$.
  In particular, this gives
  \begin{equation*}
    \hat{W}^{\star\star}(u)
    =
    \begin{cases}
      \frac{6}{\pi^2}\hat{W}(0)&\text{when $|u|\le 1/\upm$},\\
      \frac{6}{\pi^2}\hat{W}(0)-{W}(1/u)/u
      &\text{when $1/\upm<u\le 2/\upm$}.
    \end{cases}
  \end{equation*}
\end{lem}

\begin{proof}
  We only treat the case $C=\infty$. Lemma~\ref{formulacheckW}
  gives us
  \begin{equation*}
    \hat{W}^{\star\star}(u)=\frac{-1}{2i\pi}
    \int_{-i\infty}^{i\infty}
    \frac{\check{W}(s)\zeta(s)}{u^{1-s}}\sum_{c\le C}\frac{\mu(c)}{c^{1+s}}
     ds.
  \end{equation*}
  We shift the line of integration to $\Re s=2$ (since we move to the right,
  the contribution of the pole at $s=1$ is multiplied with a coefficient $-1$),
  use the development $\zeta(s)/\zeta(s+1)=\sum_{n\ge1}\phi(n)/n^{1+s}$ and
  the reverse Mellin transform to get
  \begin{equation*}
    \hat{W}^{\star\star}(u)
    =
    \frac{6}{\pi^2}\hat{W}(0)
    -
    \frac{1}{u}\sum_{n\ge1}\frac{\phi(n)}{n}W(n/u)
  \end{equation*}
  as expected.
\end{proof}

\begin{lem}
  \label{compFourierTransforms}
  We have $\hat W^{\star\star}(u)-\hat W^{\star\star}_C(u)\ll \log(|u|+2)/C$.
\end{lem}

\begin{proof}
  Indeed, by Lemma~\ref{ExplicitWhatstar}, we have $\hat
  W^{\star}(u)-\hat W^{\star\star}_C(u)\ll 1/C$ when $|u|<C/\upm$. When $u$ is larger, we use 
  \begin{equation*}
    \frac{\phi(n)}{n}-\frac{\phi_C(n)}{n}=\sum_{\substack{d|n,\\ d>C}}\frac{\mu(d)}{d}
    \ll 2^{\omega(n)}/C
  \end{equation*}
  where $\omega(n)$ is the number of prime factors of $n$.
  This implies that
  \begin{equation*}
    W^{\star\star}(u)-\hat W^{\star\star}_C(u)
    \ll
    C^{-1}+\frac{1}{u}\sum_{u\ll n\ll u}2^{\omega(n)}/C\ll \log(|u|+2)/C
  \end{equation*}
  as required.
\end{proof}

The size of $W^\star$ and $\hat W^\star$ is well controlled as shown in
the next lemma.
\begin{lem}
  \label{boundWstarC}
  Assume $W$ is at least C$\mkern1.5mu{}^3$. We have $W^{\star\star}_C(z)-J(W)\sum_{c\le C}\mu(c)/c\ll 1/(1+z^2)$. There
  exists $c_0>0$ (depending on $W$ only) such that, when $z\ge0$ and
  $\delta\in(0,1/2]$, we have $|W^{\star}(z+\delta)-W^{\star}(z)|\ll
  \exp-c_0\sqrt{-\log\delta}$ and, when $z\in(0,1]$, $W^{\star\prime}(z)\ll
  \gL(1/z)^{c_0}/z$. This shows in particular that $W^*$ is of bounded
  variations on $[0,1]$. Under the Riemann Hypothesis, we have
  $|W^{\star}(z)|\ll_\ve |z|^{\frac12-\ve}$ for any positive $\ve$. 
  
  When $z\le \lowm/C$, we have $W^{\star\star}_C(z)=0$.

  When $W$
  is four times differentiable, we have $|\hat W^\star(u)|\ll
  u^{-1}\gL(u)^{-c_0}$. Moreover $\hat
  W^\star(0)=\frac{6}{\pi^2}\int_0^\infty W(u)du$.
\end{lem}

\begin{proof}
  We split the proof is several stages.
  \par\noindent{\sl Bounding $W^{\star\star}_C$:}
  When $|z|\ge\lowm$, the first bound is a direct consequence of
  Lemma~\ref{boundWtilde}. When $|z|\le\lowm$, we write
  \begin{equation*}
    W^{\star\star}_C(z)=
    \sum_{\substack{c\le 1/|z|,\\ c\le C}}\frac{\mu(c)}{c}J(W)
    +\sum_{\substack{c> 1/|z|,\\ c\le C}}\frac{\mu(c)}{c}W(cz)
    =
    o(1)+\Ocal(1)
  \end{equation*}
  as required.
  \par\noindent{\sl Bounding the modulus of continuity of $W^{\star}$:}
  Appealing to Lemma~\ref{formulaWstar} with the change
  of variable $s\mapsto -s$, we next write
  \begin{equation*}
    W^\star(z+\delta)-W^\star(z)=
    \frac{1}{2i\pi}\int_{-i\infty}^{i\infty}\check{W}(-s)
    \frac{\zeta(1+s)}{\zeta(1-s)}sz^{s}\int_0^{\delta/z}(1+t)^{s-1}dtds.
  \end{equation*}
  Recalling that $\check{W}(-s)\ll 1/(1+|s|)^3$ and
  ${\zeta(1+it)}/{\zeta(1-it)}\ll(\log(2+|t|)^2$, this immediately gives us the bound
  $|W^\star(z+\delta)-W^\star(z)|\ll \delta/z$. This proves what we
  need (and more!) when $z\ge \sqrt{\delta}$. When $z$ is smaller, we
  proceed as in the proof of the Prime Number Theorem: when $t=\Im
  s\in[-T,T]$, we shift the line of integration to $\Re
  s=\sigma=c_1/\log T$ where $c_1>0$ is chosen so that
  $\zeta(\sigma-it)^{\pm1}\ll \log T$ when $|t|\le T$. The usual prime
  number theory gives us such a result, see e.g. \cite{Titchmarsh*51}. Skipping some
  classical steps, we reach the bound
  \begin{align*}
    W^\star(z+\delta)-W^\star(z)
    &\ll \frac{\delta}{z}z^{\frac{c_1}{\log T}}+\frac{(\log T)^2}{T}
    \\&\ll z^{\frac{c_1}{\log T}}+\frac{(\log T)^2}{T}
    \ll \delta^{\frac{c_1}{2\log T}}+\frac{(\log T)^2}{T}.
  \end{align*}
  We select $T=\exp(\sqrt{\log(1/\delta)}$. The reader will easily
  conclude from there.
  This is were the hypothesis $W$ C${}^3$ is needed.
  The bound for $W^{\star\prime}$ is obtained in the same manner.
  \par\noindent{\sl Some more upper bounds:}
  By Lemma~\ref{boundWtilde}, we have $\tilde W(z)=J(W)$ when $|z|\le
  \lowm$, hence $W_C^{\star\star}(z)=0$ when $|z|\le \lowm/C$.

  The bound for the Fourier transform follows by summation by parts.
  Concerning the value of the Fourier transform at~0, let $Z$ be a large parameter
  that goes to infinity. We write
  \begin{align*}
    2\int_0^Z \tilde W(z)dz
    &=
    2\int_0^1J(W)dz
    +2\int_1^Z\int_0^\infty\frac{B_1(z/u)}{z}(uW'(u)+W(u))dudz
    \\&=
    2J(W)
    +2\int_0^\infty\int_{1/u}^{Z/u}\frac{B_1(z)dz}{z}(uW'(u)+W(u))du
    \\&=
    2J(W)
    -2\int_0^\infty\biggl(
    \frac{B_1(Z/u)}{Z/u}\left(\frac{-Z}{u^2}\right)
    - \frac{B_1(1/u)}{1/u}\left(\frac{-1}{u^2}\right)
    \biggr)uW(u)du
    \\&=
    2J(W)
    +2\int_0^\infty
    B_1(Z/u)W(u)du
    - 2\int_0^\infty B_1(1/u)W(u)du
    \\&=
    2\int_0^\infty
    B_1(Z/u)W(u)du
    +\int_0^\infty W(u)du
  \end{align*}
  and the integral depending on $Z$ goes to 0 as $Z$ goes to infinity by
  Lebesgue's Lemma. This shows that $\hat{\tilde{W}}(0)=(1/2)\hat{W}(0)$. We
  next employ~\eqref{defWstarC} to deduce that
  \begin{equation*}
    \hat{W}^*_C(u)=\sum_{c\le C}\frac{\mu(c)}{c^2}\hat{\tilde{W}}(u/c)
  \end{equation*}
  hence the value at~$u=0$, whether $C<\infty$ or not.
\end{proof}
\section{Numerical aspects related to the smoothing kernel and its transforms}
\label{ExplicitSmoothing}
It is interesting to produce some numerical datas, so as to explore
our several transforms.
\subsection{An explicit family of smoothing kernels}
Let $\1_{[-1,1]}$ be the characteristic function of the interval $[-1,1]$. We
are interested in explicit formulae for the $m$-th convolution-power
$\1_{[-1,1]}^{(*m)}$, where $m$ is a positive integer. This function is even
with support within $[-m,m]$, and of class C$^{m-1}$. We readily check that
\begin{equation}
  \label{deb}
  \1_{[-1,1]}^{(*2)}(t)
  =
  \begin{cases}
    \displaystyle
    2-|t|&\text{when $ |t|\le 2$},\\
    0&\text{when $2\le |t|$}.
  \end{cases}
\end{equation}
Some more sweat brings the next formula:
\begin{equation*}
  \1_{[-1,1]}^{(*3)}(t)
  =
  \begin{cases}
    \displaystyle
    3-t^2&\text{when $|t|\le 1$},\\
    \displaystyle
    (3-|t|)^2/2&\text{when $1\le |t|\le 3$},\\
    0&\text{when $3\le |t|$}.
  \end{cases}
\end{equation*}
The general formula is given in \cite{Renyi*70} and reads
\begin{equation*}
  \1_{[-1,1]}^{(*m)}(t)
  =
  \begin{cases}
    \displaystyle
  \sum_{j=0}^{\lfloor(m+|t|)/2\rfloor}
  \frac{(-1)^j}{(m-1)!}\binom{m}{j}(m+|t|-2j)^{m-1}
    &\text{when $0\le |t|\le m$},\\
    0&\text{when $m<|t|$}.
  \end{cases}
\end{equation*}
Guessing this expression is not obvious, but verifying it by recursion is only
a matter of routine. The Fourier transform of $\1_{[-1,1]}$ is $\sin(2\pi
u)/(\pi u)$, so the one of $\1_{[-1,1]}^{(*m)}$ is $\sin(2\pi
u)^m/(\pi u)^m$. 
Since we will use the case $m=5$, it is worth giving its explicit expression:
\begin{equation}
  \label{eq:22}
  \1_{[-1,1]}^{(*5)}(t)
  =
  \begin{cases}
    \frac{115-30t^2+3t^4}{12}&\text{when $|t|\le 1$},
    \\
    \frac{55+10|t|-30t^2+10|t|^3-t^4}{6}&\text{when $1\le |t|\le 3$},
    \\
    \frac{625-500t+150t^2-20|t|^3+t^4}{24}&\text{when $3\le |t|\le 5$}
    \\0&\text{when $5\le |t|$}.
  \end{cases}
\end{equation}


Formula \eqref{explicitTransform-1} is handy for explicit computations. We  introduce
\begin{equation*}
  \mathfrak{p}_m(t)=\frac{4m}{2^m}\1^{*m}_{[-1,1]}(4mt-3m)
\end{equation*}
for some integer $m\ge5$. Its support lies inside $[1/2,1]$. We find that
\begin{equation*}
  \hat{\mathfrak{p}}_m(u)=e(3u/4)\left(\frac{\sin(\pi u/(2m))}{\pi u/(2m)}\right)^m.
\end{equation*}
Notice that  $\int_{0}^\infty \mathfrak{p}_m(t)dt=\hat{\mathfrak{p}}_m(0)=1$.
We then select
\begin{equation*}
  W(m;t)=\mathfrak{p}_m(1/t)/t.
\end{equation*}
For such a choice, we readily get
\begin{equation*}
  \Transform{W}(m;z)
  =
  2\sum_{n\ge1}\frac{\phi(n)}{n}
  \cos(3\pi nz/2)\left(\frac{\sin(\pi nz/(2m))}{\pi nz/(2m)}\right)^m.
\end{equation*}
When we truncate this series at the integer $N$, the error is bounded above by
\begin{equation}
  \label{error}
  2\left(\frac{2m}{\pi z}\right)^m\frac{1}{(m-1)N^{m-1}}.
\end{equation}
We then use the following Sage script (see \cite{sagemath}):
\begin{verbatim}
def Witself(t, m = 5):
   if abs(t) > 2 or abs(t) < 1:
      return(0)
   res = 0
   z = m*(4/t-3)
   coef = 2*m/factorial(m-1)/2^m
   asign = 1
   for j in range(0, floor(float((m + abs(z))/2)) + 1):
      res += asign*binomial(m, j)*(m + abs(z) -2*j)^(m-1)
      asign = -asign
   return(res*coef/t)

plot(lambda t:Witself(t, 5), (1, 2))
\end{verbatim}
\subsection{A specific kernel}
\label{ExplicitKernel}
In this section, we specify $m=5$.
\subsubsection*{\sl On $W(5;t)$:}
Here is a plot of our function.
\begin{figure}[!h]
\includegraphics[scale=0.45]{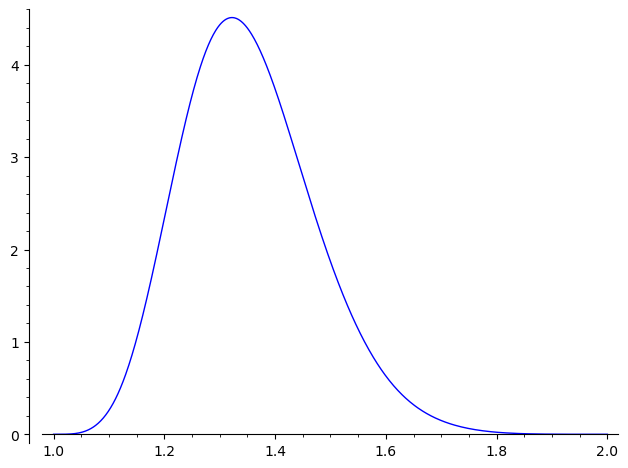}
\caption{$W(5;t)$}
\end{figure}
\FloatBarrier
\noindent
The command \texttt{integral\_numerical(lambda
  t:Witself(t,5), (1,2))}
gives us
\begin{equation*}
  \frac{6}{\pi^2}\int_0^\infty W(5;t)dt=0.816\cdots
\end{equation*}
\subsubsection*{\sl On $\Transform{W}(5;t)$:}

We get the following plot on $[0.0001, 3]$:

\begin{figure}[!h]
\subfigure[]{\includegraphics[scale=0.375]{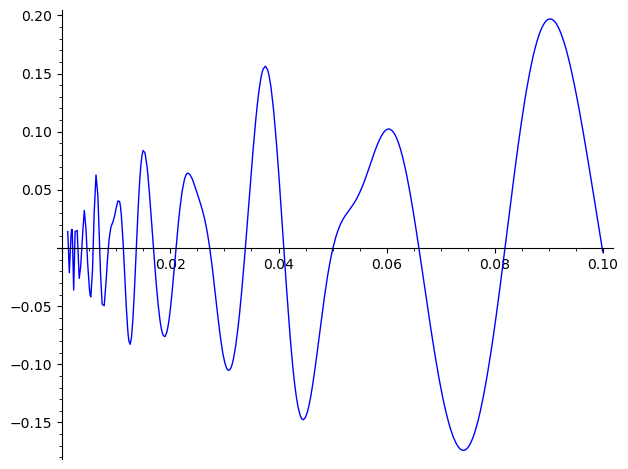}}
\subfigure[]{\includegraphics[scale=0.375]{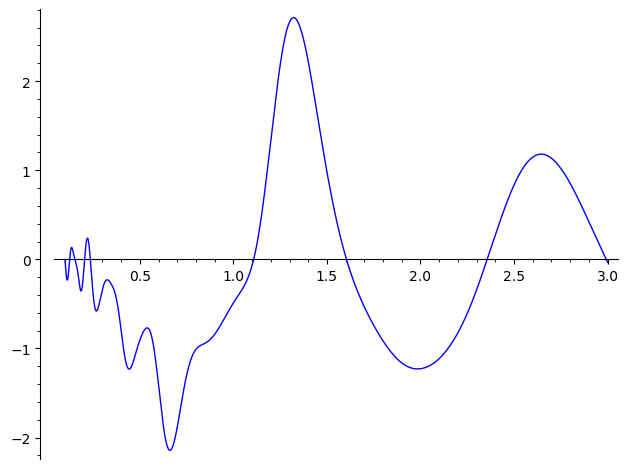}}
\caption{$\Transform{W}(5;z)$ for $0.0001\le z\le 0.1$ and for $0.1\le z\le 3$}
\end{figure}
\FloatBarrier


And here is a plot of $\hat{\Transform{W}}(5;t)$.
\begin{figure}[!h]
\includegraphics[scale=0.45]{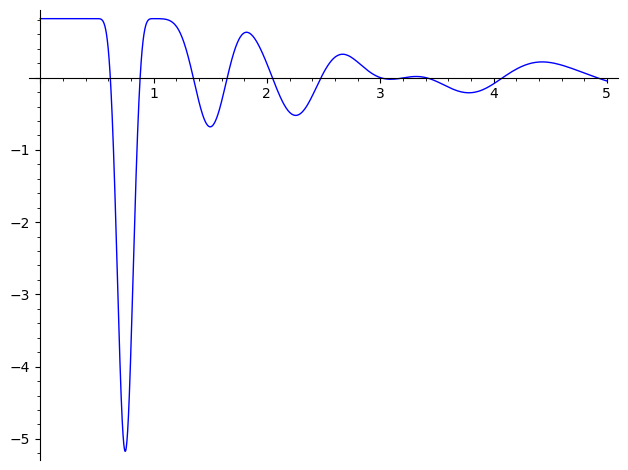}
\caption{$\hat{\Transform{W}}(5;t)$}
\end{figure}
It is worth noticing that
$\hat{\Transform{W}}(5;1)=\hat{\Transform{W}}(5;0)$. After $u=1$, we
indeed find that $\hat{\Transform{W}}(5;u)<\hat{\Transform{W}}(5;0)$.
\FloatBarrier

\section{A general formula, first step in the proof of Theorem~\protect\ref{Precise}}
\label{gene}

In analytic number theory, when we want to detect an equality, the quantity we
really study is of the shape
$\sum_{\substack{m,n}}\varphi_m\overline{\psi_m}\delta_{m=n}$ and that what we
use in an approximation of the $\delta$-symbol.  This is not only a tautology,
it also imposes a framework which decides of what are the ``trivial''
estimates and of what can be expected or not. It also splits the problem in
two parts: a combinatorial part, where one uses the fact $m$ and $n$ are
integers, possibly in certain subsequences, and an analytical part where the
quantities arising are to be estimated. There is of course an interplay
between both parts and a ``good'' decomposition is a decomposition that leads
to quantities that we know how to estimate.  It is difficult to give a precise
historical date, but the contributions of M. Jutila in \cite{Jutila*92} (see also
\cite{Harcos*03} and \cite[Theorem 2]{Jutila*07}) and of H. Iwaniec in
\cite{Duke-Friedlander-Iwaniec*93} (see also \cite[]{Duke-Iwaniec*92} and
\cite[Chapter 20]{Iwaniec-Kowalski*04}, in particular Proposition 20.16
therein) seem to be prominent. One can say rapidly that in some sense,
Iwaniec's way is to analyze the large sieve quantity to extract a diagonal
contribution, under some hypotheses, while Jutila's way is to start from the
diagonal contribution and to modify the circle to keep only the rationals one
knows how to handle, with a possible weight.

The present study is centered on the quantity
\begin{equation}
  \label{defmyS}
  \myS(Q,W)=
  \sum_{q}\frac{W(q/Q)}{q}\sum_{a\mode  q}
  |S(\varphi, a/q)|^2.
\end{equation}
Moebius inversion readily yields
\begin{equation*}
  \myS(Q,W)=
  \sum_{d}\sum_{d|q}\frac{\mu(q/d)W(q/Q)}{q}\sum_{a\mod d}
  |S(\varphi, a/d)|^2.
\end{equation*}
We expand the square, shuffle the terms around and get
\begin{equation}
  \myS(Q,W)=
  \sum_{\substack{m,n}}\varphi_m\overline{\varphi_n}
  \Delta(m-n)
\end{equation}
where we have use the notation (on setting $cd=q$)
\begin{equation}
  \label{defDelta}
  \Delta(v)=\sum_{\substack{c,d,\\d|v}}\frac{\mu(c)W(cd/Q)}{c}.
\end{equation}
Here is the decomposition of the $\Delta$-symbol we use.
\begin{lem}[Iwaniec's decomposition]
  \label{dec}
  Let $C,E, H\ge1$ be parameters that satisfy $E\le \min(\lowm Q, \upm Q/C)$.
  We have
  \begin{equation*}
    \Delta(v)
    =
    U(v)+U^\sharp(v)
    + L_0(v)
    +L(v)+L^\sharp(v)
  \end{equation*}
  where $L_0(v)$ is the diagonal contribution
    \begin{equation*}
    L_0(v)=\sum_{\substack{c\le C,\\ d\ge1}}\frac{\mu(c)W(cd/Q)}{c}\1_{v=0},
  \end{equation*}
  and $U(v)$ and $U^\sharp(v)$ are the ``direct divisor'' part:
  \begin{equation*}
    \left\{
      \begin{aligned}
        U(v)&=-\sum_{e\le E}\sum_{\substack{c\le C,\\ f\ge1}}
        \frac{\mu(c)W(cef/Q)}{cef} \sum_{a\mode e}e(av/e),
        \\
        U^\sharp(v)&=\sum_{e> E}\sum_{\substack{c> C,\\ f\ge1}}
        \frac{\mu(c)W(cef/Q)}{cef}\sum_{a\mode e}e(av/e),
      \end{aligned}
    \right.
  \end{equation*}
  while $L(v)$ and $L^\sharp(v)$ are the ``complementary divisor'' part:
  \begin{equation*}
    \left\{
      \begin{aligned}
        L(v)&=\sum_{h\le H}\sum_{\substack{h|g,\\ c\le C}}\frac{\mu(c)}{gc}
        \sum_{a\mode h}W(cv/(gQ))e(av/h),
        \\
        L^\sharp(v)&=\sum_{h> H}\sum_{\substack{h|g,\\ c\le C}}\frac{\mu(c)}{gc}
        \sum_{a\mode h}W(cv/(gQ))e(av/h).
      \end{aligned}
    \right.
  \end{equation*}
\end{lem}

\begin{proof}
  We start by splitting the range for the variable $c$:
 \begin{align*}
  \Delta(v)
  &=
  \sum_{\substack{c\le C,\\d|v}}\frac{\mu(c)W(cd/Q)}{c} 
  +
  \sum_{\substack{c>C,\\d|v}}\frac{\mu(c)W(cd/Q)}{c} 
  \\&=
  L(v)+U(v)
\end{align*}
say. When $v=0$, the term $L(v)$ restricts to $L_0(v)$. Otherwise, we
switch to the complementary divisor by setting $gd=|v|$ (and $g\ge1$ since
$v\neq 0$). We detect the divisibility condition by using additive characters:
\begin{align*}
  L(v)
  &=
  \sum_{\substack{c\le C,\\g|v}}\frac{\mu(c)W(c|v|/(gQ))}{c} 
  \\&=
  \sum_{\substack{c\le C,\\g\ge1}}\frac{1}{g}\sum_{b\mod
    g}\frac{\mu(c)W(c|v|/(gQ))}{c}
  e(bv/g)
  \\&=
  \sum_{\substack{c\le C,\\g\ge1}}\frac{1}{g}\sum_{h|g}\sum_{b\mode 
    h}\frac{\mu(c)W(c|v|/(gQ))}{c}
  e(bv/h)
\end{align*}
which amounts to
\begin{equation*}
  L(v)
  =
  \sum_{h\ge 1}\sum_{\substack{c\le C,\\ h|g}}\frac{\mu(c)}{gc}
  \sum_{b\mode  h}W(cv/(gQ)) e(bv/h).
\end{equation*}
Note that we do not need the condition $v\neq0$ since $W(cv/(gQ))=0$ when $v=0$. 
We then simply split the summation over $h$ according to whether $h\le H$ or
not, getting the two quantities $L(v)$ and $L^\sharp(v)$.

Concerning $U(v)$ we again detect the divisibility condition by using
additive characters. This gives us
\begin{equation*}
  U(v)
  =\sum_{\substack{c>C,\\ d\ge1}}\frac{\mu(c)W(cd/Q)}{cd}
  \sum_{e|d}
  \sum_{a\mode e}e(av/e).
\end{equation*}
Note that $cd/Q\le \upm$. We set $d=ef$ and thus $e\le \upm Q/C$. We continue by
splitting the range for~$e$:
\begin{multline*}
  U(v)
  =\sum_{e\le E}\sum_{\substack{c>C,\\ f\ge1}}\frac{\mu(c)W(cef/Q)}{cef}
  \sum_{a\mode e}e(av/e)
  \\+
  \sum_{ e>E}\sum_{\substack{c>C,\\ f\ge1}}\frac{\mu(c)W(cef/Q)}{cef}
  \sum_{a\mode e}e(av/e).
\end{multline*}
We recognize $U^\sharp(v)$ in the last quantity. The first one needs a
transformation.
We note that
\begin{eqnarray*}
  \sum_{\substack{c>C,\\ f\ge 1}}\frac{\mu(c)W(cef/Q)}{cef}
   &=&
   \sum_{\substack{c\ge1,\\ f\ge 1}}\frac{\mu(c)W(cef/Q)}{cef}-
   \sum_{\substack{c\le C,\\ f\ge 1}}\frac{\mu(c)W(cef/Q)}{cef}
   \\&=&
   \sum_{j\ge1}\sum_{cf=j}\frac{\mu(c)W(je/Q)}{je} -
   \sum_{\substack{c\le C,\\ f\ge 1}}\frac{\mu(c)W(cef/Q)}{cef}
   \\&=&
   \frac{W(e/Q)}{e}-
   \sum_{\substack{c\le C,\\ f\ge 1}}\frac{\mu(c)W(cef/Q)}{cef}
\end{eqnarray*}
and the first term vanishes because of the assumption $E\le Q$.
\end{proof}

The diagonal term is easily handled.

\begin{lem}
  \begin{equation*}
    \sum_{m,n}\varphi_m\overline{\varphi_n}
    L_0(m-n)
    =
    \left(\sum_{q}\frac{\phi(q)W(q/Q)}{q}+\Ocal(Q C^{-1})\right)\|\varphi\|_2^2.
  \end{equation*}
\end{lem}

\begin{proof}
  The
contribution is
\begin{equation*}
  \sum_{c\le C,d}\frac{\mu(c)W(cd/Q)}{c}\|\varphi\|_2^2
  =\sum_{q}\sum_{c|q,c\le C}\frac{\mu(c)}{c}W(q/Q)\|\varphi\|_2^2.
\end{equation*}
Since
\begin{equation*}
  \sum_{c> C,d}\frac{\mu(c)W(cd/Q)}{c}\ll \sum_{c> C}\frac{Q}{c^2}\ll Q/C
\end{equation*}
we get that this diagonal term has value:
\begin{equation*}
  \left(\sum_{q}\frac{\phi(q)W(q/Q)}{q}+\Ocal(Q C^{-1})\right)\|\varphi\|_2^2
\end{equation*}
as announced.
\end{proof}

The large sieve inequality yields an efficient bound for the
contribution of~$U^\sharp(m-n)$.

\begin{lem}
  We have 
  \begin{equation*}
    \sum_{m,n}\varphi_m\overline{\varphi_n}
  U^\sharp(m-n)
  \ll
  \sum_{m}|\varphi_m|^2 (NE^{-1}+QC^{-1})\log Q.
  \end{equation*}
\end{lem}

\begin{proof} 
  We use the bound (where $c$ and $e$ are fixed)
  \begin{equation*}
    \sum_f\frac{W(cef/Q)}{f}\ll
    \sum_{Q/(ce)\le f\le 3Q/(ce)}\frac{1}{Q/(ce)}\ll 1
  \end{equation*}
  to get:
  \begin{align*}
    \sum_{m,n}\varphi_m\overline{\varphi_n}
    U^\sharp(m-n)
    &\ll (\log Q)
    \sum_{E< e\le 3Q/C}e^{-1}
    \sum_{a\mode  e}|S(\varphi, a/e)|^2
    \\&\ll
    \sum_{m}|\varphi_m|^2 (NE^{-1}+QC^{-1})\log Q.
  \end{align*}
\end{proof}

The contribution of $L^\sharp(m-n)$ is somewhat more difficult to handle but also
relies on the large sieve inequality. We shall most of the time employ the
next lemma with a set $I$ reduces to one element. It is only in the final
applications that it is better to use the summation over some~$i\in I$.

\begin{lem}\label{intromax}
  Let $w$ be an even and $C^1$ function that vanishes when the variable is
  larger than~1. We further assume that $w$ is piecewise $C^2$. 
  Let $I$ be a finite set. We have
  \begin{equation*}
    \sum_{i\in I}
    \Bigl|\sum_{m,n}\psi_{m,i}\overline{\psi_{n,i}}w(\alpha (m-n))\Bigr|
    \le 3 \|w''\|_1 (N\alpha+1)\max_{u<v\le u+3/\alpha}
    \sum_{i\in I}
    \Bigl|\sum_{u<m\le v}\psi_{m,i}\Bigr|^2.
  \end{equation*}
\end{lem}

\begin{proof}
  The problem is twofold: localizing the variables $m$ and $n$ and separating
  these two variables. The first problem is met by a subdivision argument: we
  cover the interval $[1,N]$ by at most $N\alpha+1$ disjoint intervals
  $(a,a+\alpha^{-1}]=I_a$ of length $\alpha^{-1}$ and localize $n$ within such an
  interval. As a result we can assume that $m$ lies in
  $[a-\alpha^{-1},a+3\alpha^{-1}]=J_a$.
  We handle the separation of variables by a summation by parts and the formula
  \begin{align*}
    w(\alpha(n-m))
    &=
    -\alpha\int_{m-\alpha^{-1}}^{n}w'(\alpha(t-m))dt
    \\&=
    \alpha^2\int_{m-\alpha^{-1}}^{n}
    \int_{m}^{t+\alpha^{-1}}w''(\alpha(t-s))dsdt
  \end{align*}
  from which we infer that $\sum_{n\in I_a,m\in J_a}\psi_{m,i}\overline{\psi_{n,i}}
  w(\alpha (m-n))$ equals
  \begin{equation}
    \label{identite}
    \alpha^2 \int_{a-2\alpha^{-1}}^{a+\alpha}\int_{a-\alpha^{-1}}^{t+\alpha^{-1}}
    \sum_{\substack{s\le m\le  t-\alpha^{-1},\\  t\le n\le a+\alpha^{-1}}}
    \psi_{m,i}\overline{\psi_{n,i}}
    w''(\alpha(t-s))dsdt.
  \end{equation}
  We find that $t-3\alpha^{-1}\le a-\alpha^{-1}\le s\le m\le t-\alpha^{-1}$
  and that $t\le n\le a+\alpha^{-1}\le t+3\alpha^{-1}$, hence the inner sum
  over $m$ and $n$ is bounded above (after introducing the summation over
  $i$, by
  \begin{equation*}
    \max_{u<v\le u+3/\alpha}
    \sum_{i\in I}
    \Bigl|\sum_{u<m\le v}\psi_{m,i}\Bigr|^2.
  \end{equation*}
  A change of variables readily shows that
  \begin{equation*}
    \alpha^2
    \int_{a-2\alpha^{-1}}^{a+\alpha}\int_{a-\alpha^{-1}}^{t+\alpha^{-1}}
     |w''(\alpha(t-s))|dsdt\le 3 \|w''\|_1,
  \end{equation*}
  clearing out any uniformity problem in applications.
\end{proof}

\begin{lem}
  We have 
  \begin{equation*}
    \sum_{m,n}\varphi_m\overline{\varphi_n}
  L^\sharp(m-n)
  \ll
  \bigl( NH^{-1} + NC Q^{-1}\bigr)
  \|\varphi\|_2^2\log^4(QN).
  \end{equation*}
\end{lem}

\begin{proof}
We have to control
\begin{equation}
  \sum_{m,n}\psi_m\overline{\psi_n}W(\alpha (m-n))
\end{equation}
where $\psi_m=\varphi_m e(ma/h)$ and $\alpha=c/(gQ)\neq0$. Note that the truncation
in $c$ ensures that $|\alpha|$ is small; this truncation has been introduced 
for this very purpose. 
Practically, we appeal to Lemma~\ref{intromax} and get
\begin{multline*}
  S=\sum_{c\le C}\frac{\mu(c)}{c}
  \sum_{h>H}\sum_{h|g}g^{-1}
  \sum_{a\mode h}
  \sum_{m,n}\varphi_me(ma/h)\overline{\varphi_ne(na/h)}W(c|m-n|/(gQ))
  \\\ll
  \sum_{\substack{c\le C,h>H,\\ h|g\le cN/Q}}
  \frac{1}{gc}\Bigl(\frac{Nc}{gQ}+1\Bigr)
  \sum_{a\mode h}
  \max_{u<v\le u+9gQ/c}\Bigl|\sum_{u<m\le v}\varphi_me(am/h)\Bigr|^2
\end{multline*}
The condition $v>u$ is automatically satisfied.
We continue with $c$ fixed by localizing $h$ and using
$k=g/h$. Lemma~\ref{LSM} gives us:
\begin{multline*}
  S\ll\mkern-20mu
  \sum_{\substack{1\le k\le Q/H,\\\log H\le \ell\le \log\frac{cN}{Q}}}
  \sum_{e^{\ell-1}<h \le e^\ell}
  \Bigl(\frac{N/Q}{k^2e^{2\ell}}+\frac{1}{ke^\ell c}\Bigr)
  \sum_{\substack{a\mode  h}}
  \max_{v\le u+\frac{9ke^\ell Q}{c}}
  \Bigl|\sum_{u<m\le v}\mkern-8mu\varphi_me\Bigl(\frac{am}{h}\Bigr)\Bigr|^2
  \\\ll
  \sum_{1\le k\le Q/H}
  \sum_{\log H\le \ell\le \log(cN/Q)}
  \Bigl(\frac{N/Q}{k^2e^{2\ell}}+\frac{1}{ke^\ell c}\Bigr)
  \Bigl(\min(N, k e^\ell Q c^{-1})+\ell e^{2\ell}\Bigr)
  \|\varphi\|_2^2
  \\\ll
  \sum_{1\le k\le Q/H}
  \sum_{\log H\le \ell\le \log(cN/Q)}
  \Bigl(\frac{N}{ kce^\ell}+\frac{N\ell}{Qk^2}+\frac{\ell e^\ell}{kc}\Bigr)
  \|\varphi\|_2^2
  \\\ll
  \sum_{\log H\le \ell\le \log(cN/Q)}
  \Bigl(\frac{N}{ ce^\ell}+\frac{N\ell}{Q}+\frac{\ell e^\ell}{c}\Bigr)
  \|\varphi\|_2^2\log(Q/H)
  \\\ll
  \Bigl(N H^{-1} c^{-1}+N Q^{-1}\Bigr)
  \|\varphi\|_2^2\log^3(QN)
\end{multline*}
so this contribution is at most (on summing over $c$), up to a
multiplicative constant:
\begin{equation}
  \label{eq:1}
  Q \bigl( N(HQ)^{-1} + NC Q^{-2}\bigr)
  \|\varphi\|_2^2\log^4(QN).
\end{equation}

\end{proof}

This approximation provided by Lemma~\ref{evalWflat} together with the
large sieve inequality leads to the following formula (recall the
definition \eqref{defWsharp} of $W^\sharp$): \smallskip

\noindent\fbox{\vbox{%
\begin{multline}
  \label{formula2}
  \sum_{q}\frac{W(q/Q)}{qQ\ooWo}\sum_{a\mode q}
  |S(\varphi, a/q)|^2
  =
  \sum_m|\varphi_m|^2
  \\
  -\frac{J(W)}{Q\ooWo}
  \sum_{\substack{c\le C,\\ e\le E}}\frac{\mu(c)}{ec}
  \sum_{a\mode  e}|S(\varphi, a/e)|^2
  \\+
  \sum_{\substack{c\le C,\\ h\le H}}\frac{\mu(c)}{chQ\ooWo}
  \sum_{a\mode h}
  \sum_{m,n}\varphi_m\overline{\varphi_n}W^\sharp(c|m-n|/(hQ))e((n-m)a/h)
  \\+\Ocal\biggl( \Bigl(\frac{N}{ EQ}+ \frac{N}{HQ} + \frac{1}{C}+ \frac{NC}{ Q^{2}}\Bigr)
  \|\varphi\|_2^2\log^5(QN)\biggr)
\end{multline}
}}
\smallskip

\noindent
The first main term comes from $L_0$, the second one from $U$ and the
third one from $L$.

\section{Proof of Theorem~\ref{Precise}}
\subsection{From $W^\sharp$ to $\tilde W$: cancellation of the two main terms}
\label{Cancel}

We introduce $\tilde{W}$ by appealing to~\eqref{defWtildebis}.
The choice $E=H$ ensures that, in \eqref{formula2}, the second main term is
canceled out by the contribution of the factor linked with the $J(W)/h$
above, getting
\smallskip

\noindent\fbox{\vbox{%
\begin{multline}
  \label{formula3}
  \sum_{q}\frac{W(q/Q)}{qQ\ooWo}\sum_{a\mode q}
  |S(\varphi, a/q)|^2
  =
  \|\varphi\|_2^2
  \\
  -
  \sum_{\substack{c\le C,\\ h\le H}}\frac{\mu(c)}{chQ\ooWo}
  \sum_{a\mode h}
  \sum_{m,n}\varphi_m\overline{\varphi_n}\tilde W(c|m-n|/(hQ))e((n-m)a/h)
  \\+\Ocal\biggl( \Bigl(\frac{N}{HQ} + \frac{1}{C}+ \frac{NC}{ Q^{2}}\Bigr)
  \|\varphi\|_2^2\log^5(QN)\biggr).
\end{multline}
}}

\smallskip
The same cancellation of the main term is what presides to the introduction of
$\Delta_c(u)$ in \cite[section 20.5]{Iwaniec-Kowalski*04}, see the proof of
Lemma~20.17 therein.

\subsection{Sharpening the error term in its $H$-dependence}

One of the error term in Eq.~\eqref{formula3} is
$\Ocal(\frac{N}{HQ}\|\varphi\|^2_2\log^5(QN))$ and we want to (and need to!)
remove the $\log^5(QN)$. We have to consider
\begin{equation}
  \label{eq:12}
  \Sigma(H_1,H_2)=
  \sum_{\substack{H_1<h\le H_2}}\frac{1}{h}
  \sum_{a\mode h}
  \sum_{m,n}\varphi_m\overline{\varphi_n} W^\star_C(|m-n|/(hQ))e((n-m)a/h).
\end{equation}
We somehow go backwards and use $W^{\star\star}_C$ from~\eqref{defWstarstar} to write
\begin{equation*}
  \Sigma(H_1,H_2)
  =
  W^\star_C(0)\sum_{\substack{H_1<h\le H_2}}\frac{1}{h}
  \sum_{a\mode h}
  \Bigl|S\Bigl(\varphi,\frac{a}{h}\Bigr)\Bigr|^2
  +
  \Sigma'(H_1,H_2)
\end{equation*}
with
\begin{equation*}
  \Sigma'(H_1,H_2)=
  \int_{-\infty}^{\infty}
  \sum_{\substack{H_1<h\le H_2}}\frac{1}{h}
  \sum_{a\mode h}
  \Bigl|S\Bigl(\varphi,\frac{a}{h}+\frac{u}{Qh}\Bigr)\Bigr|^2
  \hat{W}^{\star\star}_C(u)du.
\end{equation*}
The large sieve inequality readily yields (since $W^\star_C(0)\ll 1$)
\begin{equation*}
  \Sigma(H_1,H_2)-\Sigma'(H_1,H_2)
  \ll \Bigl(\frac{N}{H_1}+H_2\Bigr)\|\varphi\|^2_2.
\end{equation*}
The treatment of $\Sigma'(H_1,H_2)$ is somewhat more difficult.
When $|u/Q|\le 1/2$, by combining a summation by parts together with the
large sieve inequality, we find that
\begin{align*}
  \sum_{\substack{H_1<h\le H_2}}\frac{1}{h}
  \sum_{a\mode h}
  \Bigl|S\Bigl(\varphi,\frac{an}{h}+\frac{un}{hQ}\Bigr)\Bigr|^2
  \le
  \|\varphi\|_2^2\Bigl(
  \frac{N}{H_1}+8H_2
  \Bigr)
\end{align*}
since the points $(\frac{a}{h}+\frac{u}{hQ})_{a,h}$ are $\tfrac12
H_2^{-2}$-well spaced.  When $|u/Q|\ge 1/2$, we use the large sieve inequality
for every~$h$. In this case the shift by $u/(hQ)$ is constant and the points
are $h^{-1}$-well-spaced, giving
  \begin{equation*}
    \sum_{\substack{H_1<h\le H_2}}\frac{1}{h}
    \sum_{a\mode h}
    \Bigl|S\Bigl(\varphi,\frac{an}{h}+\frac{un}{hQ}\Bigr)\Bigr|^2
    \ll
    \|\varphi\|_2^2\Bigl(
    N\log\frac{2H_2}{H_1}+H_2
    \Bigr).
  \end{equation*}
  As a consequence
  \begin{equation*}
    \Sigma'(H_1,H_2)/\|\varphi\|_2^2\ll
    \frac{N}{H_1}+8H_2
    +
    \frac{C}{Q}\Bigl(
    N\log\frac{2H_2}{H_1}+H_2
    \Bigr)
  \end{equation*}
  on using the bound $|\hat W^{\star\star}_C(u)|\ll
  C/(1+|u|^2)$ from Lemma~\ref{formulacheckW} when
  $|u|\ge Q/2$.
  This implies that
  \begin{equation*}
    \Sigma'(H_0,\sqrt{N})
    \ll \frac{N}{H_0}+\sqrt{N}+\frac{CN}{Q}\log N.
  \end{equation*}
  We can use formula~\eqref{formula3} with $H=\sqrt{N}$ and shorten the
  summation by the process above. On renaming $H_0=H$, we have reached:
  \smallskip

\noindent\fbox{\vbox{%
\begin{multline}
  \label{formula3b}
  \sum_{q}\frac{W(q/Q)}{qQ\ooWo}\sum_{a\mode q}
  |S(\varphi, a/q)|^2
  =
  \|\varphi\|_2^2
  \\
  -
  \sum_{\substack{ h\le H}}\frac{1}{hQ\ooWo}
  \sum_{a\mode h}
  \sum_{m,n}\varphi_m\overline{\varphi_n}W^\star_C(|m-n|/(hQ))e((n-m)a/h)
  \\+\Ocal\biggl( \Bigl(\frac{N}{HQ} + \frac{\log^5(QN)}{C}+ \frac{NC\log^5(QN)}{ Q^{2}}\Bigr)
  \|\varphi\|_2^2\biggr).
\end{multline}
}}
\smallskip

\noindent
The effect of the previous treatment is neat: the $\log$-factor attached to $N/(HQ)$ has
disappeared while the rest of the remainder term is still of the same order of
magnitude.

\subsection{Direct extension of the $c$-variable}

We handle the sum over~$c$
essentially trivially. The contribution from the diagonal term $m=n$ is
bounded above by $\sum_{c> C}\frac{\mu(c)}{c}H\|\varphi\|_2^2 /Q$. When
$|m-n|\le hQ/c$, we bound $\tilde{W}(c|m-n|/(hQ))$ by $\Ocal(1)$, getting a
contribution bounded above, up to a multiplicative constant, by
\begin{multline*}
  \sum_{\substack{h\le H,\\ c> C}}\frac{1}{chQ\ooWo}
  \sum_{a\mode h}
  \sum_{|m- n|\le hQ/c}|\varphi_m\varphi_n|
  \\
  \ll
  \sum_{\substack{h\le H,\\ c> C}}\frac{1}{cQ}
  \sum_{m}|\varphi_m|^2\frac{hQ}{c}
  \ll \|\varphi\|_2^2H^2/C.
\end{multline*}
We
use $\tilde W(z)\ll 1/(1+z^2)$ when
$|m-n|> hQ/c$, getting a
contribution bounded above, up to a multiplicative constant, by
\begin{multline*}
  \sum_{\substack{h\le H,\\ c> C}}\frac{1}{chQ\ooWo}
  \sum_{a\mode h}
  \sum_{|m- n|>hQ/c}\frac{|\varphi_m\varphi_n|}{1+ c^2(m-n)^2/(h^2Q^2)}
  \\\ll
  \sum_{\substack{h\le H,\\ c> C}}\frac{Q^2h^3}{c^3hQ}
  \sum_{m}|\varphi_m|^2\frac{c}{hQ}
  \ll 
  \|\varphi\|_2^2 H^2/C.
\end{multline*}
We thus get, for any $C'\ge C$:
\begin{multline*}
  \sum_{q}\frac{W(q/Q)}{qQ\ooWo}\sum_{a\mode q}
  |S(\varphi, a/q)|^2
  =
  \|\varphi\|_2^2
  \\
  -
  \sum_{\substack{ h\le H}}\frac{1}{hQ\ooWo}
  \sum_{a\mode h}
  \sum_{m,n}\varphi_m\overline{\varphi_n}W^\star_{C'}(|m-n|/(hQ))e((n-m)a/h)
  \\+\Ocal\biggl( \Bigl(\frac{N}{HQ} 
  +\frac{H^2}{C}
  + \frac{\log^5(QN)}{C}+ \frac{NC\log^5(QN)}{ Q^{2}}\Bigr)
  \|\varphi\|_2^2\biggr).
\end{multline*}
The optimal
choice $C=QH/N^{1/2}$ (provided that $H\le N^{1/4}$; Indeed we recall that
Lemma~\ref{dec} asks for $E\le \min( Q,\upm Q/C)$ and that we have chosen $E=H$)
may be too large.
Instead we select 
\begin{equation}
  \label{defC}
  C= \min\Bigl(\frac{QH}{\sqrt{N}},\frac{\upm Q}{H}, C'\Bigr)
  = \min\Bigl(\frac{QH}{\sqrt{N}}, C'\Bigr)
\end{equation}
and get

\noindent\fbox{\vbox{%
\begin{multline}
  \label{formula4}
  \sum_{q}\frac{W(q/Q)}{qQ\ooWo}\sum_{a\mode q}
  |S(\varphi, a/q)|^2
  =
  \|\varphi\|_2^2
  \\
  -
  \sum_{\substack{ h\le H}}\frac{1}{hQ\ooWo}
  \sum_{a\mode h}
  \sum_{m,n}\varphi_m\overline{\varphi_n}W^\star_{C'}(|m-n|/(hQ))e((n-m)a/h)
  \\+\Ocal\biggl( \Bigl(\frac{N}{HQ} 
  + \frac{H^2+\log^5(QN)}{C}\Bigr)
  \|\varphi\|_2^2\biggr).
\end{multline}
}}
\smallskip

\noindent
We may reformulate this equality by using the Fourier transform of
$W^\star$:
\begin{multline*}
  \sum_{q}\frac{W(q/Q)}{qQ\ooWo}\sum_{a\mode q}
  |S(\varphi, a/q)|^2
  =
  \|\varphi\|_2^2
  \\
  -
  \sum_{\substack{h\le H}}\frac{1}{hQ\ooWo}
  \sum_{a\mode h}
  \int_{-\infty}^{\infty}
  \hat
  W^\star_{C'}(u)\Bigl|S\Bigl(\varphi,
  \frac{an}{h}+\frac{un}{hQ}\Bigr)\Bigr|^2
  du
  \\+\Ocal\biggl( \Bigl(\frac{N}{HQ} + \frac{H^2+\log^5(QN)}{ C}\Bigr)
  \|\varphi\|_2^2\biggr).
\end{multline*}
Later, to prove~\eqref{formula11b}, it will be better to restrict the range of integration (note that the Fourier
transform has two parts: a Dirac mass and a regular part; only the
regular part is concerned, as the Dirac mass is concentrated at $u=0$). We use
the large sieve inequality with $u$ and $h$ fixed to infer that

\smallskip
\noindent\fbox{\vbox{%
\begin{multline}
  \label{formula5}
  \sum_{q}\frac{W(q/Q)}{qQ\ooWo}\sum_{a\mode q}
  |S(\varphi, a/q)|^2
  =
  \|\varphi\|_2^2
  \\
  -
  \sum_{\substack{h\le H}}\frac{1}{hQ\ooWo}
  \sum_{a\mode h}
  \int_{-U}^{U}
  \hat
  W^\star_{C'}(u)\Bigl|S\Bigl(\varphi,
  \frac{a}{h}+\frac{u}{hQ}\Bigr)\Bigr|^2
  du
  \\+\Ocal\biggl( \Bigl(\frac{N}{HQ} 
  + \frac{NC'\log H}{UQ}
  + \frac{H^2+\log^5(QN)}{ C}\Bigr)
  \|\varphi\|_2^2\biggr).
\end{multline}
}}
\smallskip
\noindent
We can however proceed in a different fashion: majorize $|\hat
W^*(u)|$ when $|u|\ge U$ by $\Ocal(1/U)$, uniformly in $C$, and use
$\int_{-\infty}^{\infty}|S(\alpha+u/(hQ)|^2du=hQ\|\varphi\|_2^2$ by Parseval.
This leads to

\smallskip
\noindent\fbox{\vbox{%
\begin{multline}
  \label{formula5b}
  \sum_{q}\frac{W(q/Q)}{qQ\ooWo}\sum_{a\mode q}
  |S(\varphi, a/q)|^2
  =
  \|\varphi\|_2^2
  \\
  -
  \sum_{\substack{h\le H}}\frac{1}{hQ\ooWo}
  \sum_{a\mode h}
  \int_{-U}^{U}
  \hat
  W^\star_{C'}(u)\Bigl|S\Bigl(\varphi,
  \frac{a}{h}+\frac{u}{hQ}\Bigr)\Bigr|^2
  du
  \\+\Ocal\biggl( \Bigl(\frac{N}{HQ} 
  + \frac{H^2}{U}
  + \frac{H^2+\log^5(QN)}{ C}\Bigr)
  \|\varphi\|_2^2\biggr).
\end{multline}
}}
\smallskip
\noindent
The difference from $\hat W^\star_{C'}$ to $\hat W^{\star\star}_{C'}$
is $J(W)\sum_{c\le C'}\mu(c)/c\cdot\delta_{u=0}$ by
Lemma~\ref{formulacheckW}. On using that $J(W)\ll1$, that $\sum_{c\le
  C'}\mu(c)/c\ll1$ and the large sieve inequality, we get a
contribution which is $\ll NH^{-1}\|\varphi\|^2_2$, thus incorporable in
the already existing error term.  We have obtained:

 \smallskip
\noindent\fbox{\vbox{%
\begin{multline}
  \label{formula51}
  \sum_{q}\frac{W(q/Q)}{qQ\ooWo}\sum_{a\mode q}
  |S(\varphi, a/q)|^2
  =
  \|\varphi\|_2^2
  \\
  -
  \sum_{\substack{h\le H}}\frac{1}{hQ\ooWo}
  \sum_{a\mode h}
  \int_{-U}^{U}
  \hat
  W^{\star\star}_{C'}(u)\Bigl|S\Bigl(\varphi,
  \frac{a}{h}+\frac{u}{hQ}\Bigr)\Bigr|^2
  du
  \\+\Ocal\biggl( \Bigl(\frac{N}{HQ} 
  + \min\biggl(\frac{NC'\log H}{UQ},\frac{H^2}{U}\biggr)
  + \frac{H^2+\log^5(QN)}{ C}\Bigr)
  \|\varphi\|_2^2\biggr).
\end{multline}
}}

\smallskip
\noindent
We can send $U$ to infinity and Theorem~\ref{Precise} follows by
keeping $U=\infty$ and sending also $C'$ to infinity.

\section{A case of large sieve equality. Proof of Theorem~\ref{Vic}}

We prove a first result that is suited for some
applications.
\begin{thm}
  \label{betterVic}
  When $\frac12\le H\le \sqrt{N}/(\log N)^5$ and $\log Q\ll \log N$, we have
  \begin{multline*}
  \sum_{q}\frac{W(q/Q)}{qQ}\sum_{a\mode q}
  |S(\varphi, a/q)|^2
  =
  \bigl(\ooWo
  + \Ocal(N(HQ)^{-1})\bigr)
  \|\varphi\|_2^2
  \\
  +\Ocal\Bigl(
  \sum_{\substack{h\le H}}\frac{N+hQ}{h^2Q^2}
  \max_{u<v<u+2hQ}
  \sum_{a\mode h}
  \Bigl|\sum_{u<n\le v}\varphi_n e(na/h)\Bigr|^2
  \Bigr).
\end{multline*}
\end{thm}

\begin{proof}
  Ideally, we would simply combine Theorem~\ref{Precise} (but we
  convert back $\hat W^\star$ in $W^\star$ as in~\eqref{formula4})
  together with Lemma~\ref{intromax} applied to $W^\star$, the set $I$
  being $\{a\mode h\}$. The function $W^\star$ is however not regular
  enough, and we have to revert to $W^\star_C$ and more precisely to
  Eq.~\eqref{formula4}. We select $C=QH/\sqrt{N}$. When $z\le 1/C$, we
  have $(W^\star_C)''=0$ while Lemma~\ref{formulaWstar} with $\ve=0$
  implies that $(W^\star_C)''(z)\ll 1$ in general. The theorem follows readily.
\end{proof}

\begin{proof}[Proof of Theorem~\ref{Vic}]
  We employ Theorem~\ref{betterVic} and simplify the remainder term by
  appealing to
  \begin{align*}
    \sum_{a\mode h}
    \Bigl|\sum_{u<n\le v}\varphi_n e(na/h)\Bigr|^2
    &\le
    \sum_{a\mod h}
    \Bigl|\sum_{u<n\le v}\varphi_n e(na/h)\Bigr|^2
    \\&\le
    h\sum_{c\mod h}
    \biggl|\sum_{\substack{u<n\le v,\\ n\equiv c[h]}}\varphi_n\biggr|^2.
  \end{align*}
  Such an extension of the variable $a$ may look a weak step, but since this
  theorem is aimed at sequences oscillating highly in small arithmetic
  progressions, the loss is not noticeable (at least in the examples I could
  think of).
\end{proof}

\section{A refinement for primes}
\label{Refinement}
When the sequence $\varphi$ is supported on integers prime to every
integer $h\le H$, we may refine Theorem~\ref{Precise} further, thanks
to the next improved large sieve inequality. This is \cite[Theorem
5.3]{Ramare*06}. See also \cite[Corollary 1.5]{Ramare*22}.
\begin{lem}\label{EIR}
  If $(\varphi_n)_{n\le N}$ is such that $\varphi_n$ vanishes as soon as $n$ has a prime
  factor 
  less than $\sqrt{N}$, then
  \begin{equation*}
    \sum_{q\le Q_0}\sum_{a\mode q}
    \bigl|S(\varphi,a/q)\bigr|^2
    \le 7\,\frac{N\log Q_0}{\log N}\|\varphi\|_2^2
  \end{equation*}
  for any $Q_0\le \sqrt{N}$ and provided $N\ge100$.
\end{lem}
This lemma enables us to improve Theorem~\ref{Precise} into the next result.
\begin{thm}
  \label{PrecisePrimes}
  When $1/2\le H\le \sqrt{N}/(\log N)^5$, $Q\le 10 N$ and $\varphi_n$
  vanishes when $n$ has a prime factor below $\sqrt{N}$, we have
  \begin{multline*}
    \sum_{q}\frac{W(q/Q)}{qQ}\sum_{a\mode q} |S(\varphi, a/q)|^2 = \biggl(\ooWo +
    \Ocal\biggl(\frac{N\log 3H}{QH\log N}\biggr)\biggr) \sum_m|\varphi_m|^2
    \\
    - \sum_{\substack{h\le H}}\frac{1}{h} \sum_{a\mode h}
    \int_{-\infty}^{\infty} \hat W^\star(u)\Bigl|S\Bigl(\varphi,
    \frac{a}{h}+\frac{u}{hQ}\Bigr)\Bigr|^2 du.
\end{multline*}
\end{thm}

\begin{proof}
  We start from Theorem~\ref{Precise}, but with say $H'$ rather than $H$
  and now shorten the sum over $h$. To do so, we write
  \begin{multline*}
    \sum_{\substack{h\sim H_1}}\frac{1}{hQ} \sum_{a\mode h}
    \int_{-\infty}^{\infty} \hat W^\star(u)
    \Bigl|S\Bigl(\varphi,
    \frac{a}{h}+\frac{u}{hQ}\Bigr)\Bigr|^2du
    \\=
    \frac{1}{H_1Q} 
    \int_{-\infty}^{\infty} \max_{h\sim H_1}\biggl|\hat W^\star\biggl(\frac{hv}{H_1}\biggr)\biggr|
    \sum_{\substack{h\sim H_1}}\sum_{a\mode h}\Bigl|S\Bigl(\varphi,
    \frac{a}{h}+\frac{v}{H_1Q}\Bigr)\Bigr|^2dv.
  \end{multline*}
  Lemma~\ref{EIR} tells us that this quantity is $\ll \frac{N\log
    3H_1}{QH_1\log N}\|\varphi\|_2^2$ from which, after noticing the
  bound for $\hat W^\star$ from Lemma~\ref{boundWstarC}, the theorem follows readily.
\end{proof}

\part{Operator Decomposition of the Large Sieve}

\section{A local geometrical space}
\label{infinity}

We consider $X_h=\Z{h}\times[0,1]$, equipped with the product of the
probability measures.
We denote by $L^2_*(X_h)$ the space of functions from $L^2(X_h)$ whose Fourier
transform with respect to the first variable is supported by
$(\Z{h})^*\times[0,1]$, i.e. functions~$f$~such that
\begin{equation*}
  \forall y\in[0,1],\forall d\in\Z{h}\ /\gcd(d,h)>1,\quad
  \sum_{b\mod h}f(b,y)e(-db/h)=0.
\end{equation*}
It is maybe simpler to say that this is the space generated by the functions
$(c,y)\mapsto e(ac/q)f(y)$ for all $f\in L^2([0,1])$ and (this is where a
restriction occurs)~$a$ prime to~$q$.
We reproduce rapidly the theory developed in \cite[Chapter~4]{Ramare*06}. Let
$k|h$ be two moduli. We consider
\begin{equation}
  \label{defLkh}\arraycolsep=1.4pt
  \begin{array}{rcl}
    L^k_h:L^2(X_k)&\rightarrow&
    \displaystyle L^2(X_h)\\[0.5em]
    \displaystyle F&\mapsto&
    \begin{array}[t]{rcl}
      L^k_h(F):\Z{h}\times[0,1]&\rightarrow&\mathbb{C}\\
      (b,y)&\mapsto&\displaystyle F(\sigma_k(b),y)
    \end{array}
  \end{array}
\end{equation}
and correspondingly
\begin{equation}
  \label{defJkh}\arraycolsep=1.4pt
  \begin{array}{rcl}
    J^h_k:L^2(X_h)&\rightarrow&
    \displaystyle L^2(X_k)\\[0.5em]
    \displaystyle F&\mapsto&
    \begin{array}[t]{rcl}
      J^h_k(F):\Z{k}\times[0,1]&\rightarrow&\mathbb{C}\\
      (b,y)&\mapsto&\displaystyle 
      \frac{1}{h/k}\sum_{\substack{c\mod h,\\ c\equiv b[k]}}F(\sigma_h(c),y).
    \end{array}
  \end{array}
\end{equation}
We finally define
\begin{equation}
  \label{defUktildevahtilde}
  U_{\tilde h\goes\tilde k}=L^k_hJ^h_k,\quad
  U_{\tilde h\goes k}=\sum_{d|k}\mu(k/d)U_{\tilde h\goes\tilde d}.
\end{equation}
Here is the structure theorem we need.\footnote{These results are easily
  proved. Details may be found in \cite[Chapter~4]{Ramare*06}, though with no
  $y$-component. This component is inert here, so the proofs carry through
  mutatis mutandis.}
\begin{thm}
  The maps $L^k_h$ and $J^h_k$ are adjoined one to the other. The
  collection $(U_{\tilde h\goes\tilde k})_{k|h}$ is a family of commuting
  orthogonal projectors. Furthermore
  \begin{equation*}
    U_{\tilde h\goes\tilde k}=\sum_{d|k}U_{\tilde h\goes d}
  \end{equation*}
  while, for any two divisors $k_1$ and $k_2$ of $h$, we have
  $
    U_{\tilde h\goes k_1}U_{\tilde h\goes k_2}=\delta_{k_1=k_2}U_{\tilde h\goes k_2}$.
  We have $L^2_*(X_h)=U_{\tilde h\goes h}L^2(X_h)$.
\end{thm}

\subsubsection*{An explicit expression}

At the heart of this matter are the Gauss sums
\begin{equation}
  \label{defGauss}
  \tau_h(\chi,\cdot)=\sum_{b\mod h}\chi(b)e(b\cdot/h).
\end{equation}

\begin{thm}
  \label{ExplicitUhtildeh}
  For any $h\ge1$, any class $b$ modulo $c$, any real number $y$ and any
  function $F\in L^2(X_h)$, the orthonormal projection $U_{\tilde h\goes h}$ on $L^2_*(X_h)$
  has the following explicit form:
  \begin{equation*}
    U_{\tilde h\goes h}F(b,y)=\frac{1}{h}\sum_{c\mod h}c_h(b-c)F(c,y).
  \end{equation*}
  Given a hilbertian orthonormal basis $(f_k)_k$ of $L^2([0,1])$,
  the family $(\Base_{h,\chi}\otimes f_k)_{\chi,k}$ where
  $\Base_{h,\chi}=\tau_h(\chi,\cdot)/\sqrt{\phi(h)}$ and  $\chi$ ranges the 
  Dirichlet characters modulo~$h$ is a hilbertian orthonormal basis of  $L^2_*(X_h)$.
\end{thm}

\begin{proof}
  We first check that
  \begin{align*}
    \sum_{b\mod h}c_h(b-c)e(bd/h)
    &=
    \sum_{a\mode h}e(-ac/h)\sum_{b\mod h}e(b(a+d)/h)
    \\&=
    \begin{cases}
      h e(dc/h)&\text{when $(d,h)=1$},\\
      0&\text{else}.
    \end{cases}
  \end{align*}
  and since $(b\mapsto e(bd/h))_{d\mod h}$ generates the whole space of
  functions over $X_h$, this proves our first assertion.
  The introduction of the Dirichlet character may be arbitrary, but in fact
  $(\tau_h(\chi,\cdot))_\chi$ is the full set of eigenfunctions of 
  $
  f\mapsto
  \sum_{c\mod h}c_h(b-c)f(c)/h
  $
  that are associated to a non-zero eigenvalue. We simply have
  \begin{equation}
    \label{deftauhchib}
    \forall b\in\Z{h},\quad\tau_h(\chi,b)=\frac{1}{h}
    \sum_{c\mod h}c_h(b-c)\tau(\chi,c).
  \end{equation}
  Note finally that
  \begin{align*}
    \frac{1}{h}\sum_{c\mod h}
    \tau_h(\chi_1,c)\overline{\tau_h(\chi_2,c)}
    &=
    \sum_{a,b\mod h}\chi_1(a)\overline{\chi_2(b)}
    \frac{1}{h}\sum_{c\mod h}e\Bigl(\frac{c(a-b)}{h}\Bigr)
    \\&=
    \1_{\chi_1=\chi_2}\varphi(h)
  \end{align*}
  as required.
\end{proof}

\section{Analysis of a class of difference operators}
\label{ACDO}

We treat here the analysis of the intervening family of operators in
an abstracted setting.
Let $V$ be a function satisfying the following assumptions:
\begin{itemize}
\item[($R_1$) $\bullet$] $V$ is a continuous real-valued even function of bounded
  variations and integrable over $\mathbb {R}$.
\item[($R_2$) $\bullet$] $V(0)=0$.
\item[($R_3$) $\bullet$] There exist $B\ge \|V\|_\infty$,  $c\in(0,1]$ and
  $A>0$ such that, for every $\delta\in(0,1)$ and $x\in[0,1-\delta]$, we have
  $|V(x+\delta)-V(x)|\le B \exp(-c\sqrt{-\log \min(1,A\delta)})$. 
\end{itemize}
Recall that we defined
\begin{equation}
  \tag{\ref{defV0}}
  \Vscr_0:\quad G\in L^2([0,1])\mapsto \biggl(y\mapsto\int_0^1 G(y')V(y-y')dy'\biggr)
\end{equation}
It is classical theory that $\Vscr_0$ is a compact Hilbert-Schmidt
operator, see for instance \cite[Theorem
7.7]{Gohberg-Goldberg-Krupnik*00}. Let $(\lambda_\ell,G_\ell)_\ell$ be
a complete orthonormal system of eigenvalues / eigenfunctions, ordered
with non-increasing $|\lambda_\ell|$.  The Fredholm equation
$\lambda G(y')=\int_0^1 K(y',y)G(y)dy$ has been intensively
studied. It is not the purpose of this paper to introduce to this
theory, a task for which it is better to read the complete and
classical book~\cite{Gohberg-Krein*69} by I. Gohberg, I. C. \&
M.G. Kre{\u i}n, or the more modern \cite{Gohberg-Goldberg-Krupnik*00}
by I. Gohberg, S. Goldberg \& N. Krupnik.  Kernel of type $V(y'-y)$
are often called \emph{difference kernel}, and lead to operators that
are distinct from convolution operators as the integration and
definition interval is \emph{not} the whole real line. The book
\cite{Sakhnovich*15} by L. Sakhnovich is dedicated to the operators
built from such kernels. The book \cite{Cochran*72} by J. Cochran
contains also many useful informations.

\subsection{L${}^2$-norm}

We readily find that
\begin{equation}
  \int_0^1\int_0^1 |V(y'-y)|^2 dy
  =
  \int_{-1}^1|V(z)|^2(1-|z|)dz.
\end{equation}
Hence
\begin{equation*}
  \sum_{\ell\ge 1}|\lambda_\ell|^2=\int_{-1}^1|V(z)|^2(1-|z|)dz
  \le 2 \int_{0}^1|V(z)|^2dz.
\end{equation*}
As a consequence, and enumerating the eigenvalues in such a way that
$|\lambda_\ell|$ is non-increasing, we find that
\begin{equation}
  \label{localbound1}
  |\lambda_\ell|\le \sqrt{2\int_{0}^1|V(z)|^2dz}/\sqrt{\ell}.
\end{equation}
Theorem~\ref{nuclear} will enable us to replace~$\sqrt{\ell}$
by~$\ell$, but it uses the above bound.

\subsection{Properties of the eigenvectors}
%

The eigenvectors of~$\Vscr_0$ attached to non-zero eigenvalues are
classically shown to be continuous. Since the $L^1$-norm is not more than the
$L^{2}$-norm squared here, we have $\|G\|_1\le 1$. Each of them thus satisfies
\begin{equation}
  \label{propeigenvectors}
  |\lambda|\|G\|_\infty\le 2\int_0^1 |V(z)|(1-z)dz.
\end{equation}
Furthermore, we find that
\begin{equation}
  \label{regulationcontinuity}
  |\lambda||G(y+\delta)-G(y)|\le \|G\|_1 \omega(V,\delta)
  = \|G\|_1 \max_{-1\le z\le 1}|V(z+\delta)-V(z)|.
\end{equation}
These functions are also of bounded variation. Indeed, with obvious notation,
we find that
\begin{align*}
  |\lambda|\sum_{1\le i\le n}|G(y_{i+1})-G(y_i)|
  &\le \int_0^1 |G(y)|\sum_{1\le i\le n}
    \bigl|V(y_{i+1}-y)-V(y_i-y)\bigr|dy
  \\&\le \int_{-1}^1|V'(y)|dy\int_0^1 |G(y)|dy\le \int_{-1}^1|V'(y)|dy
\end{align*}
since $\|G\|_1\le 1$.

\subsection{Nuclearity}

A consequence of a theorem of Fredholm
from \cite{Fredholm*03} is that, when $y\mapsto V(y)$ is H\"older of
exponent~$\alpha$, then the eigenvalues verify
$\sum_{\ell\ge1}|\lambda_\ell|^{p}<\infty$ for every
$p>2/(1+2\alpha)$. This proof is reproduced in the book~\cite[Chapter IV, Theorem
  8.2]{Gohberg-Goldberg-Krupnik*00} by I. Gohberg, S. Gohberg \& N. Krupnik. This is too strong a condition for us if
we are to avoid the Riemann Hypothesis (in which case $\alpha=1/2+\ve$ would
be accessible). D. Swann in \cite{Swann*71} considered the effect of bounded
variation on a general kernel, but his theorem asks again for too strong
hypotheses since the function $(y',y)\mapsto V(y'-y)$ is a priori not of
bounded variation. However, each function $y\mapsto V(y'-y)$ is uniformly
of bounded variation (i.e. its total variation is, as function of $y'$
integrable; in our case, it is even bounded), a case that is mentioned
(with more generality) in the paragraph preceding \cite[Theorem 3]{Swann*71}
and more formally in \cite[Theorem 16.2]{Cochran*72} in the monograph
of J. Cochran. We follow this approach.

In this subsection, we use 
\begin{equation}
  \label{deflog-}
  \log^-t=\log\min(1,t);
\end{equation}
We consider the coefficients of the Carleman determinant, see \cite[Chapter 4, (3)]{Cochran*72}, for $\nu\ge2$:
\begin{equation}
  \label{defdnu}
  d_\nu=\frac{(-1)^\nu}{\nu!}\int_0^1\mkern-10mu\cdots\mkern-8mu\int_0^1
  \begin{vmatrix}
0 & V(y_1-y_2) &\mkern-7mu \cdots\mkern-6mu & V(y_1-y_\nu) \\ 
V(y_2-y_1) & 0 &\mkern-7mu \cdots\mkern-6mu & V(y_2-y_\nu) \\ 
\vdots & \vdots &  & \vdots \\ 
V(y_\nu-y_1) & V(y_\nu-y_2) &\mkern-7mu \cdots\mkern-6mu &  0 \\
\end{vmatrix}dy_1dy_2\cdots dy_\nu.
\end{equation}
As $V(y-y)=0$, this is also the Fredholm determinant, see
\cite[Chapter VI, (1.5)]{Gohberg-Goldberg-Krupnik*00}.
The above determinant, say $K(y_1,\cdots,y_\nu)$, can be rewritten as
\begin{equation*}
  \begin{vmatrix}
0 & V(y_1-y_2)-V(y_1-y_1) &\mkern-7mu \cdots\mkern-6mu &
V(y_1-y_\nu)-V(y_1-y_{\nu-1}) \\  
V(y_2-y_1) & V(y_2-y_2)-V(y_2-y_1)  &\mkern-7mu \cdots\mkern-6mu &
V(y_2-y_\nu)-V(y_2-y_{\nu-1}) 
\\\vdots & \vdots &  & \vdots \\ 
V(y_\nu-y_1)&V(y_2-y_\nu)-V(y_2-y_{\nu-1})   
 &\mkern-7mu \cdots\mkern-6mu &  V(y_\nu-y_\nu)-V(y_\nu-y_{\nu-1})  \\
\end{vmatrix}.
\end{equation*}
We use the symmetry of the integral and now assume that $0\le y_1< y_2<\cdots
<y_\nu\le1$ (when an equality occurs between these variables, the determinant
vanishes). We define $\delta_i=y_{i+1}-y_i$ so that $\sum_{1\le i\le
  \nu-1}\delta_i\le 1$. We divide the second column by
$\sqrt{B}\exp(-(c/2)\sqrt{-\log^-(A\delta_1)})$, the third one by
$\sqrt{B}\exp(-(c/2)\sqrt{-\log^-(A\delta_2)})$ and so on, getting a factor
\begin{equation*}
  B^{(\nu-1)/2}\prod_{1\le i\le \nu-1}\exp(-(c/2)\sqrt{-\log^-(A\delta_i)})
\end{equation*}
in front of our determinant. We first note the following lemma.
\begin{lem}
  We have $\sum_{1\le i\le n}\sqrt{-\log^-(A \delta_i)}\ge n\sqrt{\log n}$
  when the $\delta_i$'s are positive real numbers such that
  $\sum_{1\le i\le n}\delta_i\le 1$.
\end{lem}
\begin{proof}
  Given an $n$-tuple $(\delta_1,\cdots,\delta_n)$, we note that the $n$-tuple
  obtained by replacing each $\delta_i$ by $\min(A^{-1},\delta_i)$ satisfies the same
  constraint with an equal sum of $\sqrt{-\log^-(A\cdot)}$. In order to find
  the minimum required, we may thus restrict our attention to variables that
  verify $\delta_i\le 1/A$.  Set $x_i=(-\log^- (A\delta_i))^{1/4}$. This
  variable ranges possibly $(0,\infty)$. The condition on $\delta_i$ now reads
  $\sum_{1\le i\le n}e^{-x_i^4/A}=\delta$ for some $\delta\in(0,1]$, while we
  seek to minimize $\sum_{1\le i\le n}x_i^2$ and we forget the condition
  $e^{-x_i^4/A}\le 1/A$. We use the Lagrange method and
  consider
  \begin{equation*}
    Y(x_1,\dots,x_n,\lambda) =
    \sum_{1\le i\le n}x_i^2-\lambda\bigr(\sum_{1\le i\le n}e^{-x_i^4/A}-\delta\bigl).
  \end{equation*}
  Its critical points, obtained by equating all the partial derivatives to 0, satisfy:
  \begin{equation*}
    \left\{
      \begin{aligned}
        &\forall i\le n,\quad 2x_i+4A^{-1}x_i^3\lambda e^{-x_i^4/A}=0,\\
        &\sum_{1\le i\le n}e^{-x_i^4/A}=\delta.
      \end{aligned}
    \right.
  \end{equation*}
  This implies that\footnote{Any choice $x_i=0$ means that
    $\delta_i=1$, which implies that any other $\delta_j$ vanishes,
    leading to the maximum being $\infty$ when $n\ge2$.}
  $\lambda/A=-2e^{-x_i^4/A}/x_i^2$. The function
  $y\mapsto 2e^{-y^4/A}/y^2$ is decreasing, from which we conclude
  that all $x_i$'s are equal, which in turn implies that all
  $\delta_i$'s are equal, and equal to $\delta/n$. The choice
  $\delta=1$ is also optimal.
\end{proof}

Next we use Hadamard's inequality (as in all such proofs!) together
with the previous lemma (and $(\nu-1)\ge \nu/2$) and get
\begin{multline*}
  \frac{|K(y_1,\cdots,y_\nu)|}{B^{\nu/2}e^{-(c/4)\nu\sqrt{\log(\nu-1)}}}
  \le \prod_{i}\Bigl(\sum_{j}|a_{i,j}|^2\Bigr)^{1/2}
  \\\le 
  \prod_{i}
  \Bigl(
  \|V\|_\infty\bigl|V(y_i-y_1)-V(y_i-y_{i})\bigr|
  +B\sum_{2\le j\le \nu}\bigl|V(y_i-y_{j-1})-V(y_i-y_{j})\bigr|\Bigr)^{1/2}
  \\\le 
  \Bigl(2c\int_{-1}^1\bigl|V'(y)\bigr|dy\Bigr)^{\nu/2}
\end{multline*}
since $B\ge \|V\|_\infty$.
As a consequence, we find that the Carleman determinant
\begin{equation}
  \label{defmyDcal}
  \Dcal(V,z)=1+\sum_{\nu\ge2}d_\nu z^\nu
  =\prod_{\ell\ge1}(1-\lambda_\ell z)e^{\lambda_\ell z}
\end{equation}
satisfies, with $M=\sqrt{2B\int_{0}^1\bigl|V'(y)\bigr|dy}$,
\begin{align*}
  |\Dcal(V,z)|
  &\le 1+\sum_{\nu\ge2} \frac{M^\nu|z|^\nu e^{-\frac{c}{4}\nu\sqrt{\log(\nu-1)}}}{\nu!}
  \\&
  \le
  (M|z|)^{N+1}e+e^{M|z|e^{-\frac{c}{4}\sqrt{\log N}}}
  \le
  H^{N+1}e+e^{He^{-\frac{c}{4}\sqrt{\log N}}}
\end{align*}
with $H=M|z|\ge1$ and for any real valued parameter $N\ge1$ that we
may choose.  When $H\le e^2$, we use the upper bound
$|\Dcal(V,z)|\le e^{e^2}$. When $\log H\ge 2$, we select
\begin{equation}
  \label{choiceN}
  N=He^{-\frac{c}{4}\sqrt{\log (H+1)}}/\log H.
\end{equation}
When $\log H\ge 2$,
we check that (recall that we have assumed that $c\le1$)
\begin{align*}
  \log N
  &=\log H-\frac{c}{4}\sqrt{\log H}-\log\log H
  \\&=\log (H+1) \Bigl(\frac{\log H}{\log(H+1)}
      -\frac{1}{4\sqrt{\log (H+1)}}-\frac{\log\log H}{\log (H+1)}\Bigr)
  \\&
      \ge \log (H+1) \Bigl(\frac{\log 2}{\log 3}-\frac{1}{4\sqrt{\log 3}}-\frac{1}{e}\Bigr)
      \ge\frac{\log (H+1)}{49}.
\end{align*}
We thus find that, in this case, we have
\begin{align*}
  |\Dcal(V,z)|
  &\le 
  H e^{He^{-\frac{c}{4}\sqrt{\log (H+1)}}}e+e^{He^{-\frac{c}{28}\sqrt{\log (H+1)}}}
  \\&\le 
  6H e^{He^{-\frac{c}{28}\sqrt{\log (H+1)}}}.
\end{align*}
Next, $He^{-\frac{c}{28}\sqrt{\log (H+1)}}+\log H$ is certainly not
more than $He^{-\frac{c}{30}\sqrt{\log (H+1)}}$ provided $H$ be larger
than some constant depending on~$c$. So, in general, we find that
$He^{-\frac{c}{28}\sqrt{\log (H+1)}}+\log H\le
He^{-\frac{c}{30}\sqrt{\log (H+1)}}+c''$,
where $c''$ is a constant depending solely on $c$. We have proved that
\begin{equation*}
  |\Dcal(V,z)|\le 6e^{c''}e^{He^{-\frac{c}{30}\sqrt{\log (H+1)}}}
\end{equation*}
when $H\ge e^2$. The minimum of
$6e^{c''}e^{He^{-\frac{c}{30}\sqrt{\log (H+1)}}}$ when $H$ ranges
$[0,e^2]$ is some positive constant, say $c'''$, depending only on $c$
(we have introduced $\log(H+1)$ rather tha $\log H$ earlier for this
very purpose). As a consequence, we have, for any $H\ge0$,
\begin{equation*}
  |\Dcal(V,z)|\le 
  e^{e^2}\frac {6e^{c''}e^{He^{-\frac{c}{30}\sqrt{\log (H+1)}}}}{\min(1,c''')}.
\end{equation*}
Here is the lemma we have proved.
\begin{lem}
  \label{majD}
  There exists a positive constant $c'=c'(c)$ such that we have 
  \begin{equation*}
    |\Dcal(V,z)|
    \le c' e^{M|z|e^{-\frac{c}{30}\sqrt{\log(M|z|+1)}}}
  \end{equation*}
  with $M=\sqrt{2B\int_{0}^1\bigl|V'(y)\bigr|dy}$.
\end{lem}

We continue with the following general lemma.
\begin{lem}
  \label{Jensenplus}
  Let $f$ be an entire function of finite order and such that $f(0)=1$ and let
  $(\rho_\ell)$ be an enumeration of its zeroes with non-decreasing
  $|\rho_\ell|$. Let $g$ be a $C^2$-function over $(0,\infty)$. Assume that,
  as $t$ goes to infinity,
  \begin{equation*}
    \frac{1}{2\pi}\int_0^{2\pi}\log|f(te^{i\theta})|d\theta \, tg'(t)\goes 0.
  \end{equation*}
  Then, provided the RHS converges absolutely, we have
  \begin{equation*}
    \sum_{\ell\ge1}g(|\rho_\ell|)
    =\frac{1}{2\pi}\int_a^\infty\int_0^{2\pi}
    \log|f(te^{i\theta})|d\theta
    (tg''(t)+g'(t))dt
  \end{equation*}
  for any $a\in[0,|\rho_1|]$.
\end{lem}
The reader may want to read \cite{Boas*54}, for instance Theorem~8.4.1,
for general results on entire functions having only real zeroes. 

\begin{proof}
  We denote by $n(t)$ the number of zeroes of $f$ (counted with
  multiplicities) that are of modulus not more than $t$.
  We use an integration by parts to write
  \begin{align*}
    \sum_{\ell\ge1}g(|\rho_\ell|)
    &=
    -\sum_{\ell\ge1}\int_{|\rho_\ell|}^\infty g'(t)dt
    =-\alpha\int_{a}^\infty \frac{n(t)}{t}tg'(t)dt
    \\&=
    -\bigg[\int_0^t\frac{n(u)du}{u}tg'(t)\biggr]_a^\infty
    +\int_a^\infty\int_0^t\frac{n(u)du}{u}(tg''(t)+g'(t))dt
  \end{align*}
  We only have to introduce Jensen's formula in the RHS and use our
  hypothesis to get our lemma.
\end{proof}

When used with $g(t)=1/t$ and appealing to Lemma~\ref{majD}, we get the
following important result. 
\begin{thm}
  \label{nuclear}
  The hypothesis on $V$ being as above, the operator $\Vscr_0$ is nuclear. Furthermore, 
  it satisfies $\sum_{\ell\ge1}\lambda_\ell=0$ and
  \begin{equation*}
    \sum_{\ell\ge1}|\lambda_\ell|
    \ll
    \sqrt{\int_0^1|V(z)|^2dz}
    \ e^{-c_3\sqrt{\log(1+\|V\|_\infty\int_0^1|V'(t)|dt/\int_0^1|V(t)|^2dt)}}
  \end{equation*}
  for some positive constant $c_3$ that depends only on $B$
  and~$c$. In particular, we have
  \begin{equation}
    \label{localbound}
    |\lambda_\ell|\ll \sqrt{\int_0^1|V(z)|^2dz}/\ell.
  \end{equation}
\end{thm}
In our case of application, the $L^2$-norm of $V$ is controlled
by Lemma~\ref{norm}.

\begin{proof}
  On combining Lemma~\ref{Jensenplus} together with
  Lemma~\ref{majD}, we readily find that
  \begin{align*}
    \sum_{\ell\ge1}|\lambda_\ell|
    &\ll
    M\int_{1/|\lambda_1|}
    t\,e^{-\frac{c}{30}\sqrt{\log(Mt+1)}}\frac{dt}{t^2}
    \\&\ll
    M\int_{M/|\lambda_1|}
    e^{-\frac{c}{30}\sqrt{\log(t+1)}}\frac{dt}{t}
    \ll
    |\lambda_1|e^{-\frac{c}{30}\sqrt{\log(M|\lambda_1|^{-1}+1)}}.
  \end{align*}
  By \eqref{localbound1} with $\ell=1$, we find that $|\lambda_1|\le
  \sqrt{2\int_0^1|V(z)|^2dz}$, hence the bound for
  $\sum_{\ell\ge1}|\lambda_\ell|$.  Lidskii's Theorem  then applies
  giving us that $\sum_{\ell\ge1}\lambda_\ell=\int_0^1V(y-y)dy=0$.
\end{proof}

\subsection{Oscillation of the eigenvalues}

Let us consider the eigenvalues of~$\Vscr_0$. At least one of them is
positive and at least one of them is negative because
\begin{equation*}
  \sum_{\ell\ge1}\lambda_\ell=0
\end{equation*}
and $V$ is not identically 0.  Proving that infinitely many of them are
positive or negative seems to be more difficult, if true.

\subsection{A Mercer Theorem}

Let us select a complete system of non-zero eigenvectors
$(G_\ell)_{\ell\ge1}$ associated with the eigenvalues
$(\lambda_\ell)_{\ell}$ that are repeated according to multiplicity and
arranged in non-increasing order of their absolute values.
\begin{thm}
  \label{Mercer}
  For every positive integer $N$, we have
  \begin{equation*}
    \max_{y,y'\in[0,1]}\Bigl|
    V(y'-y)-\sum_{\ell\le N}\lambda_\ell
   \overline{G_\ell}(y)G_\ell(y')
    \Bigr|\le |\lambda_{N+1}|.
  \end{equation*}
\end{thm}

This theorem contains the value of the trace. Indeed,
on selecting $y'=y$, we get
$\sum_{\ell\ge1}\lambda_\ell|{G_\ell}(y)|^2=0$; we then integrate
this equality over~$y$ and recover the trace $\sum_{\ell\ge1}\lambda_\ell=0$.

\begin{proof}
 We have, for
any $y'$ in $[0,1]$ and any $L^2$-function $f$:
\begin{align*}
  \int_0^1V(y'-y)f(y)dy
  &=\sum_{\ell\ge1}\lambda_\ell(f|G_\ell)G_\ell(y')
  \\&=\sum_{\ell\le N}\lambda_\ell
  \int_0^1 f(y)\overline{G_\ell}(y)dyG_\ell(y')
  +\sum_{\ell\ge N+1}\lambda_\ell(G_\ell|f)G_\ell(y').
\end{align*}
This implies that, for any test function $h$, we have
\begin{equation}
  \label{eq:16}
  \int_0^1 h(y')\Bigl|
  \int_0^1 
  \Bigl(V(y'-y)-\sum_{\ell\le N}\lambda_\ell
   \overline{G_\ell}(y)G_\ell(y')\Bigr)f(y)dy\Bigr|dy'
   \ll |\lambda_{N+1}|\|f\|\|h\|
\end{equation}
by using Cauchy's inequality and
\begin{align*}
  \int_0^1\Bigl|\sum_{\ell\ge N+1}\lambda_\ell(G_\ell|f)G_\ell(y')\Bigr|^2dy
  &\le
  \int_0^1\sum_{\ell\ge N+1}|\lambda_\ell(G_\ell|f)|^2|G_\ell(y')|^2dy
  \\&\le
  |\lambda_{N+1}|^2\sum_{\ell\ge N+1}|(G_\ell|f)|^2
  \le |\lambda_{N+1}|^2\|f\|^2.
\end{align*}
Select a point $y_0$ from $(0,1)$ and a positive $\ve$ such that
$[y_0-\ve,y_0+\ve]\subset[0,1]$. We take $f=\1_{[y_0-\ve,y_0+\ve]}$ and get
\begin{equation*}
  \int_0^1 h(y)
  \biggl|\frac{1}{2\ve}\int_{y_0-\ve}^{y_0+\ve}
  \Bigl(V(y'-y)-\sum_{\ell\le N}\lambda_\ell
   \overline{G_\ell}(y)G_\ell(y')\Bigr)dy\biggr|dy'
   \ll |\lambda_{N+1}|\|h\|.
\end{equation*}
However we have
\begin{align*}
  V(y'-y_0)-\frac{1}{2\ve}\int_{y_0-\ve}^{y_0+\ve}
  V(y'-y)dy
  &\ll
  \frac{1}{2\ve}\int_{y_0-\ve}^{y_0+\ve}\exp\bigl(-c\sqrt{\min(1,A\ve)}\bigr)dy
  \\&\ll
  \exp\bigl(-c\sqrt{\min(1,A\ve)}\bigr)
\end{align*}
which tends to zero with $\ve$. The same applies to $y\mapsto
\sum_{\ell\le N}\lambda_\ell \overline{G_\ell}(y)G_\ell(y')$.
In case of the two endpoints $y_0=0$
and $y_0=1$, we simply select $f=\1_{[0,\ve]}$ in the first case and
$f=\1_{[1-\ve,1]}$ in the second one.
We then employ the same trick regarding the variable $y'$.
 We leave the details to the reader.
\end{proof}

\subsection{Influence of the Riemann Hypothesis}

As we already mentioned, under the Riemann Hypothesis, the function
$y'\mapsto V(y'-y)$ is uniformly H\"older with exponent $1/2-\ve$ for any
$\ve>0$. In which case, \cite[Theorem 16.3-1]{Cochran*72} gives us that
\begin{equation*}
  \sum_{\ell\ge1}|\lambda_\ell|^{p}\ll_p 1
\end{equation*}
for every $p>4/5$. This implies that the number of eigenvalues below $t$, say
$n(t)$, satisfies $n(t)\ll_\ve t^{4/5+\ve}$ under the Riemann
Hypothesis.



\subsection{Bounds from Fourier analysis and non-negativity}

Since the function $V$ is even over $\mathbb{R}$ its Fourier transform is (a
cosine transform and hence) real valued. In practice, we will use $V(u)=W^\star(\tau
u/h)$ where  $\hat{W}^\star$ is also given by~\eqref{neat}; hence we can bound above
the values of the eigenvalues when $W$ is assumed to be non-negative.

\begin{thm}
\label{bounds}
  Assume that $\hat{V}(u)\le M_1$ when
$u\in\mathbb{R}$. Then the eigenvalues of $\Vscr_0$ are not more than $M_1$.
There exists a positive constant $c_4$ such that, if we further assume 
that $\hat{V}(u)\le M_2$ when $|u|\ge U_2$ 
for some positive
parameters $M_1>M_2$ and $U_2$,
then the eigenvalues of $\Vscr_0$ are not more than $M_1-ce^{-c_4 U_2}$ for some positive
constant $c$ depending on $M_1$ and $M_2$ (but not on $V$ nor on~$U_2$).
\end{thm}

The proof uses
F.I. Nazarov's form \cite{Nazarov*92}, \cite{Nazarov*93} of the Amrein-Berthier
Theorem \cite{Amrein-Berthier*77} (see also \cite[Section
4.11]{Havin-Joricke*94} in the monograph of V. Havin \& B. J\"oricke) that we now recall.
\begin{thm}[Nazarov]
  \label{Nazarov}
  There exist two positive constants $c_4,c_8$ such that, for any
  measurable subsets $E$ and $\Sigma$ of $\mathbb{R}$ of finite
  measure, and for any $f\in L^2(\mathbb{R})$, we have
  \begin{equation*}
    \|f\|_2^2\le c_8\, e^{c_4|E||\Sigma|}\biggl(
    \int_{x\in\mathbb{R}\setminus E}|f(x)|^2dx
    +
    \int_{u\in\mathbb{R}\setminus \Sigma}|\hat{f}(u)|^2du
    \biggr).
  \end{equation*}
\end{thm}

We thank P. Jaming for giving some advice on this result, for pointing
out that a theorem of V.N. Logvinenko \and Ju.F. Sereda
\cite{Logvinenko-Sereda*74} would be enough here (since we
consider only the case when $E$ and $\Sigma$ are intervals), and for
giving us the reference to the paper \cite{Kovrijkine*00} of
O. Kovrijkine that gives a simpler proof. P. Jaming also told us that he
believes $c_8=300$ and $c_4=120$ to be an admissible choice.

\begin{proof}[Proof of Theorem~\ref{bounds}]
  We write
  \begin{equation*}
    V(y-y')=\int_{-\infty}^\infty \hat{V}(u)e(-u(y-y'))du
  \end{equation*}
  and thus, for any $G\in L^2([0,1])$, we have
  \begin{equation}
    \label{squares}
    [G,\Vscr_0(G)]=\int_{-\infty}^\infty \hat{V}(u)|\hat{G}(u)|^2du.
  \end{equation}
  Some comments are called for here. We have
  \begin{equation*}
    \hat{G}(u)=\int_{-\infty}^\infty G(v)e(-uv)dv
  \end{equation*}
  i.e. we have extended $G$ from $[0,1]$ to $\mathbb{R}$ by 0
  outside. By the result of Nazarov cited above, its Fourier transform
  is not accumulated on an interval. More precisely, on selecting
  $E=[0,1]$ and $\Sigma=[-U_2,U_2]$ in Theorem~\ref{Nazarov}, we find
  that
  \begin{equation*}
    \int_{|u|\ge U_2}|\hat{G}(u)|^2du\ge e^{-c_4 U_2}\|G\|_2^2/c_8
  \end{equation*}
  and thus
  \begin{equation*}
    [G,\Vscr_0(G)]\le (M_1-e^{-c_4 U_2}c_8^{-1}(M_1-M_2))\|G\|_2^2.
  \end{equation*}
  The theorem follows readily.
\end{proof}
In between, \eqref{squares} implies the following.
\begin{lem}
  The eigenvalues of $\Vscr_0$ lie inside $[-\min \hat{V}(u), \max \hat{V}(u)]$.
\end{lem}

\subsection{Spectral decomposition of $\Vscr$ from the one of $\Vscr_0$}

Now that we have the spectral decomposition of $\Vscr_0$ with couples
$(\lambda_\ell, G_\ell)$, we recover a spectral decomposition of
$\Vscr|_{L^2_*(\Z{h})}$ (the restriction of $\Vscr$ to $L^2_*(\Z{h})$), by
considering the eigenvectors $\Base_{h,\chi}\otimes G_\ell$, where
$\Base_{h,\chi}$ comes from Theorem~\ref{ExplicitUhtildeh}. These eigenvectors
are of norm~1 and are associated with the eigenvalues
$\lambda_\ell$. When we want to refer to the eigenvalues of
$\Vscr|_{L^2_*(\Z{h})}$, we use the notation $\lambda_\ell$ and we add the
${}$ superscript for $\Vscr_0$. We go from the latter to the former by
repeating $\phi(h)$ times each eigenvalue.


\section{From global to local: two embeddings}


The hermitian product on
$\N{N}$ is given by~\eqref{defscalarN}.

\subsection*{From the sequence $\varphi$ to a local function}

We explore the embedding defined in~\eqref{defFNh}.

Concerning~\eqref{defN'}, we specify here that we could select a uniform value
for $N'$, typically $N+H$ where $H$ is a bound to be chosen (like $\exp
c_1\sqrt{\log N}$). Since $N$ is supposed to be much larger than~$H$, the
introduction of this parameter in the next definition is only to correct some
effects on the border of our domain, see the proof of Lemma~\ref{iso}
below. There are several ways to handle this situation, we could have
considered $[0,2]$ rather than $[0,1]$ in the definition of $X$ or we could
also have kept $N$ and $[0,1]$ and simply replaced the equality of
Lemma~\ref{iso} by an equality with an error term and carried this error term
throughout the proofs. The choice above has the advantage of being independent
of an external upper bound (but is not henceforth canonical).

As a consequence, we note directly\footnote{Indeed, under the stated condition
on $y$, we have $[N'y/h]\ge [N'/h]\ge N/h$ and thus the index $b+h[N'y/h]$ is
strictly larger than $N$.} here that
\begin{equation}
  \label{border}
  \Gamma_{N,h}(\varphi)(b,y)=0\quad\text{when $y\ge [N'/h]h/N'$.}
\end{equation}
The fundamental property of $\Gamma_{N,h}$ is that it preserves the hermitian product
up to a multiplicative constant (but is {\sl not} isometric as it is not onto).
\begin{lem}
  \label{iso}
  For any positive integer $h\le N'-N$, we have 
  \begin{equation*}
    \tfrac{N}{N'}[\varphi,\psi]_N=\langle \Gamma_{N,h}(\varphi),\Gamma_{N,h}(\psi)\rangle_h.
  \end{equation*}
\end{lem}
\noindent\emph{The reader should notice a notational difficulty here: the norm
  $\|\varphi\|_2$ that we have used up to now corresponds to the scalar
  product \emph{only} up to the scalar $1/N$. We will thus refrain from using
  $\|\varphi\|_N^2$ as a shortcut to $[\varphi,\varphi]_N$.}
\begin{proof}
  Indeed, we have
  \begin{align*}
    \langle \Gamma_{N,h}(\varphi),&\Gamma_{N,h}(\psi)\rangle_h
    \\&=
    \frac{1}{h}\sum_{1\le b\le h}
    \int_0^1 \Gamma_{N,h}(\varphi)(b,y)\overline{\Gamma_{N,h}(\psi)(b,y)}dydz
    \\&=
    \frac{1}{h}
    \sum_{1\le b\le h}
    \sum_{0\le k\le \frac{N'}{h}-1}
    \int_{kh/N'}^{(k+1)h/N'}
    \Gamma_{N,h}(\varphi)(b,y)\overline{\Gamma_{N,h}(\psi)(b,y)}
    dydz
    \\&\qquad\qquad+\frac{1}{h}\sum_{1\le b\le h}
    \int_{[\frac{N'}{h}]h/N'}^{1} 
    \Gamma_{N,h}(\varphi)(b,y)\overline{\Gamma_{N,h}(\psi)(b,y)}dy
    \\&=
    \frac{1}{N'}\sum_{1\le b\le h}\sum_{n\equiv b[h]}\varphi_n\overline{\psi_n}
  \end{align*}
  on employing~\eqref{border}.
  Hence the result.
\end{proof}

\subsection*{Local adjoint}

For every $\varphi$, the linear functional $f\mapsto\langle
\Gamma_{N,h}(\varphi),f\rangle$ can be uniquely represented in the form
$[\varphi, \LL_{N,h}(f)]_N$, i.e. we have
\begin{equation}
  \label{eq:11}
  \langle \Gamma_{N,h}(\varphi),f\rangle=[\varphi, \LL_{N,h}(f)]_N.
\end{equation}
The functional $f\mapsto\LL_{N,h}(f)$ is of course linear.
We find that
\begin{align*}
  [\varphi, \LL_{N,h}(f)]_N
  &=
  \frac{1}{h}
  \sum_{c\mod h}
  \int_0^1
  \varphi_{c+h[Nx/h]}
  \overline{f(c,y)}dy
  \\&=
  \frac{1}{h}
  \sum_{c\mod h}
  \sum_{\substack{n\le N, \\ n\mod h=c}}\varphi_n
  \int_{[\frac{n-c}{N'},\frac{n-c+h}{N'}]}
  \overline{f(c,y)}dy
\end{align*}
and thus, for any integer $n\le N$, we deduce the following explicit expression:
\begin{equation}
  \label{explicitLL}
  \LL_{N,h}(f)(n)=\frac{N}{h}\int_{[\frac{n-\sigma_h(n)}{N'},\frac{n-\sigma_h(n)+h}{N'}]}
  f(\sigma_h(n),y)dy.
\end{equation}
We conclude from $\tfrac{N}{N'}[\varphi,\psi]_N=\langle
\Gamma_{N,h}(\varphi),\Gamma_{N,h}(\psi)\rangle=[\LL_{N,h}\Gamma_{N,h}\varphi,\psi]_N$ that 
\begin{equation}
  \label{eq:15}
  \LL_{N,h}\Gamma_{N,h}=\tfrac{N}{N'}\Id.
\end{equation}
And some easy manipulations tell us that $\Gamma_{N,h}\LL_{N,h}=\tfrac{N}{N'}P_h$
where $P_h$ is the orthogonal projector on $\text{Im}\,
\Gamma_{N,h}=\Gamma_{N,h}\bigl(L^2(\N{N})\bigr)$.
\begin{proof}
  Indeed, we find that, for any $\varphi$ and $\psi$, we have
  \begin{align*}
    \langle \Gamma_{N,h}\LL_{N,h}\Gamma_{N,h}\LL_{N,h}\varphi,\psi\rangle
    &=
    [\LL_{N,h}\Gamma_{N,h}\LL_{N,h}\varphi,\LL_{N,h}\psi]_N
    \\&=
    \tfrac{N}{N'}[\LL_{N,h}\varphi,\LL_{N,h}\psi]_N
    =
    \tfrac{N}{N'}[\Gamma_{N,h}\LL_{N,h}\varphi,\psi]_N.
  \end{align*}
  We conclude from these equalities that
  $(\Gamma_{N,h}\LL_{N,h})^2=\tfrac{N}{N'}\Gamma_{N,h}\LL_{N,h}$. The conclusion is easy.
\end{proof}

\subsection*{Pure embeddings}

It will be clear in a moment that, if $\Gamma_{N,h}(\varphi)$ is easier to
grasp from a geometrical viewpoint, our object is in fact
$\Pure_{N,h}=U_{\tilde h\va h}\circ \Gamma_{N,h}$ as already defined in~\eqref{defPure},
i.e.  the orthonormal projection of $\Gamma_{N,h}$ on the space $L^2_*(X_h)$
(see section~\ref{infinity})\footnote{The choice of notation $\Gamma_{N,h}^*$
  would lead to confusion since adjoints are present in the latter theory.}.
We call the function $\Pure_{N,h}$ the \emph{pure embedding}. From Theorem~\ref{ExplicitUhtildeh}, we get
\begin{equation}
  \label{ortho}
  \Pure_{N,h}(\varphi)(b,y)
  =
  \frac{1}{h}\sum_{c\mod h}c_h(b-c)
  \Gamma_{N,h}(\varphi)(c,y)
\end{equation}
from which we readily compute that
\begin{equation}
  \label{ortho*}
  \Pure_{N,h}^*(\varphi)(n)
  =
  \frac{1}{h^2}\sum_{b\mod h}c_h(b-n)
  \int_{\frac{n-\sigma_h(n)}{N'}}^{\frac{n-\sigma_h(n)+h}{N'}}
  \varphi(b,z)dz.
\end{equation}
Note that
\begin{align}
  \|\Pure_{N,h}(\varphi)\|^2
  &=\notag
  \frac{1}{h^2}
  \sum_{b\mod h}
  \int_0^1\biggl|
  \sum_{c\mod h}c_h(b-c)\Gamma_{N,h}(\varphi)(c,y)\biggr|^2dy
  \\&=\label{cute}
  \frac{1}{h^2}
  \sum_{b\mode h}
  \int_0^1
  \biggl|
  \sum_{0\le a< h}
  \Gamma_{N,h}(\varphi)(a,y)e(-ab/h)
  \biggr|^2dy.
\end{align}
Eq.~\eqref{cute} shows immediately (by extending the summation in $b$ to
all of $\Z{h}$) that $\|\Pure_{N,h}(\varphi)\|\le
\|\Gamma_{N,h}(\varphi)\|$, a fact that could have been more easily obtained
by noticing that the norm of an orthogonal projection is surely not more than
the initial norm. We can also get 
an explicit expression of $\|\Pure_{N,h}(\varphi)\|^2$ in terms of $\varphi$:
\begin{align}
  \|\Pure_{N,h}(\varphi)\|^2
  &=\notag
  \frac{1}{h^2}
  \sum_{b\mode h}
  \sum_{k\ge0}\int_{kh/N'}^{(k+1)h/N'}
  \biggl|\sum_{0\le a<h}\varphi_{a+kh}e(-ba/h)\biggr|^2dy
  \\&=\label{ExplicitG}
  \frac{1}{hN'}
  \sum_{b\mode h}
  \sum_{k\ge0}
  \biggl|\sum_{n/ [n/h]=k}\varphi_{n}e(-bn/h)\biggr|^2.
\end{align}

\section{Theorem~\ref{Precise} in functional form}

We start with an easy lemma.
\begin{lem}
  \label{repSah}
  We have
  \begin{equation*}
    S\Bigl(\varphi;\frac{a+\vartheta}{h}\Bigr)
    =
    \frac{N'}{h}
    \sum_{1\le b\le h}\int_0^1
    \Gamma_{N,h}(\varphi)(b,y)e\Bigr(\frac{ab}{h}\Bigr)
    e\Bigl(\Bigl(\frac{b}{h}+[N'y/h]\Bigr)\vartheta\Bigr)dy.
  \end{equation*}
\end{lem}

\begin{proof}
  When $m\equiv b[h]$ with $1\le b$, we have
  \begin{equation*}
    \varphi_m=\frac{N'}{h}\int_{\frac{m-b}{N'}}^{\frac{m-b}{N'}+\frac{h}{N'}}
    \Gamma_{N,h}(\varphi)(b,y)dy.
  \end{equation*}
  It is straightforward to get the lemma from this expression.
\end{proof}

When $H\le N^{1/8}(\log N)^{-3/2}$, $N\ll QH$ and $Q\le N^2$ (this condition is only to
control $\log Q$ in the error term. In practice, $Q$ is not more than $N$, but
we may want to select $Q=\text{constant}\times N$), we have the following.
\smallskip

\noindent\fbox{\vbox{%
\begin{multline}
  \label{formula11b}
  \sum_{q}\frac{W(q/Q)}{qQ\ooWo}\sum_{a\mode q}
  |S(\varphi, a/q)|^2
  =
  \|\varphi\|_2^2\Bigl(1+\Ocal\Bigl(\frac{N}{QH}\Bigr)\Bigr)
  \\
  -\mkern-6mu
  \sum_{h\le H}
  \sum_{\substack{1\le b\le h} }
  \int_0^1\mkern-9mu\int_0^1\mkern-5mu
  \frac{\Pure_{N,h}(\varphi)(b,y)
  \overline{\Pure_{N,h}(\varphi)(b, y')}}{h^2}  
  \frac{W^\star\bigl(\tau\frac{y-y'}{h}\bigr)dydy'}{Q\ooWo/N^2}
  .
\end{multline}
}}

\begin{remark}
  Most of the work below is to allow $H$ to be a power of $N$. If one
  can control the continuity of $W^\star$, like under the Riemann
  Hypothesis, then the proof is much simpler. We instead rely heavily
  on the bilinear structure.
\end{remark}

\begin{proof}
We start from~\eqref{formula51} and Lemma~\ref{repSah} to get that:
\begin{multline}
  \label{ini}
  \sum_{a\mode h}
  \int_{-U}^{U}
  \hat W_{C'}^{\star\star}(u)
  \Bigl|S\Bigl(\varphi, \frac{a}{h}+\frac{u}{hQ}\Bigr)\Bigr|^2
  du
  =\\
  \frac{N^{\prime2}}{h^2}
  \sum_{1\le b_1,b_2\le h}\int_0^1\int_0^1
  \Gamma_{N,h}(\varphi)(b_1,y)
  \Gamma_{N,h}(\varphi)(b_2,y')
  \sum_{a\mode h}e\Bigl(\frac{(b_1-b_2)a}{h}\Bigr)
  \\
  \int_{-U}^{U}\hat W_{C'}^{\star\star}(u)
  e\Bigl(\frac{b_1-b_2}{h}\frac{u}{Q}
  +\Bigl([N'y/h]-[N'y'/h]\Bigr)\frac{u}{Q}\Bigr)dudydy'.
\end{multline}
In the inner integration, we replace
  \begin{equation*}
    e\Bigl(\frac{b_1-b_2}{h}\frac{u}{Q}
    +\Bigl([N'y/h]-[N'y'/h]\Bigr)\frac{u}{Q}\Bigr)
  \end{equation*}
  by $e((y-y')Nu/(hQ))$. We call $\Delta_h(b_1,b_2,y,y')$ the difference
  of the two, integrated against $\hat
  W_{C'}^{\star\star}(u)$. We have
  \begin{equation*}
    \Delta_h(b_1,b_2,y,y')\ll 
    \int_{-U}^{U}\Bigl|\hat
    W^{\star\star}_{C'}(u)\Bigr|\min(1,|u|/Q)du.
  \end{equation*}
  This gives rise to the error term
  \begin{equation*}
    \frac{N^{\prime2}}{h^2}\sum_{1\le b_1,b_2\le h}
    \mkern-8mu
    c_h(b_1-b_2)\int_0^1\int_0^1
  \Gamma_{N,h}(\varphi)(b_1,y)
  \Gamma_{N,h}(\varphi)(b_2,y')
  \Delta_h(b_1,b_2,y,y')dydy'.
  \end{equation*}
  We get $\max |\Delta_h(b_1,b_2,y,y')|$ out, separate
  $\Gamma_{N,h}(\varphi)(b_1,y)$ from $\Gamma_{N,h}(\varphi)(b_2,y')$
  by using $2|z_1z_2|\le |z_1|^2+|z_2|^2$ and have to bound
  \begin{equation*}
    \Sigma=\frac{N^{\prime2}}{h^2}\sum_{1\le b_1,b_2\le h}
    |c_h(b_1-b_2)|\int_0^1
  |\Gamma_{N,h}(\varphi)(b_1,y)|^2dy
  \max_{y,y', b_1,b_2}
  |\Delta_h(b_1,b_2,y,y')|.    
  \end{equation*}
  We use
  \begin{align*}
    \sum_{1\le b_2\le h} |c_h(b_1-b_2)|
    &=
    \sum_{1\le b\le h} |c_h(b)|
    =\sum_{d|h}\sum_{\substack{b\mod h,\\ \gcd(b,h)=h/d}}\mu^2(d)\frac{\phi(h)}{\phi(d)}
    \\&=
    \sum_{d|h}\mu^2(d)\phi(d)\frac{\phi(h)}{\phi(d)}
    =2^{\omega(h)}\phi(h).
  \end{align*}
  This and the isometrical property of $\Gamma$ leads to
  \begin{equation*}
    \Sigma=\frac{2^{\omega(h)}\phi(h)N^{\prime}}{h}
    \|\varphi\|_2^2
      \max_{b_1,b_2,y,y'}
      |\Delta_h(b_1,b_2,y,y')|.
  \end{equation*}
  Next by using Lemma~\ref{formulacheckW}, we check that $|\Delta_h(b_1,b_2,y,y')|\ll
  C'/Q$. The total error term is $\ll \sum_{h\le
    H}(Q h)^{-1}C'2^{\omega(h)}N^{3/2}\|\varphi\|_2^2/Q\ll
  C'N^{3/2}\|\varphi\|_2^2(\log H)^{2}/Q^2$ which we call $E_1$.
  Thus we have reduced the right-hand side de \eqref{ini} to
  \begin{multline}
    \label{eq:20}
  \frac{N^{\prime2}}{h^2}
  \sum_{1\le b_1,b_2\le h}\int_0^1\int_0^1
  \Gamma_{N,h}(\varphi)(b_1,y)
  \Gamma_{N,h}(\varphi)(b_2,y')
  c_h(b_1-b_2)
  \\
  \int_{-U}^{U}\hat W^{\star\star}_{C'}(u)
  e\Bigl((y-y')\frac{Nu}{hQ}\Bigr)dudydy'.
  \end{multline}
  By \eqref{ortho}, this is also
  \begin{multline*}
    \frac{N^{\prime2}}{h}
  \sum_{1\le b_2\le h}\int_0^1\int_0^1
  \Pure_{N,h}(\varphi)(b_2,y)
  \overline{\Gamma_{N,h}(\varphi)(b_2,y')}
  \\
  \int_{-U}^{U}\hat W^{\star\star}_{C'}(u)
  e\Bigl((y-y')\frac{Nu}{hQ}\Bigr)dudydy',
  \end{multline*}
  which, by orthogonality, is also
  \begin{multline*}
    \frac{N^{\prime2}}{h}
  \sum_{1\le b\le h}\int_0^1\int_0^1
  \Pure_{N,h}(\varphi)(b,y)
  \overline{\Pure_{N,h}(\varphi)(b,y')}
  \\
  \int_{-U}^{U}\hat W^{\star\star}_{C'}(u)
  e\Bigl((y-y')\frac{Nu}{hQ}\Bigr)dudydy'.
  \end{multline*}
  We want to replace $\hat W^{\star\star}_{C'}(u)$ by $\hat W^{\star}$.
We assume $U\le C'/2$, hence $\hat W^{\star}(u)=\hat
W^{**}_{C'}(u)+\text{Constant}$ when $|u|\le U$ and this constant is
$\Ocal(1/C')$. Again using $2|z_1z_2|\le |z_1|^2+|z_2|^2$ on
$\Pure_{N,h}(\varphi)$, and noting that (with $s=NU/hQ$)
\begin{equation*}
  \int_0^1
  \frac{\sin((y'-y)s)}{(y'-y)s}dy
  =
  \int_{y's-y}^{y's}
  \frac{\sin x}{x}dx\ll1
\end{equation*}
uniformly in $s$ and $y'$, we get an error term of size $\Ocal((\log
H)N'\|\varphi\|_2^2/(QC'))$.  We finally want to extend the path of
integration in $u$ to infinity. Again using $2|z_1z_2|\le
|z_1|^2+|z_2|^2$, this means bounding
  \begin{equation*}
    A=\int_0^1 \biggl|\int_{U}^\infty\hat W^\star(u)
  e\Bigl((y-y')\frac{Nu}{hQ}\Bigr)du\biggr|dy
  \end{equation*} 
  and similarly with $y'$. We employ Cauchy's inequality and open the square, getting:
  \begin{equation*}
    A^2\ll 
    \int_{U}^\infty\int_{U}^\infty
    \hat W^\star(u_1)
    \overline{\hat W^\star(u_2)}
    \int_0^1
    e\Bigl((y-y')\frac{N(u_1-u_2)}{hQ}\Bigr)du.
  \end{equation*}
  We employ Lemma~\ref{compFourierTransforms} on $u_1$ and $u_2$. When
  $|u_1-u_2|\le 1$, we get the contribution $\Ocal(1/U)$; When $|u_1-u_2|\ge 1$, we integrate in $y$
  and get the contribution 
  \begin{equation*}
    \int_{U}^\infty\int_{U}^\infty\frac{du_1du_2}{u_1u_2(1+|u_1-u_2|)}.
  \end{equation*}
  On splitting the path of integration on $u_2$ in
  $[U,\max(U,u_1/2)]$, followed by $[\max(U,u_1/2),2u_1]$ and
  finally by $[2u_1,\infty)$, we readily see that this integral is 
$\Ocal((\log
  U)/U)$. Summing over $h$ gives the contribution
  \begin{equation*}
    \frac{N}{Q\ooWo}\sum_{h\le H}\frac{1}{h}
    \|\varphi\|_2^2\sqrt{\frac{hQ(\log U)}{NU}}
      \ll
      \frac{\sqrt{NH\log U}}{\sqrt{UQ}}\|\varphi\|_2^2=E_2.
  \end{equation*}
In total, we get the error term bounded above by a constant multiple of
\begin{equation*}
  \biggl(
  \frac{N}{QH}
  +
  \frac{N\log H}{QC'}
  +
  \sqrt{\frac{NH\log U}{UQ}}
  +
  \frac{C'N^{3/2}}{Q^2}(\log H)^2
  +
  \frac{H^2}{U}
  +
  \frac{H^2+(\log N)^5}{C}
  \biggr)\|\varphi\|_2^2.
\end{equation*}
It is best to take $U$ as large as possible, so we select $U=C'/2$.
In turn, we select $C'=QH^{1/3}/N^{2/3}$ and we check that $C'=C$ (see~\eqref{defC}).
The error term becomes not more than a constant multiple times
\begin{equation*}
  \biggl( 
  \frac{N}{QH}
  +
  \frac{N^{5/3}\log H}{Q^2H^{1/3}}
  +
 \frac{N^{5/6}H^{1/3}}{Q}(\log N)^2
  +
  \frac{H^2+(\log N)^5}{QH^{1/3}}N^{2/3}
  \biggr)\|\varphi\|_2^2.
\end{equation*}
We then check that this reduces to
\begin{equation*}
  \biggl( 
  \frac{N}{QH}
  +
  \frac{N}{QH}\frac{NH^{5/3}}{QHN^{1/3}}
  +
 \frac{N^{5/6}H^{1/3}}{Q}(\log N)^2
  \biggr)\|\varphi\|_2^2
\end{equation*}
when $H\le N^{1/8}$. And we check further that
${N^{5/6}H^{1/3}}Q^{-1}(\log N)^2\ll N/(QH)$ when $H\le N^{1/8}(\log N)^{-3/2}$. The second term equally disappears, as $N\ll QH$.
\end{proof}
Herv\'e Queff\'elec has kindly pointed out to me that when $q=1$, this process
bears similarities with the one devised independently by \cite{Shampine*66}
and \cite{Widom*66}, and which is nicely presented in \cite[Section 3]{Bottcher*13}.

On recalling the definition of the operator $\Vscr_{\tau,h}$
in~\eqref{defVscrtauh}, here is another manner of writing~\eqref{formula11b}:
\smallskip

\noindent\fbox{\vbox{%
\begin{multline}
  \label{formula11bb-1}
  \sum_{q}\frac{W(q/Q)}{qQ}\sum_{a\mode q}
  |S(\varphi, a/q)|^2
  =
  \ooWo\|\varphi\|_2^2(1+\Ocal(\tau/H))
  \\
  -
  N\sum_{h\le H}\frac{\tau}{h}
  \bigl[\Pure_{N,h}(\varphi)|\Vscr_{\tau,h} \Pure_{N,h}(\varphi)\bigr]_{h\times[0,1]}
  \\
  (H\ll N^{1/8}(\log N)^{-3/2},N\ll QH, Q\ll N^2)
  .
\end{multline}
}}
\smallskip


\section{Using spectral analysis}
Formula~\eqref{formula11bb-1} involves the operators
$\Vscr_{\tau,h}\circ U_{\tilde{h}\va h}$. In this section, we first
diagonolize them as local operators (i.e. on a space that depends on $h$), and control the dependance in $h$ and $\tau$.
We then lift this diagonalization to the global space (where the
sequence $\varphi$ lives) and show that the resulting family of
eigenvectors, $h$ varying, is near-orthonormal (see
Lemma~\ref{beauty}). We encounter a problem (that may be only
technical): the control we have of the modulus of continuity of these
eigenfunctions is weak when they are associated with very small
eigenvalues. But then, their total contribution is small, and we then
introduce a trade-off point with the condition $|\lambda_{h,\ell}|\ge
\xi\eta_0(N)^{1/8}$. We conclude this part with another consequence of
the near-orthonormality which enables us to control the quadratic form
resulting from taking some upper bound for the eigenvalues. This is
required because, when using \eqref{eq:7} to simplify our statement,
the near-orthogonality is not apparent anymore.

\subsection{Decomposing the implied operators}
\label{divertimento}

The operator $\Vscr_{\tau,h} $ does not touch the $b$-variable, from
which we infer that (recall the definition of the rothonormal
projector $U_{\tilde{h}\va h}$ in \eqref{defUktildevahtilde})
\begin{equation*}
  U_{\tilde{h}\va h}\circ \Vscr_{\tau,h} = \Vscr_{\tau,h} \circ U_{\tilde{h}\va h}.
\end{equation*}
This has two consequences: first the image of $\Vscr_{\tau,h}$ lies
inside $L^2_*(X_h)$ and second, its couples eigenvalues / eigenvectors
are simply (tensor) products of the respective couples coming from the
two operators:
\begin{equation*}
  F\in L^2(\Z{h})\mapsto\biggl(b\mapsto\frac{1}{h}\sum_{c\mod
    h}c_h(b-c)F(c)\biggr)
\end{equation*}
where the only difference with the operators $U_{\tilde{h}\va h}$ and
$\Vscr_{\tau,h}$ are the spaces. The first operator is covered by
Theorem~\ref{ExplicitUhtildeh}. We are left with the second one which
belongs to the class described in Section~\ref{ACDO} (if we ignore the
first variable, as we may). The regularity assumptions $(R_1)$,
$(R_2)$ and $(R_3)$ are met by Lemma~\ref{boundWstarC}.

\subsection{Diagonalisation in the local spaces}
\label{SpectralAnalysis}

We use the eigenvectors / eigenvalues
$(G_{h,\ell,\chi},\lambda_{h,\ell})_{\chi,\ell}$ of 
$\Vscr_{\tau,h}$ as well as the ones of $\Pure_{\tau,h}$ (see
Theorem~\ref{ExplicitUhtildeh}) to write
\begin{equation*}
  \bigl[\Pure_{N,h}(\varphi)|\Vscr_{\tau,h}
  \Pure_{N,h}(\varphi)\bigr]_{h\times[0,1]}
  =
  \sum_{\ell\ge1}\lambda_{h,\ell}
  \sum_{\chi \mod h}\bigl[\Pure_{N,h}(\varphi)|\Base_{h,\chi}\otimes G_{h,\ell}\bigr]_{h\times[0,1]}^2.
\end{equation*}
We then divide this quantity by $h$ and sum that over $h$.
 Before proceeding,
let us note the following lemma.
\begin{lem}
  \label{norm}
  \begin{equation*}
    \int_0^1 \Bigl|W^\star\Bigl(\frac{\tau z}{h}\Bigr)\Bigr|^2dz\ll
    \begin{cases}
      \exp-2c_0\sqrt{\log(2+h/\tau)},\\
      h/\tau.
    \end{cases}
  \end{equation*}
  We will use the latter when $h\le 2\tau$ and the former otherwise. It is
  however better for questions of uniformly to state them in general
\end{lem}
\begin{proof}
    When $h\ge 2\tau$, we use Lemma~\ref{boundWstarC} and bound the value
  $|W^\star({\tau z}/{h})|$ by
  $\Ocal(\exp-c_0\sqrt{\log(h/\tau)})$. 
  When $h\le 2\tau$, we use
  \begin{equation*}
    \int_0^1
    W^\star\Bigl(\frac{\tau z}{h}\Bigr)^2dz
    \le 2
    \int_0^1
    W^\star\Bigl(\frac{\tau w}{h}\Bigr)^2dw
    \le \frac{2 h}{\tau}
    \int_0^{\tau/h}
    W^\star(w)^2dw\ll h/\tau.
  \end{equation*}
  The lemma is proved.
\end{proof}

Since $|\lambda_{h,\ell}|\ll 1/\ell$ by~\eqref{localbound} and
Lemma~\ref{norm}, we can explicitly shorten the spectral decomposition in
(recall also Lemma~\ref{iso})
\begin{multline*}
  \bigl[\Pure_{N,h}(\varphi)|\Vscr_{\tau,h}
  \Pure_{N,h}(\varphi)\bigr]_{h\times[0,1]}
  =
  \\
  \sum_{\ell\le L}\lambda_{h,\ell}
  \sum_{\chi \mod h}\bigl[\Pure_{N,h}(\varphi)|\Base_{h,\chi}\otimes
  G_{h,\ell}\bigr]_{h\times[0,1]}^2
  +\Ocal(N^{-1}\|\varphi\|_2^2/L).
\end{multline*}
We can similarly restrict the summation to $|\lambda_{h,\ell}|\ge
\eta_0(N)^{1/4}$ (with $\eta_0(x)=\exp-\frac{c_0}{2}\sqrt{\log x}$) and get,
for any $\xi\in[0,1]$:
\begin{multline*}
  \bigl[\Pure_{N,h}(\varphi)|\Vscr_{\tau,h}
  \Pure_{N,h}(\varphi)\bigr]_{h\times[0,1]}
  =
  \\
  \sum_{\substack{\ell\le L,\\ |\lambda_{h,\ell}|\ge
      \xi\eta_0(N)^{1/4}}}\lambda_{h,\ell}
  \sum_{\chi \mod h}\bigl[\Pure_{N,h}(\varphi)|\Base_{h,\chi}\otimes
  G_{h,\ell}\bigr]_{h\times[0,1]}^2
  \\+\Ocal\bigl(N^{-1}\|\varphi\|_2^2(\eta_0(N)^{1/4}+1/L)\bigr).
\end{multline*}
The parameter $\xi$ is here for flexibility, in case we want the sum not to
depend on the parameter~$N$.
We may rewrite formula~\eqref{formula11bb-1} by introducing the
adjoint $\Pure_{N,h}^*$ of $\Pure_{N,h}$, as follows.
\smallskip

\noindent\fbox{\vbox{%
\begin{multline}
  \label{formula11bb}
  \sum_{q}\frac{W(q/Q)}{qQ}\sum_{a\mode q}
  |S(\varphi, a/q)|^2
  =
  \ooWo\|\varphi\|_2^2
  \\
  -
  N\sum_{h\le H}\frac{\tau}{h}
  \sum_{\substack{\ell\le L,\\ |\lambda_{h,\ell}|\ge
      \xi\eta_0(N)^{1/4}}}\lambda_{h,\ell}
  \sum_{\chi \mod h}\bigl[\varphi|\Pure_{N,h}^*\Base_{h,\chi}\otimes
  G_{h,\ell}\bigr]_N^2
  \\+\Ocal\biggl( 
  \Bigl(
  \frac{\xi\log H}{\exp\frac{c_0}{8}\sqrt{\log N}}
  +
  \frac{\log H}{L}+\frac{1}{H} 
  \Bigr)
  \tau
  \|\varphi\|_2^2\biggr)
  .
\end{multline}
}}
\smallskip

\noindent
Our task is now to replace $\Pure_{N,h}^*\Base_{h,\chi}\otimes
  G_{h,\ell}$ by a simpler expression.

\subsection{Approximate diagonalization in the global space}
\label{QF}

We define
\begin{equation}
  \label{defg}
  g_{h,\ell, \chi,N,\tau}=\Pure_{N,h}^* \Base_{h,\chi}\otimes
  G_{h,\ell}=\Gamma_{N,h}^*\Base_{h,\chi}\otimes G_{h,\ell}, 
\end{equation}
as well as
\begin{equation}
  \label{defgflat}
  g_{h,\ell, \chi,N,\tau}^\flat(n)=
  \frac{\tau_h(\chi,n)}{\sqrt{\phi(h)}}G_{h,\ell}\Bigl(\frac{n}{N}\Bigr).
\end{equation}

The function $g_{h,\ell, \chi,N,\tau}$ inherits from $\Base_{h,\chi}\otimes
G_{h,\ell}$ a similar separation of behaviour between arithmetic and size characters.
\begin{lem}
  \label{splitvariations}
  When $|t-n|\le N^{1/2}$ and $Q\le N$, we have
  \begin{equation*}
    g_{h,\ell,\chi,N,\tau}(n)
    =
    \frac{\tau_{h}(\chi,n)}{\sqrt{\phi(h)}}
    G_{h,\ell}\Bigl(\frac{t}{N}\Bigr)
    +\Ocal\biggl(
    \frac{ \sqrt{\phi(h)}}{|\lambda_{h,\ell}(\tau)|}
    \exp\Bigl(-\frac{c_0}{2}\sqrt{\log N}\Bigr)
    \biggr)
  \end{equation*}
  where $c_0$ is defined in Lemma~\ref{boundWstarC}. Moreover, we have
  \begin{equation*}
    |\lambda_{h,\ell}(\tau)|\|g_{h,\ell, \chi,N,\tau}\|_\infty\ll \sqrt{\phi(h)}.
  \end{equation*}
 \end{lem}
 In particular, $g_{h,\ell, \chi,N,\tau}^\flat$ approximates $g_{h,\ell, \chi,N,\tau}$.

\begin{proof}
  We
  have by \eqref{explicitLL}:
  \begin{equation*}
    g_{h,\ell,\chi,N,\tau}(n)
    =\frac{N}{h}
    \frac{1}{\sqrt{\phi(h)}}\tau_{h}(\chi,n)
    \int_{\frac{n-\sigma_{h}(n)}{N'}}^{\frac{n-\sigma_{h}(n)+h}{N'}}G_{h,\ell}(y)dy.
  \end{equation*}
  We next use \eqref{regulationcontinuity} together with
  Lemma~\ref{boundWstarC} to infer that, when $\delta\in[0,1]$, we
  have, for any $y\in[0,1-\delta]$,
  \begin{equation}
    \label{regGhell}
    |\lambda_{h,\ell}(\tau)||G_{h,\ell}(y+\delta)-G_{h,\ell}(y)|
    \ll  \exp-c_0\sqrt{-\log\min\Bigl(1,\frac{\tau\delta}{h}\Bigr)}.
  \end{equation}
  We note that $\tau\le 1$ and that $h\ge1$.
  Hence, for any $t$ such that $|t-n|\le \sqrt{N}$, we have
  \begin{align*}
    g_{h,\ell,\chi,N,\tau}(n)
    &=
    \tau_{h}(\chi,n)\frac{N}{\sqrt{\phi(h)}N'}
    G_{h,\ell}\Bigl(\frac{t}{N'}\Bigr)
    \\&\qquad
    +\Ocal\biggl(
    \frac{ \tau_{h}(\chi,n)}{|\lambda_{h,\ell}(\tau)|\sqrt{\phi(h)}}
    \exp\Bigl(-\frac{c_0}{4}\sqrt{\log N}\Bigr)
    \biggr)
  \end{align*}
  from which the stated estimate readily follows, up to two blemishes:
  the factor $N/N'=1+\Ocal(N^{-1/2})$ and the $G(t/N')$ instead of
  $G(t/N)$. This last modification follows from~\eqref{regGhell}, the
  former one being trivial.
  For the $L^\infty$-norm, note
  that (see~\eqref{propeigenvectors})
  \begin{equation*}
    \|g_{h,\ell, \chi,N,\tau}\|_\infty \le 
    2\|W^\star(\tau\cdot/h)\|_1\sqrt{\phi(h)}/|\lambda_{h,\ell}(\tau)|.
  \end{equation*}
\end{proof}

\begin{lem}
  \label{beauty}
  When $h,h'\le H\le N^{1/5}$ and $N'\le N+\sqrt{N}$, we have
  \begin{equation*}
    [g_{h,\ell, \chi,N,\tau},g_{h',\ell', \chi',N,\tau}]_N
    =
    \delta_{h=h'}
    \delta_{\ell=\ell'}
    \delta_{\chi=\chi'}
    +
    \Ocal\biggl(
    \frac{1}{\sqrt{N}}+
    \frac{H\exp\bigl(-\frac{c_0}{4}\sqrt{\log N}\bigr)}{|\lambda_{h,\ell}(\tau)||\lambda_{h'\ell'}(\tau)|}
    \biggr)
  \end{equation*}
  where $c_0$ is defined in Lemma~\ref{boundWstarC}.  The same applies
  when replacing $g_{h,\ell, \chi,N,\tau}$ and $ g_{h',\ell',
    \chi',N,\tau}$ respectively by $g_{h,\ell, \chi,N,\tau}^\flat$ and $ g_{h',\ell',
    \chi',N,\tau}^\flat$.
\end{lem}

\begin{proof}
  In order to compute $[g_{h,\ell, \chi,N,\tau},g_{h',\ell',
    \chi',N,\tau}]_N$, we split the interval $[1,N]$ in $\Ocal(N/(hh'))$
  sub-intervals containing $hh'$ consecutive integers and a remaining one. 
  We employ Lemma~\ref{splitvariations} on each sub-interval, selecting a $t$
  that is independent on the point~$n$, for instance choosing it at the origin
  of such a segment, but we shall use the freedom on choice in $t$ to shorten
  the argument below. We bound the $L^\infty$-norm of the other factor by
  Lemma~\ref{splitvariations}.
  The
  error term for each interval is
  \begin{equation*}
    \ll hh'
    \frac{ \max(\|W^\star(\tau\cdot/h)\|_1,
      \|W^\star(\tau\cdot/h')\|_1)
      \sqrt{\phi(h')\phi(h)}}{|\lambda_{h',\ell'}(\tau)|
       |\lambda_{h,\ell}(\tau)|
    }
    \exp\Bigl(-\frac{c_0}{4}\sqrt{\log N}\Bigr)
  \end{equation*}
  which we have to sum over all intervals and divide by $N$ (since the scalar
  product $[,]_N$ is scaled in this manner). The total error term incurred is
  thus
  \begin{equation*}
    \frac{ H\max(\|W^\star(\tau\cdot/h)\|_1,
      \|W^\star(\tau\cdot/h')\|_1)}{|\lambda_{h',\ell'}(\tau)|
       |\lambda_{h,\ell}(\tau)|
    }
    \biggl(\frac{H^2}{N}+
    \exp\Bigl(-\frac{c_0}{4}\sqrt{\log N}\Bigr)\biggr).
  \end{equation*}
  The summand $H^2/N$ comes from the end interval. Concerning this end
  interval, we should have had $\|W^\star(\tau\cdot/h')\|_1\cdot
  \|W^\star(\tau\cdot/h)\|_1$ rather than the maximum, but each norm is
  bounded (uniformly in $\tau$), which legitimates the bound above.

  Whenever $h\neq h'$ or $\chi\neq \chi'$, the summation over the remaining
  intervals vanishes by orthogonality. We are left with the case when $h=h'$
  and $\chi=\chi'$, in which case we have to evaluate
  \begin{equation*}
    \frac{1}{N}\sum_{n\mod h^2}
    \frac{|\tau_{h}(\chi,n)|^2}{\phi(h)}
    \sum_{t}G_{h,\ell}\Bigl(\frac{t}{N}\Bigr)
    \overline{G_{h,\ell'}\Bigl(\frac{t}{N'}\Bigr)}.
  \end{equation*}
  The sum upon $n$ is $h^2\|\Base_{h,\chi}\|_2^2=h^2$.  Concerning the sum
  upon $t$, we employ the following trick: given any interval we can use any
  $t$ from within, hence we can integrate over $t$ and divide by the length
  $h^2$ of the interval. Concerning the final interval, the reader will check
  that the contribution to include it is not more than what we already paid for
  discarding it. As a result, we get as a main term
  \begin{equation*}
    \int_0^1G_{h,\ell}(u)
    \overline{G_{h,\ell'}(u)}du
  \end{equation*}
  which is $\delta_{\ell=\ell'}$. 
\end{proof}

\subsection{External control of the eigenvectors}
Let us recall an inequality due to Selberg (given in \cite[Proposition 1]{Bombieri*71}
or in extended form in \cite[Lemma 1.1-1.2]{Ramare*06}).
\begin{lem}
  \label{Bessel}
  Let $(g_i)_{i\in I}$ be a finite family of vectors in the Hilbert
  space $\mathcal{H}$, and $f$ be some fixed vector in this same
  space. We have 
  \begin{equation*}
    \sum_{i\in I}|[f|g_i]|^2/\sum_{j\in I}|[g_i|g_j]|\le \|f\|^2.
  \end{equation*}
\end{lem}
We apply Lemma~\ref{Bessel} to the family
\begin{equation*}
  \bigl\{g_{h,\ell, \chi,N,\tau}^\flat: h\le H, \chi\mod h,\ell\le L,
  |\lambda_{h,\ell}(\tau)|\ge \eta_0(N)^{1/4}\bigr\}.
\end{equation*}
By Lemma~\ref{beauty}, we infer that
\begin{equation}
  \label{NearOptimalUp}
  N
  \mathop{\sum_{h\le H}\sum_{\ell\le L}\sum_{\chi \mod h}}_{|\lambda_{h,\ell}|\ge
      \exp-\frac{c_0}{8}\sqrt{\log N}}
  \bigl[\varphi|g_{h,\ell, \chi,N,\tau}^\flat\bigr]_N^2\le 
  \|\varphi\|_2^2\Bigl(1+ H^2L\exp-\frac{c_0}{8}\sqrt{\log N}\Bigr). 
\end{equation}

Finally we use the identity:
\begin{equation}
  \label{eq:7}
  \sum_{\chi \mod h}\biggl|
  \sum_{n\le
    N}\varphi(n)\frac{\tau_h(\chi,n)}{\sqrt{\phi(h)}}
  G\Bigl(\frac{n}{N}\Bigr)
  \biggr|^2
  =
  \sum_{a\mod^*h}\biggl|
  \sum_{n\le
    n}\varphi(n)e(na/h)
  G\Bigl(\frac{n}{N}\Bigr)
  \biggr|^2.
\end{equation}

\section{Deducing Theorem~\ref{yoddle} and~\ref{lowerbound}}

\subsection{Proof of Theorem~\protect\ref{yoddle}}

The spectral decomposition is treated in
Subsection~\ref{SpectralAnalysis}. The family $g_{h,\ell, \chi,N,\tau}$ is
defined in the next subsection at~\eqref{defg} and its near orthonormal
property in proved in Lemma~\ref{beauty}. The global decomposition is
given in
\eqref{formula11bb} once $\Pure_{N,h}^*\Base_{h,\chi}\otimes G_{h,\ell}$ is
replaced by $g_{h,\ell, \chi,N,\tau}$ and the relative sizes are taken into account.
The final property is in~\eqref{NearOptimalUp}.

Note that, for each $h$, we have at a positive and a negative
eigenvalue. Recalling~\eqref{t2}, we see that
$\max_{\ell}|\lambda_{h,\ell}(\tau)|$ goes to zero.
Hence these
positive or negative values of $\lambda_{h,\ell}(\tau)$ cannot be the same one save
for finitely many $h$'s. This is how we prove that infinitely many of them
are positive (resp. negative).
\subsection{Proof of Theorem~\ref{lowerbound}}
To prove Theorem~\ref{lowerbound}, we first introduce a smooth
non-negative function $W$ verifying $(W_1)$, $(W_2)$ and $(W_3)$ stated in the
introduction and write
\begin{equation*}
  \sum_{\lowm<q/Q\le \upm}\sum_{a\mode q}
  |S(\varphi, a/q)|^2\ge
  \sum_{q\ge 1}\frac{W(q/Q)}{q}\sum_{a\mode q}
  |S(\varphi,a/q)|^2.
\end{equation*}
We then use Theorem~\ref{Precise}.
Theorem~\ref{yoddle} is our next step, with $\xi=1$.
We select $H=L=\exp c\sqrt{\log N} \tau$  for some small but positive $c$. Given
$h\le H$, we may first employ the first statement of Theorem~\ref{bounds} together
with~\eqref{neat} and~\eqref{defIW} to get that 
$\frac{\tau}{h}\lambda_{h,\ell}\le \ooWo+\Ocal(1/\sqrt{Q})$.
This already ensures us that
\begin{multline*}
  N\sum_{h\le H}\frac{\tau}{h}
  \sum_{\substack{\ell\le L,\\ |\lambda_{h,\ell}|\ge
      \xi\exp-{c_3}\sqrt{\log N}}}\lambda_{h,\ell}
  \sum_{a\mode q} |S(\varphi, a/q)|^2
  \\
  \le 
  \ooWo N\|\varphi\|_2^2\Bigl(1+ H^2L\exp-c_3\sqrt{\log N}+Q^{-1/2}\Bigr).
\end{multline*}
This is not quite enough.
The full strength of Theorem~\ref{yoddle} uses the non-negativity of
$W$. We employ this theorem with $U_2=\tau/h$, and this gives us that 
\begin{equation*}
  \frac{\tau}{h}\lambda_{h,\ell}\le \ooWo(1-ce^{-c_4\tau/h})
  +\Ocal(1/\sqrt{Q}).
\end{equation*}
Theorem~\ref{lowerbound} readily follows.




\begin{thebibliography}{10}

\bibitem{Amrein-Berthier*77}
W.~O. Amrein and A.~M. Berthier.
\newblock On support properties of {$L^{p}$}-functions and their {F}ourier
  transforms.
\newblock {\em J. Functional Analysis}, 24(3):258--267, 1977.

\bibitem{Boas*54}
Ralph~Philip Boas, Jr.
\newblock {\em Entire functions}.
\newblock Academic Press Inc., New York, 1954.

\bibitem{Boca-Radziwill*16}
Florin~P. Boca and Maksym Radziwi\l\l.
\newblock Limiting distribution of eigenvalues in the large sieve matrix.
\newblock {\em J. Eur. Math. Soc. (JEMS)}, 22(7):2287--2329, 2020.

\bibitem{Bombieri*71}
E.~Bombieri.
\newblock A note on the large sieve.
\newblock {\em Acta Arith.}, 18:401--404, 1971.

\bibitem{Bottcher*13}
Albrecht {B}{\"o}ttcher.
\newblock {Best constants for Markov type inequalities in Hilbert space norms.}
\newblock In {\em {Recent trends in analysis. Proceedings of the conference in
  honor of Nikolai Nikolski on the occasion of his 70th birthday, Bordeaux,
  France, August 31 -- September 2, 2011}}, pages 73--83. Bucharest: The Theta
  Foundation, 2013.

\bibitem{Chan-Kumchev*12}
T.~H. Chan and A.~V. Kumchev.
\newblock On sums of {R}amanujan sums.
\newblock {\em Acta Arith.}, 152(1):1--10, 2012.

\bibitem{Cochran*72}
James~Alan Cochran.
\newblock {\em The analysis of linear integral equations}.
\newblock McGraw-Hill Book Co., New York-D\"usseldorf-Johannesburg, 1972.
\newblock McGraw-Hill Series in Modern Applied Mathematics.

\bibitem{Conrey-Iwaniec*03}
B.~Conrey and H.~Iwaniec.
\newblock Asymptotic large sieve.
\newblock {\em Private communication}, 2003.

\bibitem{Conrey-Iwaniec-Soundararajan*12}
J.~B. Conrey, H.~Iwaniec, and K.~Soundararajan.
\newblock Small gaps between zeros of twisted {L}-functions.
\newblock {\em Acta Arith.}, 155(4):353--371, 2012.

\bibitem{Dress*99}
F.~Dress.
\newblock Discr{\'e}pance des suites de {F}arey.
\newblock {\em J. Th\'eor. Nombres Bordx.}, 11(2):345--367, 1999.

\bibitem{Duke-Friedlander-Iwaniec*93}
W.~Duke, {J.B.} Friedlander, and H.~Iwaniec.
\newblock Bounds for automorphic $l$-functions.
\newblock {\em Invent. Math.}, 112(1):1--8, 1993.

\bibitem{Duke-Iwaniec*92}
{W.} Duke and {H.} Iwaniec.
\newblock Estimates for coefficients of {$L$}-functions. ii.
\newblock In Bombieri E., editor, {\em Proceedings of the Amalfi conference on
  analytic number theory}, pages 71--82, held at Maiori, Amalfi, Italy, from 25
  to 29 September, 1989. Salerno, 1992.

\bibitem{Elliott*85b}
P.D.T.A. {Elliott}.
\newblock Additive arithmetic functions on arithmetic progressions.
\newblock {\em Proc. London Math. Soc.}, 54(3):15--37, 1985.

\bibitem{Elliott*91}
P.D.T.A. {Elliott}.
\newblock On maximal variants of the {L}arge {S}ieve.
\newblock {\em J. Fac. Sci. Univ. Tokyo, Sect. IA}, 38:149--164, 1991.

\bibitem{Elliott*97}
P.D.T.A. {Elliott}.
\newblock {\em Duality in analytic number theory}, volume 122 of {\em Cambridge
  Tracts in Mathematics}.
\newblock Cambridge University Press, Cambridge, 1997.

\bibitem{Erdelyi-Magnus-Oberhettinger-Tricomi*54}
A.~Erd{\'e}lyi, W.~Magnus, F.~Oberhettinger, and F.~G. Tricomi.
\newblock {\em Tables of integral transforms. {V}ol. {I}}.
\newblock McGraw-Hill Book Company, Inc., New York-Toronto-London, 1954.
\newblock Based, in part, on notes left by Harry Bateman.

\bibitem{Erdos-Renyi*68}
P.~Erd{\H{o}}s and A.~R{\'e}nyi.
\newblock Some remarks on the large sieve of {Y}u. {V}. {L}innik.
\newblock {\em Ann. Univ. Sci. Budapest. E\"otv\"os Sect. Math.}, 11:3--13,
  1968.

\bibitem{Fredholm*03}
Ivar Fredholm.
\newblock Sur une classe d'\'equations fonctionnelles.
\newblock {\em Acta Math.}, 27(1):365--390, 1903.

\bibitem{Friedlander-Iwaniec*92}
J.~Friedlander and H.~Iwaniec.
\newblock A mean-value theorem for character sums.
\newblock {\em Mich. Math. J.}, 39(1):153--159, 1992.

\bibitem{Gallagher*67}
{P.X.} Gallagher.
\newblock The large sieve.
\newblock {\em Mathematika}, 14:14--20, 1967.

\bibitem{Gohberg-Krein*69}
I.~C. Gohberg and M.~G. Kre{\u i}n.
\newblock {\em Introduction to the theory of linear nonselfadjoint operators}.
\newblock Translated from the Russian by A. Feinstein. Translations of
  Mathematical Monographs, Vol. 18. American Mathematical Society, Providence,
  R.I., 1969.

\bibitem{Gohberg-Goldberg-Krupnik*00}
Israel Gohberg, Seymour Goldberg, and Nahum Krupnik.
\newblock {\em Traces and determinants of linear operators}, volume 116 of {\em
  Operator Theory: Advances and Applications}.
\newblock Birkh\"auser Verlag, Basel, 2000.

\bibitem{Harcos*03}
Gergely Harcos.
\newblock An additive problem in the {F}ourier coefficients of cusp forms.
\newblock {\em Math. Ann.}, 326(2):347--365, 2003.

\bibitem{Havin-Joricke*94}
Victor Havin and Burglind J{\"o}ricke.
\newblock {\em The uncertainty principle in harmonic analysis}, volume~28 of
  {\em Ergebnisse der Mathematik und ihrer Grenzgebiete (3) [Results in
  Mathematics and Related Areas (3)]}.
\newblock Springer-Verlag, Berlin, 1994.

\bibitem{Iwaniec-Kowalski*04}
H.~Iwaniec and E.~Kowalski.
\newblock {\em Analytic number theory}.
\newblock American Mathematical Society Colloquium Publications. American
  Mathematical Society, Providence, RI, 2004.
\newblock xii+615 pp.

\bibitem{Jutila*92}
M.~Jutila.
\newblock Transformations of exponential sums.
\newblock In {\em Proceedings of the {A}malfi {C}onference on {A}nalytic
  {N}umber {T}heory ({M}aiori, 1989)}, pages 263--270. Univ. Salerno, Salerno,
  1992.

\bibitem{Jutila*07}
Matti Jutila.
\newblock Distribution of rational numbers in short intervals.
\newblock {\em Ramanujan J.}, 14(2):321--327, 2007.

\bibitem{Kobayashi*73}
{I.} Kobayashi.
\newblock A note on the {S}elberg sieve and the large sieve.
\newblock {\em Proc. Japan Acad.}, 49(1):1--5, 1973.

\bibitem{Kovrijkine*00}
Oleg Kovrijkine.
\newblock Some results related to the {L}ogvinenko-{S}ereda theorem.
\newblock {\em Proc. Amer. Math. Soc.}, 129(10):3037--3047 (electronic), 2001.

\bibitem{Logvinenko-Sereda*74}
V.~N. Logvinenko and Ju.~F. Sereda.
\newblock Equivalent norms in spaces of entire functions of exponential type.
\newblock {\em Teor. Funkci\u \i\ Funkcional. Anal. i Prilo\v zen.}, (Vyp.
  20):102--111, 175, 1974.

\bibitem{Nazarov*92}
F.~L. Nazarov.
\newblock On the theorems of {T}ur\'an, {A}mrein and {B}erthier, and {Z}ygmund.
\newblock {\em Zap. Nauchn. Sem. S.-Peterburg. Otdel. Mat. Inst. Steklov.
  (POMI)}, 201(Issled. po Linein. Oper. Teor. Funktsii. 20):117--123, 191,
  1992.

\bibitem{Nazarov*93}
F.~L. Nazarov.
\newblock Local estimates for exponential polynomials and their applications to
  inequalities of the uncertainty principle type.
\newblock {\em Algebra i Analiz}, 5(4):3--66, 1993.

\bibitem{Niederreiter*73}
H.~Niederreiter.
\newblock The distribution of {F}arey points.
\newblock {\em Math. Ann.}, 201:341--345, 1973.

\bibitem{Ramare*07a}
O.~{Ramar\'e}.
\newblock Eigenvalues in the large sieve inequality.
\newblock {\em Funct. Approximatio, Comment. Math.}, 37:7--35, 2007.

\bibitem{Ramare*06}
O.~Ramar{\'e}.
\newblock {\em Arithmetical aspects of the large sieve inequality}, volume~1 of
  {\em Harish-Chandra Research Institute Lecture Notes}.
\newblock Hindustan Book Agency, New Delhi, 2009.
\newblock With the collaboration of D. S. Ramana.

\bibitem{Ramare*09a}
O.~{Ramar\'e}.
\newblock Eigenvalues in the large sieve inequality, {II}.
\newblock {\em J. Th\'eorie N. Bordeaux}, 22(1):181--196, 2010.

\bibitem{Ramare*22}
O.~Ramar\'e.
\newblock Notes on restriction theory in the primes.
\newblock {\em Israel J. of Math.}, page 21pp, 2022.

\bibitem{Renyi*70}
A.~R{\'e}nyi.
\newblock {\em Probability theory}.
\newblock North-Holland Publishing Co., Amsterdam, 1970.
\newblock Translated by L{\'a}szl{\'o} Vekerdi, North-Holland Series in Applied
  Mathematics and Mechanics, Vol. 10.

\bibitem{Sakhnovich*15}
Lev~A. Sakhnovich.
\newblock {\em Integral equations with difference kernels on finite intervals},
  volume~84 of {\em Operator Theory: Advances and Applications}.
\newblock Birkh\"auser/Springer, Cham, 2015.
\newblock Second edition, revised and extended.

\bibitem{SchlagePuchta*09}
Jan-Christoph Schlage-Puchta.
\newblock Lower bounds for expressions of large sieve type.
\newblock {\em Arch. Math. (Brno)}, 45(2):79--82, 2009.

\bibitem{Shampine*66}
Lawrence~F. Shampine.
\newblock An inequality of {E}. {S}chmidt.
\newblock {\em Duke Math. J.}, 33:145--150, 1966.

\bibitem{Swann*71}
Dale~W. Swann.
\newblock Some new classes of kernels whose {F}redholm determinants have order
  less than one.
\newblock {\em Trans. Amer. Math. Soc.}, 160:427--435, 1971.

\bibitem{sagemath}
{The Sage Developers}.
\newblock {\em {S}ageMath, the {S}age {M}athematics {S}oftware {S}ystem
  ({V}ersion 9.5)}, 2022.
\newblock {\tt https://www.sagemath.org}.

\bibitem{Titchmarsh*51}
E.C. Titchmarsh.
\newblock {\em The {T}heory of {R}iemann {Z}eta {F}unction}.
\newblock Oxford Univ. Press, Oxford 1951, 1951.

\bibitem{Uchiyama*72}
{Sabur\^o} Uchiyama.
\newblock The maximal large sieve.
\newblock {\em Hokkaido Math. J.}, 1:117--126, 1972.

\bibitem{Widom*66}
Harold Widom.
\newblock Hankel matrices.
\newblock {\em Trans. Amer. Math. Soc.}, 121:1--35, 1966.

\bibitem{Wolke*74}
D.~Wolke.
\newblock A lower bound for the large sieve inequality.
\newblock {\em Bull. London Math. Soc.}, 6:315--318, 1974.

\end{thebibliography}

\end{document}